\documentclass[10pt]{article}
\usepackage[utf8]{inputenc}
\usepackage[english]{babel}
\usepackage{amsmath,amsthm,amsfonts}
\usepackage{graphicx}
\usepackage{bm}
\usepackage{tikz}
\usepackage{multirow}
\usepackage{subcaption}

\newcommand{\grad}{\mathop{\rm grad}\nolimits}
\renewcommand{\div}{\mathop{\rm div}\nolimits}

\newcommand{\norm}[1]{\left\|#1\right\|}

\graphicspath{ {./figs/} }

\title{Preconditioning Markov Chain Monte Carlo Method for Geomechanical Subsidence using multiscale method and machine learning technique}

\author{
Maria Vasilyeva \thanks{Institute for Scientific Computation (ISC),
Texas A\&M University,
College Station, TX 77843-3368, USA \& 
Institute of Mathematics and Informatics, 
North-Eastern Federal University, 
Yakutsk, Republic of Sakha (Yakutia), Russia, 677980 
}
\and
Aleksei Tyrylgin \thanks{
Institute of Mathematics and Informatics, 
North-Eastern Federal University, 
Yakutsk, Republic of Sakha (Yakutia), Russia, 677980 
}
\and
Donald L. Brown \thanks{
School of Mathematical Sciences, 
GeoEnergy Research Center,
The University of Nottingham,
University Park
Nottingham, NG7 2RD, United Kingdom 
}
\and
Anirban Mondal
 \thanks{Department of Mathematics, Applied Mathematics, and Statistics
Case Western Reserve University
Cleveland, OH 44106}
}

\begin{document}
\maketitle

\begin{abstract}
In this paper, we consider the numerical solution of the poroelasticity problem with stochastic properties. 
We present a Two-stage Markov Chain Monte Carlo method for geomechanical subsidence. 
In this work, we study two techniques of preconditioning: (MS) multiscale method for model order reduction and (ML) machine learning technique. 
The purpose of preconditioning is the fast sampling, where a new proposal is first testes by a cheap multiscale solver or using fast prediction of the neural network and the full fine grid computations will be conducted only if the proposal passes the first step. 
To construct a reduced order model, we use the Generalized Multiscale Finite Element Method and present construction of the multiscale basis functions for pressure and displacements in stochastic fields.  
In order to construct a machine learning based preconditioning, we generate a dataset using a multiscale solver and use it to train neural networks.  
The Karhunen-Lo{\'e}ve expansion is used to represent the realization of the stochastic field. 
Numerical results are presented for two- and three-dimensional model examples. 
\end{abstract}

%\begin{keywords}
%poroelasticity problems, Two-stage Markov Chain Monte Carlo method,  fast sampling, multiscale method, machine learning
%\end{keywords}

%\begin{AMS}\end{AMS}

\section{Introduction}

Modelling geomechanical problems has important implications on the understanding of many of the physical processes. 
For example, in areas such as environmental engineering with  modelling permafrost compaction and subsidence \cite{mitchell1977mechanical}, and the oil and gas industry for reservoir geomechanics to increase production and overall reservoir life \cite{sayers2007introduction,zoback2010reservoir}. 
There are many challenges in simulating these problems. 
One challenging being the length scales involved. There is a  high variability in material parameters such as permeability and porosity for flow and Young's modulus in the case of mechanical properties. 
Further compounding the problem is lack of knowledge of these parameters due to subsurface depth. This uncertainty may arise from under resolution from seismic data or, in the case of resistivity measurements with electromagnetic sensors, unclear physical  matching of properties to data. 

% pororelasticity problem
%Problems of mechanics and flow in porous media  have wide ranging applications in  many areas of science and engineering. Particularly in geomechanical modeling and its applications to reservoir engineering
%for enhanced production and environmental safety due to overburden subsidence and compaction \cite{sayers2007introduction,zoback2010reservoir}. One of the key challenges is the multiscale nature of the geomechanical problems.
%
As noted, the computational challenge is often two fold. First, heterogeneity of subsurface properties need to  be accurately accounted for in the geomechanical model. In turn,  this requires  high resolution  with adding many degrees of freedom that can be computationally expensive.  
Second is the uncertainty in the subsurface properties such as permeability and elastic parameters.  However, in today's ``Big-Data'' world, often a plethora of information is available to help characterize the subsurface. For example, in reservoir engineering there is extensive logs on production data such as oil-cut of produced wells. %\misref{refFractionalFlowhistorymatching}. 
In addition, other metrics such as bottom hole pressure and time of flight have been extensively used in the history matching literature to integrate data into the subsurface models \cite{ma2008efficient, yin2011hierarchical}. In the areas of environmental science and engineering,  near-surface resistivity  measurements are often available to obtain constraints on subsurface moisture content reflit. 
%Moroever, there are disparate scales between the often relatively thin reservoir structure and the large overburden  surrounding the reservoir that adds more complexity to the simulation.   
% Therefore, we propose a multiscale method  to attempt overcome some of these challenges.

%Various other UQ methodologies, be politically correct here! 
There are various  techniques to integrate data into subsurface models.
The literature is quite extensive and various uncertainty quantification methodologies exist and have various advantages and disadvantages.  %\misref{reflaundrylistofUQ}.
Various data assimilations and variants are often used to integrate dynamic data information constantly into models. This is particularly useful for weather forecast modelling as well as geophysical problems. 
In this work, we will consider a Bayesian framework of updating a prior probability  distribution with information to obtain a data-integrated posterior. 
It is well known that the resulting problems are high-dimensional and suffer from the curse of dimensionality. To circumvent some of this difficulty  Monte-Carlo based techniques have been proposed.
For example, a popular and useful technique is the Ensemble Kalman Filter (EnKF) \cite{evensen2003ensemble} to integrate data into subsurface models. In this work, however, we consider a Markov-Chain Monte-Carlo method (MCMC) based on the classical Metropolis-Hastings sampling algorithm \cite{chib1995understanding}. This is an effective tool to  efficiently sample from an unknown posterior distribution that is conditioned to the data. 
%
%problems with MCMC, use coarse-scale models.
However, the MCMC algorithms suffer from the fact that many simulations must be computed \cite{he2019data, kumar2020parametric}.  This is again tied to the multiscale and high-contrast nature of the material properties as direct numerical simulations must fully resolve these scales. Moreover, the acceptance rate of proposed subsurface properties in the Metropolis-Hastings algorithm is known to be very low.

%History of Geomechanics
With the availability surface sensors and satellite observation data (InSAR) \cite{massonnet1998radar}, it is now possible to better characterize the geomechanical picture of the subsurface based on this data. 
This has proven particularly useful in the context of monitoring large scale  CO2 sequestration projects \cite{vasco2010satellite, rucci2013monitoring}. There are various subsurface mechanical models that one can consider, however, in this work we will work with a poroelastic earth model.
The  mathematical structure  of the poroelasticity models  are  coupled equations for pressure and displacements known as Biot models \cite{biot1941general}. 
%For pressure, or flow equations, we have the  parabolic equation Darcy equation with a time dependent coupling to volumetric strain. 
%The stress equation is the quasi-static elasticity equations with a coupling to the pressure gradients as a forcing. 
%
Poroelastic models of this type have been explored in the petroleum engineering  literature in the context of geomechanics for some time \cite{settari2001advances, settari1998coupled, mikelic2013convergence, minkoff2003coupled, rahrah2020moving}. 
An interesting surface observation application such as permafrost modelling would require thermal  and partially saturated moisture considerations, but we leave this to future work. 
% There are noted issues that arise. The first being heterogeneities of the reservoir and surrounding media add many complications to 
%the effective simulation due to complexity of scales. Moreover, development of flow and mechanics simulation were often considered separately. 
%Progress was made on this problem by considering various  coupling strategies \cite{settari1998coupled}. However, in the instance that the physics is not well understood a 
%fully coupled scheme may be desired. 

As mentioned prior, the MCMC algorithm suffer from large number of simulation runs and low acceptance rates. 
A useful technique to expedite this procedure is preconditioning the Metropolis-Hasting algorithm in a two-stage (or multi-stage) procedure utilizing coarse-scale or upscaled models \cite{efendiev2005efficient}. This is accomplished by using the coarser-scale simulation as a prior filtering stage in the accept-rejection procedure. 
There are many  effective multiscale frameworks that have been developed in recent years for the poroelasticity problem 
\cite{castelletto2019multiscale, sokolova2019multiscale, altmann2020computational, tyrylgin2020generalized, fu2019computational, vasilyeva2019constrained}. 
We will use the Generalized Multiscale Finite Element Method (GMsFEM) framework developed for poroelasticty  in \cite{brown2016generalized1, brown2016generalized2}, which is a generalization of the multiscale finite element method to build our course-scale models \cite{efendiev2009multiscale}. 
Utilizing GMsFEMs for this application has also been useful in Multi-Level MCMC for elliptic problems with high-contrast \cite{ginting2011application}.

% %
% % Why use MS methods for MCMC past literature and other ms methods
% %
% %Once the equations have been split in time we wish to resolve in space and will utilize a multiscale method. 
% There are many  effective multiscale frameworks that have been developed in recent years. 
% %There are rigorous approaches based on homogenization of partial differential equations, 
% %where effective equations are derived based fine-scale equations at the microstructure level \cite{brown2011,brown2014}. 
% %However, these approaches may have limited computational use and more practical multiscale methods are used. 
% Examples include the Heterogeneous Multiscale Method (HMM), where  macro-scale equations on coarse-grids are solved while local problems 
% on the fine-scale are resolved at each coarse grid nodes to build effective coefficients
%  \cite{E:Engquist:2003,Abdulle:E:Engquist:Vanden-Eijnden:2012}. 
%  An approach based on the Variational Multiscale Method (see \cite{MR2300286}), utilizing coarse-grid quasi-interpolation 
% operators  to build an orthogonal splitting into a multiscale space and a fine-scale space \cite{MP11}.
% Fine-scale space corrections are then localized to create a computationally tractable. 
 
%HMM
%VMS
%MsFEM

%  our method
The GMsFEM has the advantage of being able to capture small scale features from the heterogeneities into coarse-grid basis functions and offline spaces, as well as having a unified computational grids for both mechanics and flow solves. 
The offline multiscale basis construction may proceed in both fluid and mechanics in parallel and both constructions are comparable. 
First  a coarse-grid is generated  and in each grid block a local static problem with varying boundary conditions is solved to construct the snapshot spaces. 
We then perform a dimension reduction of the snapshot space by solving auxiliary eigenvalue problems. Taking the corresponding smallest eigenpairs,  and multiplying by a multiscale partition of unity we are able to construct our offline basis. 
In this greatly reduced dimension offline basis, the online solutions may be calculated for pressure and displacements for any viable boundary condition or forcing.  
Given a set or material properties such as permeability a set of standard MsFEM basis functions must be computed, however, utilizing the GMsFEM's ability to handle parameters (as is used in nonlinear GMsFEM) as well as scales we are able to compute a single set of enriched basis functions for many possible realizations of physical parameters. 
% ------------------------------------------

% ml
For further reduction of the computational time of the first stage in the two-stage MCMC method, we present a machine learning technique \cite{vasilyeva2020learning, vasilyeva2018machine}. 
The machine learning is used to quickly predict displacements for the estimation of the proposed fields. 
We generate a dataset using a multiscale solver and use it to train neural networks and learn dependencies between heterogeneous properties and displacements in each direction on the top boundary. 
As soon as neural networks are trained on the dataset, fast calculations can be performed as preconditioning of the MCMC method.  
We use a convolutional neural network and GPU training process to construct a machine learning algorithm \cite{lecun2015deep, krizhevsky2012imagenet}.

% in this paper
The work is organized as follows.
In Section 2, we provide the mathematical background of the geomechanical subsidence problem. 
We introduce the Biot type model and highlight where the heterogeneity primarily occurs. 
In Section 3, we discuss the fine-scale and GMsFEM (coarse-scale) solution of the geomechanical problem. We utilize the nonlinear GMsFEM for poroelasticity developed in \cite{brown2016generalized1}, whereby we are
able to handle parametrization. 
In Section 4 we introduce the broad concepts of Bayesian uncertainty quantification. We discuss the parametrization of the randomness via the Karhunen-Loeve expansion so that when we search in the MCMC procedure the dimension of the space is reduced. We outline the single and two-stage Metropolis-Hastings algorithms used in the accept-reject procedure to sample from the ``data-integrated'' posterior distribution. 
In Section 5, we present a numerical algorithm based on two- and three-dimensional synthetic data at the surface to show the efficiency of the method as an expedited MCMC sampling method. 

\section{Problem Formulation} 

We start with the description of the general geomechanical model that we use in our simulations. We keep the discussion very general and abstract, but ultimately the idea is to have one surface boundary that has open surface boundary conditions, and a truncated in situ ground that has fixed motion and is in physical reality connected to a much larger (functionally infinite) domain. 

We denote computational domain $\Omega\subset\mathbb{R}^d, d=2,3,$ to be a bounded sufficiently smooth (Lipschitz) region.  
We consider linear poroelasticity problem, with random or uncertain coefficients. That we may view as parameters. 
We wish to find a pressure $p$ and displacements $u$ satisfying the following Biot effective stress poroelasticity law
\begin{equation}
\label{eq:main}
\begin{split}
- \div  \sigma(x, \theta,  u) + \alpha \grad p  &= 0, 
\quad  x \in \Omega, \quad t > 0, \\
 \alpha \frac{\partial \div  u}{\partial t} + \frac{1}{M} \frac{\partial p}{\partial t} + \div  q(x, \theta, p) &= 0, 
\quad  x \in \Omega, \quad t > 0,
\end{split}
\end{equation}
where $M$ is the Biot modulus and $\alpha$ is the Biot-Willis fluid-solid coupling coefficient. These terms are lower order derivatives, so for simplicity we will suppose that these are constants and not random.  Body forces, such as gravity, are neglected without loss of generality.

Here we suppose that the stress tensor $\sigma$ and flux $q$ depend on both space $x$ and a large dimensional random parameter $\theta$
\[
q(x, \theta, p) = -\frac{k(x,\theta)}{\nu} \grad p,
\quad
\sigma(x, \theta, u) = 2 \mu (x,\theta) \varepsilon (u) 
 + \lambda(x,\theta) \div u  \,  \mathcal{I},
\]
where $\varepsilon (u) =  ( \grad  u + \grad  u^T )/2$, 
$\nu$ is the fluid viscosity, 
$k(x,\theta)$ is the permeability,  
$\mu(x,\theta)$, $\lambda(x,\theta)$ are Lam{\'e} coefficients, 
$\mathcal{I}$ is the identity tensor.  

For Lam{\'e} coefficients $\lambda$ and $\mu$, we have following relations
\[
\mu(x,\theta) = \frac{E(x,\theta)}{2 (1 + \eta)}, \quad
\lambda(x,\theta) = \frac{E(x,\theta) \eta}{(1+ \eta) ( 1- 2 \eta)},
\]
where $E(x,\theta)$ is the random spatially varying elastic modulus and $\eta$ is a constant Poisson's ratio. 
One could choose to vary both, but for this work we consider a varying elastic modulus. 
Thus, the coefficients $k(x,\theta)$, $\mu(x,\theta)$ and $\lambda(x,\theta)$ may be highly variable and contain randomness or uncertainty.  

We denote the initial condition for pressure 
\[
p = p_0, \quad x \in \Omega, \quad t = 0.
\]
In general, we suppose the following Neumann and Robyn boundary conditions on each portion
\[
 u = 0, \quad  x \in \Gamma_u, \quad 
 \sigma \cdot n = 0, \quad  x \in \partial \Omega / \Gamma_u, 
 \]
 and 
 \[
q \cdot n = \gamma (p - p_1), \quad  x \in \Gamma_p, \quad  
q \cdot n= 0, \quad  x \in \partial \Omega / \Gamma_p,
\]
where $n$ is the unit normal to the boundary. 

Here the primary sources of the heterogeneity and uncertainty in the physical properties are from mechanical properties related to $ E(x,\theta)$ and fluid flow properties related to $k(x,\theta)$.  
%The pressure and displacements depend on space $x$, time $t$, and randomness $\theta$, but we suppress this notation unless there an ambiguity. 

%\section{Fine-Scale Discretizations}
%Here describe the numerical solution approaches for the random poroelastic system \eqref{eq:main}. We start with the fine grid approximation.

To  solve \eqref{eq:main}, we use a standard  finite element method and implicit time integration. We have following variational formulation: find $(u,p) \in V \times Q$ such that 
\begin{equation}
\label{eq:fine}
\begin{split}
a(\theta; u,  v) + g(p,  v) &= 0, \quad \forall v \in V,  \\
d \left( \frac{ u -  \check{u}}{\tau}, q \right)  
+ m \left( \frac{p - \check{p}}{\tau}, q \right)
+ b(\theta; p, q) &= l(q), \quad \forall q \in Q,
\end{split}
\end{equation}
where 
$V \in \{v \in [H_1(\Omega)]^d: v(x) = 0, x \in \Gamma_u \}$,   
$Q = H_1(\Omega)$, $\tau$ is the time step, $\check{u}$ and $\check{p}$ are the solutions from previous time layer. 

Here for bilinear and linear forms we have 
\[
a(\theta; u,  v)  = \int_{\Omega}  \sigma(x, \theta, u) : \varepsilon(v) dx, 
\quad
g(p,  v)  = \int_{\Omega}\alpha \, \grad p \, v \, dx,
\]\[
b(\theta; p, q) = \int_{\Omega} \frac{k(x, \theta)}{\nu} \grad p \cdot \grad q \, dx + \int_{\Gamma_p} \gamma \, p \, q\, ds,
\quad
l(q) = \int_{\Gamma_p} \gamma \, p_1 \, q\, ds,
\]\[
m(p, q) = \int_{\Omega} \frac{1}{M} \, p \, q \, dx,  
\quad
d( u, q) = \int_{\Omega} \alpha \, \div  u \, q \, dx.
\]

Let $\mathcal{T}^h$ be a fine grid partition of the computational domain $\Omega$ into finite elements and
\[
u = \sum u_i \phi_i, \quad 
p = \sum p_i \psi_i,
\]
where $\phi_i$ and $\psi_i$ are the linear  basis functions defined on $\mathcal{T}^h$. 

Therefore, we have following matrix form on the fine grid
\begin{equation}
\label{eq:mat}
\begin{split}
A^h(\theta) u + G^h p & = 0, \\
D^h \frac{u - \check{u}}{\tau} 
+ M^h \frac{p - \check{p}}{\tau} 
+ B^h(\theta) p & = F^h,
\end{split}
\end{equation}
where  
$M^h = [m_{ij}], \quad m_{ij} =  m(\psi_i, \psi_j)$,  
$B^h = [b_{ij}], \quad b_{ij} =  b(\theta; \psi_i, \psi_j)$,  
$A^h = [a_{ij}], \quad a_{ij} =  a(\theta; \phi_i, \phi_j)$,  
$D^h = [d_{ij}], \quad d_{ij} =  d(\phi_i, \psi_j)$, 
$G^h = [g_{ij}], \quad g_{ij} =  g(\phi_i, \psi_j)$     
and $F = \{ f_j \}$, $f_j  = l(\psi_j)$,

\section{Coarse-Scale Discretization}

To construct a reduced order model on the coarse grid, we use a Generalized Multiscale Finite Element Method (GMsFEM). 
We construct an offline multiscale space for pressure and displacements using some number of random coefficients. Therefore constructed basis functions can be used for any input parameters $k(x, \theta)$ and $E(x, \theta)$. %[Efendiev-Iliev-Kronsbein-MLMC, Efendiev-Hou-Luo].

\begin{figure}[h!]
\centering
    \begin{subfigure}[b]{0.45\textwidth}
        \includegraphics[width=1.0\linewidth]{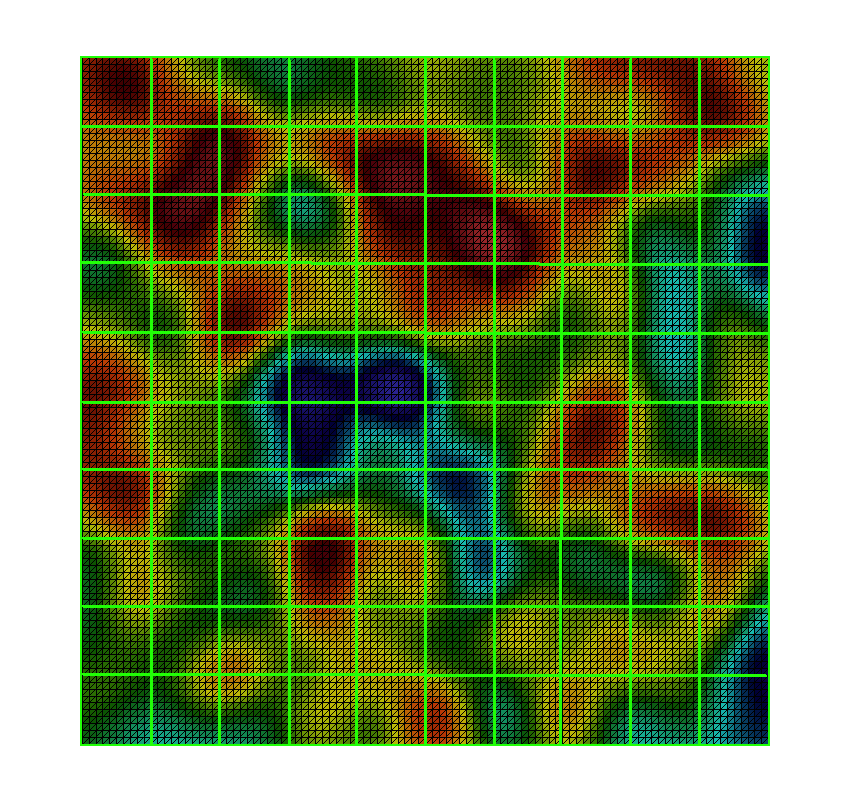}
        \caption{2D}
    \end{subfigure}
    \begin{subfigure}[b]{0.49\textwidth}
        \includegraphics[width=1.0\linewidth]{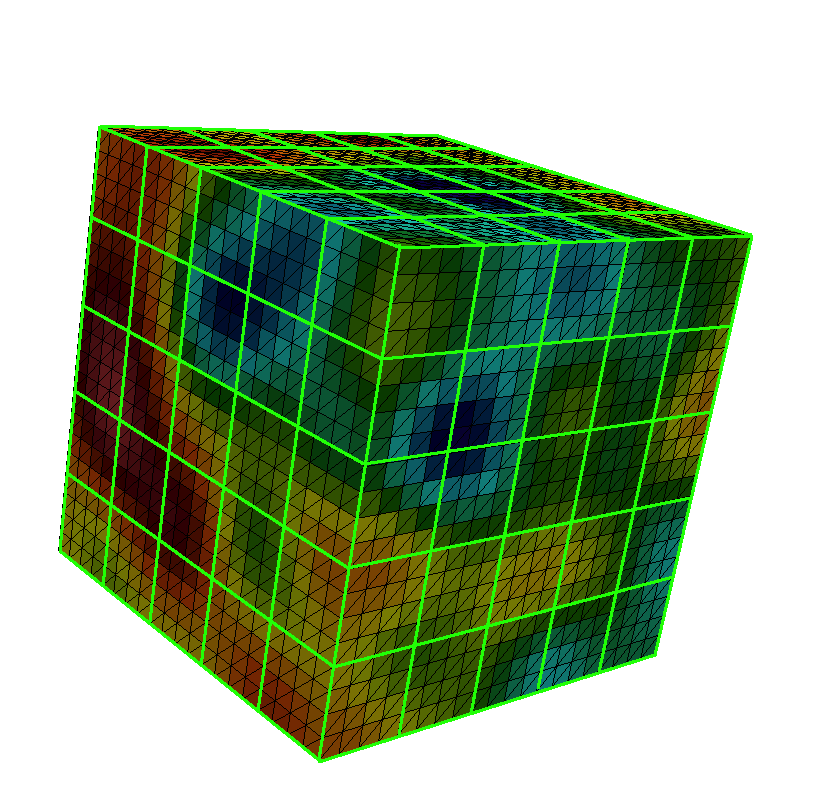}
        \caption{3D}
    \end{subfigure}
\caption{Illustration of a coarse grid and fine grid (fine grid - black color, coarse grid - green color.  
(a) 2D domain with $10 \times 10$ coarse mesh. 
(b) 3D domain with $5 \times 5 \times 5$ coarse mesh}
\label{schematic}
\end{figure}

Let $\mathcal{T}^H$ be a standard conforming partition of the computational domain $\Omega$ into finite elements (Figure~\ref{schematic}). 
We refer to this partition as the coarse-grid and assume that each coarse element is partitioned into a connected union of fine grid blocks.  
We use $\{x_i\}_{i=1}^{N_c}$ to denote the vertices of the coarse mesh $\mathcal{T}^H$, and define the neighborhood of the node $x_i$ by
\[
\omega_i=\bigcup_{j}\left\{ K_j \in\mathcal{T}^H \, | \,  x_i\in \overline{K}_j\right\},
\]
where $K_j$ is the coarse cell and $N_c$ is the number of coarse nodes. 

The main idea for solution problem in the stochastic media is to precompute multiscale basis functions using a number of permeability and elastic modulus realizations and use them for the making inexpensive coarse grid calculations. 

\subsection{Multiscale basis functions for pressure}
\label{pressuresolve}

% spashot space 
In the offline computation, we first construct a snapshot space. 
Construction involves solving the local problem for various choices of input parameters  and various boundary conditions. 
For each fixed parameter $\theta_r$ ($r = 1,...,N_r$), we propose a  snapshot space generated by harmonic extensions of $b(\theta_r; p,q)$
\begin{equation} 
\label{harmonic_ex}
\begin{split}
b(\theta_r ; \psi_{r,j}^{\omega_i, \text{snap}} , q) &= 0 \quad 
x \in \omega_i, \\
\psi_{r,j}^{\omega_i, \text{snap}} &=g_j(x) \quad 
x \in \partial \omega_i,
\end{split}
\end{equation}
where $g_j(x)=\delta_{j,k}$, $\forall k \in J_{h}(\omega_i)$ ($J_{h}(\omega_i)$ denotes the fine-grid boundary node on $\partial \omega_i$.  

We collect all solutions as a snapshot space in the local domain $\omega_i$ to define local snapshot space
\[
Q^{\omega_i, \text{snap}} = 
\text{span}\{ 
\psi_{r,j}^{\omega_i, \text{snap}}: \, 
1 \leq r \leq N_r, \, 
0 \leq j \leq N^{\partial \omega_i}
\},  
\]
where $N^{\partial \omega_i}$ is the number of the fine-grid boundary nodes on $\partial \omega_i$.

We reorder the snapshot functions using a single index to create the matrix
\[
R_p^{\omega_i, \text{snap}} = \left[ 
\psi_{1}^{\omega_i, \text{snap}}, \ldots, 
\psi_{M_{p,\text{snap}} }^{\omega_i, \text{snap}} 
\right]^T,
\]
where $M_{p,\text{snap}} = N^{\partial \omega_i} N_r$ denotes the total number of functions to keep in the snapshot construction. 

% offline space
To construct the offline space $Q_{\omega_i, \text{off}}$, we perform a dimension reduction of the space of snapshots by using an auxiliary spectral decomposition. More precisely,  we solve the  eigenvalue problem in the space of snapshots:
\begin{equation} 
\label{offeig}
B^{\omega_i, \text{off}} \Psi_j^{\omega_i, \text{off}} 
= \lambda_j  S^{\omega_i, \text{off}} \Psi_j^{\omega_i, \text{off}},
\end{equation}
where
\[
\begin{split}
B^{\omega_i, \text{off}} = \{b_{lk}\}, \quad 
b_{lk} &= \int_{\omega_i}  
\overline{k}(x) \,
\nabla \psi^{\omega_i, \text{snap}}_l , 
\nabla \psi^{\omega_i, \text{snap}}_k  \, dx 
= R_p^{\omega_i, \text{snap}} \, \overline{B} \, (R_p^{\omega_i, \text{snap}})^T,
\\
S^{\omega_i, \text{off}} = \{s_{lk}\}, \quad 
s_{lk} &= \int_{\omega_i}  
\overline{k}(x)  \,
\psi^{\omega_i, \text{snap}}_l
\psi^{\omega_i, \text{snap}}_k \, dx 
=  R_p^{\omega_i, \text{snap}} \, \overline{S} \, (R_p^{\omega_i, \text{snap}})^T.
\end{split}
\]

Here 
\[
\overline{k}(x) = \sum_{r=1}^{N_r} t_r k(x, \theta_r),
\] 
is independent of $\theta_r$ and $t_j$ are prescribed non-negative weights.  
The main objective is to use the offline space to accurately construct a set of multiscale basis functions for each $\theta_r$. At the offline stage the bilinear forms are chosen to be parameter-independent, such that there is no need to reconstruct the offline space for each $\theta_r$. 

We then  choose the smallest $M^{\omega_i}_p$ eigenvalues from Eq.~\eqref{offeig} and form the corresponding eigenvectors in the space of snapshots by setting
\[
\psi_j^{\omega_i, \text{off}} = 
(R_p^{\omega_i, \text{snap}})^T
\Psi_j^{\omega_i, \text{off}}  , 
\quad 
j = 1, \ldots, M^{\omega_i}_p.
\]  
 
% Partition of unity
Finally, we multiply the partition of unity functions $\chi_i$ by the eigenfunctions to construct the resulting basis functions
\begin{equation} 
\label{cgbasis}
\psi_{i,j} = \chi_i \psi_j^{\omega_i, \text{off}}, \quad 
1 \leq i \leq N_v, \quad  
1 \leq j \leq  M^{\omega_i}_p,
\end{equation}
where $\chi_i$ is the standard linear partition of unity function. 

Next, we define the offline space and projection matrix as 
\begin{equation} \label{cgspace}
Q_H = \text{span} \{ 
\psi_{i,j} : \, 
1 \leq i \leq N_c, \, 
1 \leq j \leq M_p 
\},  \quad 
R_p = \left[ 
\psi_{1,1} , \ldots, \psi_{N_c, M_p} 
\right]^T,
\end{equation}
where $N_c$ is number of coarse mesh nodes and $M_p = M^{\omega_i}_p$ ($i = 1,...,N_c$).

\subsection{Multiscale basis functions for displacements}

% snapshot spac
For construction of multiscale basis functions for displacements we use similar algorithm. 
We first construct a snapshot space for each parameter $\theta_r$ as a harmonic extension of $a(\theta_r; u, v)$
\begin{equation} 
\label{harmonic_ex2}
\begin{split}
a(\theta_r; \phi_{r,j}^{\omega_i, \text{snap}}, v) &= 0, \quad 
x \in \omega_i, \\
\phi_{r,j}^{\omega_i, \text{snap}} &= g_j(x), \quad 
x \in \partial \omega_i,
\end{split}
\end{equation}
where 
$g_j(x)= (\delta_{l,k}, 0, 0 )$ or $ (0, \delta_{l,k}, 0 )$ or $ (0, 0, \delta_{l,k} )$, $\forall l,k \in \textsl{J}_{h}(\omega)$ ($r = 1,...,N_r$). 

Define local snapshot space
\[
V^{\omega_i, \text{snap}} = 
\text{span}\{ 
\phi_{r,j}^{\omega_i, \text{snap}}: \, 
1 \leq r \leq N_r, \, 
0 \leq j \leq d \cdot N^{\partial \omega_i}
\},  
\]
for each subdomain $\omega_i$ and $d=2,3$. 
 
We denote the corresponding matrix of snapshot functions, again with similar notation, to be 
\[
R_u^{\omega_i, \text{snap}} = \left[ 
\phi_{1}^{\omega_i, \text{snap}}, \ldots, 
\phi_{M_{u,\text{snap}} }^{\omega_i, \text{snap}} 
\right]^T,
\]
where 
$M_{u,\text{snap}} = d \cdot N^{\partial \omega_i} \cdot N_r$ denotes the total number of functions to keep in the snapshot construction. 

% offline space
Again, we perform a dimension reduction of the space of snapshots by using an auxiliary spectral decomposition. We solve the parameter-independent eigenvalue problem in the space of snapshots
\begin{equation} 
\label{offeig2}
A^{\omega_i, \text{off}} \Phi_j^{\omega_i, \text{off}} 
= \eta_j C^{\omega_i, \text{off}} \Phi_j^{\omega_i, \text{off}},
\end{equation}
where
\[
A^{\omega_i, \text{off}} 
= R^u_{\omega_i, \text{snap}} 
\overline{A} 
(R^u_{\omega_i, \text{snap}})^T,
\quad
C^{\omega_i, \text{off}} 
= R^u_{\omega_i, \text{snap}} 
\overline{C} 
(R^u_{\text{snap}})^T,
\]
where $\overline{A}$ and $\overline{C}$ denote fine scale matrices 
\[
\overline{A} = \{a_{lk}\}, \quad 
a_{lk} = \int_{\omega_i} 
\sigma(x, \phi_l) : \varepsilon(\phi_k) \, dx, \quad
\overline{C} = \{c_{lk}\}, \quad 
c_{lk} = \int_{\omega_i} 
( \overline{\lambda} (x) + 2 \overline{\mu}(x) ) \, \phi_l \cdot  \phi_k  \, dx, 
\]\[
\sigma(x, u) = 2 \overline{\mu} (x) \varepsilon (u) 
 + \overline{\lambda}(x) \div u  \,  \mathcal{I}, 
\]
and $\phi_k$ are fine-scale basis functions. 

Further, we have 
\[
\overline{\mu}(x) = \frac{ \overline{E}(x)}{2 (1 + \eta)}, \quad
\overline{\lambda}(x) = \frac{\overline{E}(x)\eta}{(1+ \eta) ( 1- 2 \eta)},
\]
where
\[
\overline{E}(x) = \sum_{r=1}^{N_r} t_r E(x, \theta_r), 
\]
 is independent of $\theta_r$ and $t_r$ are prescribed non-negative weights.   
 
As before for the fluids flow module, at the offline stage of the mechanics the bilinear forms are chosen to be parameter-independent, such that there is no need to reconstruct the offline space for each $\theta_r$. 

We then  choose the smallest $M_u^{\omega_i}$ eigenvalues from Eq.~\eqref{offeig2} and form the corresponding eigenvectors in the space of snapshots by setting
\[
\phi_j^{\omega_i, \text{off}} = (R^u_{\text{snap}})^T \Phi_j^{\omega_i, \text{off}},  \quad 
j=1,\ldots, M_u^{\omega_i}.
\]

Finally, we multiply the linear partition of unity functions $\xi_i$ by the eigenfunctions to construct the resulting basis functions
\begin{equation} 
\label{cgbasis.mechanics}
\phi_{i,j} = \xi_i \phi_j^{\omega_i, \text{off}} , \quad
1 \leq i \leq N_c, \quad , 
1 \leq j \leq M_u^{\omega_i}.
\end{equation}

Next, we define the multiscale space and projection matrix as
\begin{equation} 
\label{cgspace.mechanics}
V_H  = \text{span} \{ 
\phi_{i,j} : \,  
1 \leq i \leq N_c \,   
1 \leq j \leq M_u \}, \quad
R_u = \left[ \phi_{1,1} , \ldots, \varphi_{N_c, M_u} \right]^T,
\end{equation}
where $M_u =  M_u^{\omega_i}$ ($\forall i = 1,...,N_c$).

\subsection{Global coupling}

The multiscale spaces  are constructed  for both the fluid and mechanics, and  we  can use them at the global level.  
Using the matrices
\[
R_p = \left[ \psi_{1,1} , \ldots, \psi_{N_c, M_p} \right]^T,
\quad \text{and} \quad 
R_u = \left[ \phi_{1,1} , \ldots, \phi_{N_c, M_u} \right]^T,
\] 
we may write matrix form for the multiscale approximation 
\begin{equation}
\label{eq:coarse}
\begin{split}
A^H(\theta) u^H + G^H p^H & = 0, \\
D^H \frac{u^H - \check{u}^H}{\tau} 
+ M^H \frac{p^H - \check{p}^H}{\tau} 
+ B^H(\theta) p^H & = F^H,
\end{split}
\end{equation}
where 
\[
A^H(\theta) = R_u A(\theta) R^T_u, \quad 
G^H = R_u G^h R^T_p,  
\]\[
B^H(\theta)  = R_p B^h(\theta)  R^T_p,  \quad
M^H = R_p M^h R^T_p,  \quad
D^H = R_p  D^h  R^T_u, \quad 
F^H = R_p F.
\]
Finally, we reconstruct solution on the fine grid by $u^{ms} = R^T_u u^H$, $p^{ms} = R^T_p p^H$.

Note that, the construction of the multiscale basis functions is performed on the offline stage and the same for any random field realization. 
Construction of basis usually contains two steps: (1) the construction of a snapshot space that will be used to compute an offline space; and (2) the construction of an offline space by performing a dimension reduction in the snapshot space. 
Multiscale basis functions for pressure and displacement described above are referred to Type 1. 
One can also use all the fine grid nodal basis as snapshots. This offline space is referred to Type 2.

\section{Bayesian Uncertainty Quantification}
%\section{Preconditioned Markov Chain Monte Carlo Method}
 
As noted earlier, in many geomechanical applications a multitude of data observations are measured, such as surface subsidence based on InSar   data collected from earth observation \cite{massonnet1998radar, burgmann2000synthetic, hanssen2001radar}. 
Due to the collection procedures and data processing times involved in the apprehension of such data, these pieces of information are often sparse in time, yet spatially global. This is because only one or two post-processed displacement fields are available per year depending on various observation and processing constraints. 
Other local sensors may stream relatively constant data of surface subsidence, but only at few sparse random spatial locations. To integrate these observations into the subsurface picture we utilize a Bayesian update framework.% refBayesLit.

We denote the set of such observations as $F_{obs}$ and want to condition the probability distribution of  random fields such as $E$ and $k$ to respect the observed data. 
We  develop an  algorithm to sample the permeability and elastic parameters given observed  data $F_{obs}$ that  include subsidence data. 
This is accomplished by the well known Bayesian formula where we may relate a new posterior probability distribution given that we have measured information from a likelihood and a prior (usually Gaussian). 
%In the context of reservoir engineering production data is extensively measured and may also be available to further constrain the flow.
 %or in the context of reservoir engineering production data is extensively measured.
 
%Next, we consider a Markov chain Monte Carlo (MCMC) algorithm \cite{robert2013monte}. 

\subsection{Subsurface Properties Parametrization} 

%\comment{ Permeability is a bit standard, but for the Lame constants $\mu$ and $\lambda$ we will have to see how they have been expanded in the KL expansions in the past. Further, is there a way to cross-correlated them since permeability is linked somewhat to rock properties.}

To parametrize the subsurface properties, we use the Karhunen-Lo{\'e}ve expansion (KLE) \cite{efendiev2005efficient, efendiev2006preconditioning}. 
Let $Y(x, \theta)$ be a stochastic process such that
\[
\mathbb{E}\left[\norm{Y}^2_{L^2(\Omega)}\right]<\infty,
\]
where $\mathbb{E}$ is expectation. 

$Y(x, \theta)$ can be expand as a general Fourier series
\[
 Y(x,\theta) = 
 \sum_{k=1}^{\infty} 
 Y_k(\theta) \varphi_k(x),
\]
where 
$Y_k(\theta) = \int_{\Omega} Y(x,\theta) \, \varphi_k(x) \, dx$ 
is the Fourier coefficient in the given $L^2$ basis $\{\varphi_k\}_{k=1}^{\infty}$. 
We wish to find an $L^2$ basis so that the Fourier coefficients are uncorrelated, or in some sense orthogonal in the expectation. 
More specifically we require  
$\mathbb{E}\left[Y_i(\theta) Y_j(\theta) \right] = 0.$ 
We denote the covariance matrix, usually assumed to be Gaussian in their correlation lengths, to be given as
$R(x,y) = \mathbb{E} \left[ Y(x,\theta) Y(y,\theta)\right].$  
Note that it is symmetric and positive definite, thus, we have that 
\[
\begin{split}
\mathbb{E} \left[Y_i(\theta) Y_j(\theta) \right] 
&=
\mathbb{E} \left[ 
\int_{\Omega} Y(x,\theta)\varphi_i(x)dx 
\int_{\Omega} Y(y,\theta)\varphi_j(y)dy
\right]\\
&=
\int_{\Omega} \int_{\Omega} R(x,y) \varphi_i(x) \varphi_j(y)\, dx \, dy 
=\delta_{ij},
\end{split}
\]
where $\{\varphi_k\}_{k=1}^{\infty}$ are mutually orthogonal eigenfunctions corresponding to $R$. Indeed, we have 
\[
\int_\Omega R(x,y) \varphi_k(y) \, dy 
= \psi_k \varphi_k(x), \quad k=1,2,\dots,
\]
where $\psi_k = \mathbb{E} \left[ Y^2_k(\theta) \right]$. 

We suppose that the Covariance structure $R(x,y)$ is of the form
\begin{equation}
R(x,y) = 
\sigma_R^2 \exp\left(
- \Delta^2
%-\frac{|x_1-x_2|^2}{2l_x}- \frac{|y_1-y_2|^2}{2l_y} 
%-\frac{|x_1-x_2|^2}{l_x^2}- \frac{|y_1-y_2|^2}{l_y^2} 
\right),
\end{equation}
with 
\[
\Delta^2 = \frac{|x_1-x_2|^2}{l_x^2} + \frac{|y_1-y_2|^2}{l_y^2},
\]
for two - dimensional case and 
\[
\Delta^2 = \frac{|x_1-x_2|^2}{l_x^2}
+ \frac{|y_1-y_2|^2}{l_y^2} 
+ \frac{|z_1-z_2|^2}{l_z^2},
\]
for three - dimensional case with correlation lengths $l_x$, $l_y$, $l_z$ and variance $\sigma_R^2$.

We denote the normalized stochastic Fourier coefficients as $\nu_k(\theta) = Y_k(\theta)/\sqrt{\psi_k}$ and write
\begin{equation}
\label{KLEnormalized}
 Y(x,\theta) = 
 \sum_{k=1}^{\infty}  \sqrt{\psi_k}  \nu_k(\theta) \varphi_k(x),
\end{equation}
where 
$\mathbb{E} \left[\nu_i(\theta) \right]=0$, 
$\mathbb{E}\left[\nu_i(\theta) \nu_j(\theta) \right] =\delta_{ij}$, 
and 
$\{\psi_k,\varphi_k(x)\}_{k=1}^{\infty}$ are the eigenpairs associated to the convolution eigenvalue problem. %\eqref{eigenR}.
We assume that the eigenvalues $\psi_k$ are ordered $\psi_1 \geq \psi_2 \geq ...$. 

In simulations to characterize stochastic process, we keep $L$-leading terms to capture most of the energy of the $Y(x, \theta)$
\begin{equation}
\label{KLEnormalizedN}
 Y_L(x,\theta) = 
 \sum_{k=1}^{L}  \sqrt{\psi_k}  \nu_k(\theta) \varphi_k(x),
\end{equation}
with following energy ratio  of the approximation
\[
e(L) = \frac{E||Y_L||^2}{E||Y||^2} 
= \frac{\sum_{k=1}^{L} \psi_k}{\sum_{k=1}^{\infty} \psi_k}. 
\]

%my_dataK[i] = math.exp(40*val)
%my_dataE[i] = 0.1*math.pow((1.0 - val)/val, 1.5)

For definition of the permeability and elastic properties, we normalize and rescale random filed to define porosity
\begin{equation}
\label{porosity}
\phi(x,\theta) = \phi(Y_L(x,\theta)), 
\end{equation}
where  $Y_L(x,\theta)$ is given by \eqref{KLEnormalizedN} with the corresponding covariance as $R$, eigenfunctions $\varphi_k$ and stochastic coefficients as $\nu_k$. 

We suppose that the permeability as a function of porosity
\begin{equation}
\label{logperm}
k(x,\theta) = \exp\left(
a \, \phi( Y_L(x,\theta) ) 
\right), 
\end{equation}
where $a > 0$.

The elasticity constants $\lambda$ and $\mu$ are given such that
\[
\mu(x,\theta) = \frac{E(x,\theta)}{2 (1 + \eta)}, \quad
\lambda(x,\theta) = \frac{E(x,\theta) \eta}{(1+ \eta) ( 1- 2 \eta)},
\]
where the elastic modulus is given by
\begin{equation}
\label{E.stoch}
E(x,\theta) = 
b \left( 
\frac{1.0 - \phi( Y_L(x,\theta) ) }{\phi( Y_L(x,\theta) ) }
\right)^n,
\end{equation}
where $n = 1.5$ and $b > 0$ \cite{yang2013simulator}.

\subsection{Observable Surface Data}

%For black-oil two-phase flow refYalchinMCMCPapers, the fractional flow is used as the observable data. 
%This is the fraction of the flow that is oil at time $t$ at the outlet boundary (or well bore in the case of oil fields).
%%
%%
%This is
%$$
%F(t)=1-\frac{\int_{\partial \Omega_{out}} f(S)v\cdot n dl }{\int_{\partial \Omega_{out}}v\cdot n dl},
%$$
%here $S$ is saturation of water, $f$ is the flux function, and $v$ is the Darcy velocity. 
%Let $F_{k}(t)$ be the fractional flow that is computed using permeability $k$, and $F(t)$ the observed fractional flow. Then, we define the data misfit
%\begin{align}
% \norm{F(\cdot)-F_k(\cdot)}^2=\int_{0}^T |F(t)-F_{k}(t)|^2 dt.
%\end{align}

For the case of surface subsidence, the surface displacement  is used as the observable data
\[
u_{obs}(x,t), \quad x\in \partial \Omega_{surf}
\]
where $\partial \Omega_{surf}$ the surface boundary.  
This data is usually dynamic, but sparse in the spatial extent, with the occasional inclusion of global deformation at sparse time interval snapshots. % refNeedToFind.

\begin{figure}[h!]
	\centering
	\includegraphics[width=0.8 \textwidth]{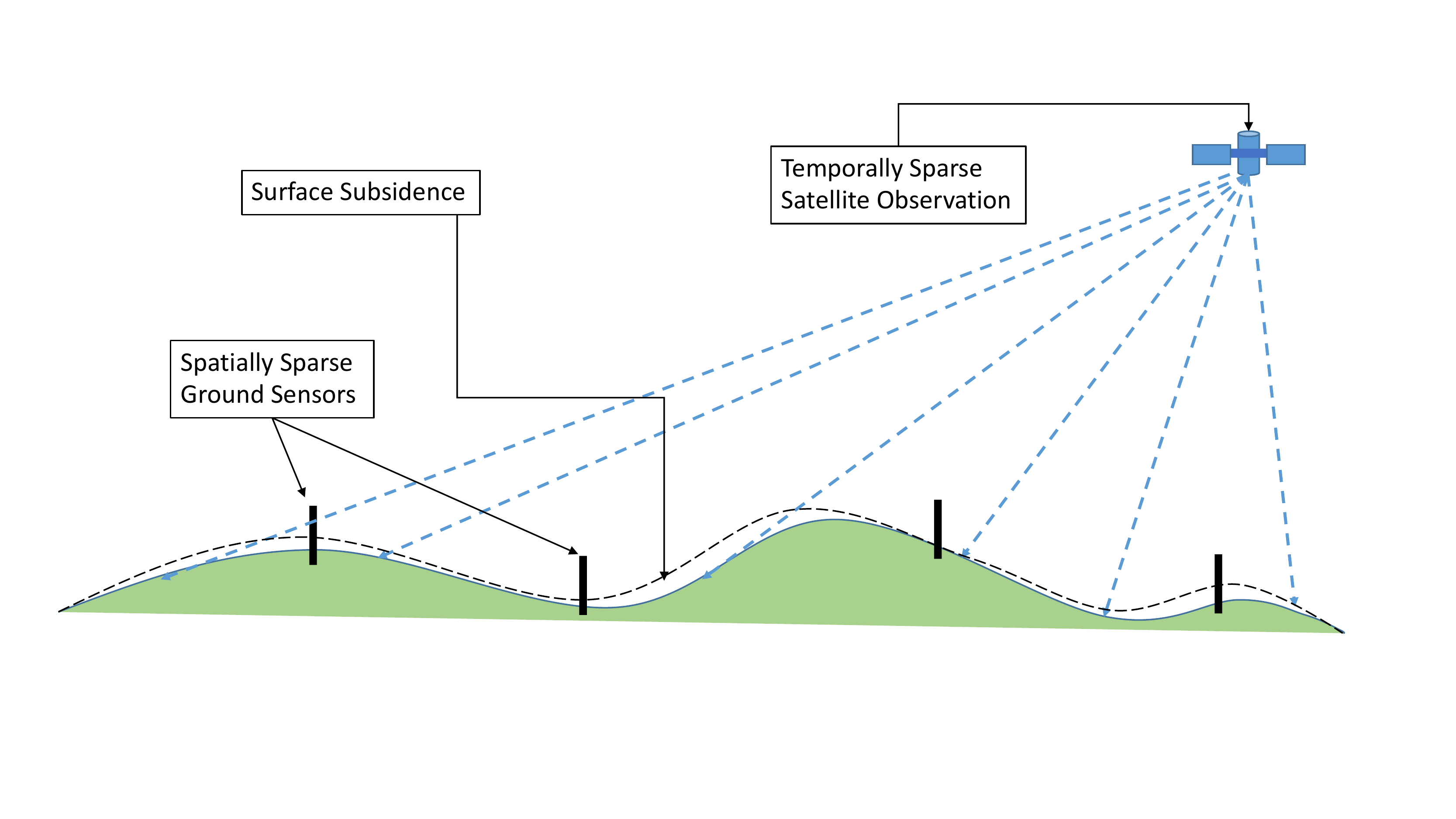}
	\caption{Observable surface data of geomechanical subsidence. Spatially sparse, but temporally constant ground motion sensors. Temporally sparse, but spatially global satellite observation. }
	\label{Data}
\end{figure}

The dynamic data comes from surface sensors at sparse locations throughout the area of interest that continuously monitor displacement in time at sparse locations. While, the global information is available from InSar satellite observation data. %refGeospatialData. 
However, due to the nature of the collection of this data and the length of time to process the large data sets involved, only snapshots of data are usually available in practical settings.  
To make this more precise, we suppose that at distinct points $x_{i}\in \partial \Omega_{surf}$,  for $i=1,\cdots N^s_{data}$, $u_{obs}(x_i,t)$ is known for all time $t\in[0,T]$. 
Further, we suppose that at snapshot times $t_j$, $j=1,\cdots N^t_{data}$, $u_{obs}(x,t_j)$ is known for all time $x\in \partial \Omega_{surf}$

Suppose that $u_{\theta}(x,t)$ is the displacement computed using properties $E(\theta)$ and $k(\theta)$, then we define our data misfit as 
\begin{equation}
\norm{F_{obs}-F_\theta}^2 = 
\sum_{i=1}^{N^s_{data}} \int_{0}^T 
|u_{obs}(x_i,t)-u_{\theta}(x_i,t)|^2 dt
+
\sum_{j=1}^{N^t_{data}} \int_{\partial \Omega_{surf}} 
|u_{obs}(x,t_j)-u_{\theta}(x,t_j)|^2 ds, 
\end{equation}
where $ds$ is the surface measure.

\subsection{Single-stage Metropolis - Hastings Algorithm}

%We denote $\theta = (E, k)$ to be the pair of random fields in \eqref{eq:main}. 
The posterior distribution $\pi(\theta)$ can be expressed as 
\begin{equation}
\label{Bayesian} 
\pi(\theta) = p(\theta|F_{obs}) 
\propto p(F_{obs}| \theta) \, p(\theta),
\end{equation}
where $ p(F_{obs}| \theta)$ is the likelihood and  $p(\theta)$ is the prior distribution. % (again normally Gaussian or log-normal). 

Given a set of observed data $F_{obs}$ and heterogeneous field parameter data $\theta$ we write $F(\theta)$ to be the corresponding simulated data (on the fine-grid) i.e. subsidence displacement or pressures. 
Due to uncertainty in the model and measurements, we suppose that the error is such that 
\[
F_{obs} - F(\theta) = \varepsilon,
\]
where $\varepsilon$ is random normal with mean $0$ and variance $\sigma_f^2$, or in standard notation  $\varepsilon$ follows ${\cal N}(0,\sigma_f^2)$. %{\cal I})$.
Thus, we will write the posterior as 
\begin{equation}
\label{likelihood}
p(\theta|F_{obs})\propto \exp\left(-\frac{\norm{F(\theta)-F_{obs}}^2}{\sigma_f^2}\right)p(\theta).
\end{equation}

%Here we outline both the standard Metropolis-Hastings MCMC algorithm and the Two-Stage MCMC algorithm.
%\subsubsection{Single-Stage Metropolis-Hastings Algorithm}
%We begin with the standard Metropolis-Hastings algorithm. 

Let $q(\theta|\theta_n)$ be the instrumental distribution that is used to choose the next fields $\theta$ given the previous properties $\theta_n$.
%This is usually the random walk or could be something more complex such as a Langevin sampler. % refLAngevinLit.
We consider a random walk samplers for the transitional probability distribution  $q(\theta|\theta_n)$ and set  $q(\theta|\theta_n) = \theta_n + \delta \cdot r$, where $r$ is a Gaussian random variable with zero mean and
variance 1.

\textit{The single - stage Metropolis - Hastings algorithm:} %proceeds as follows: % for $n=0,2,3,\dots$
\begin{enumerate}
\item Given $\theta_n$, generate new proposal $\theta$, from $q(\theta|\theta_n)$.
\item Solve  forward problem, using  $\theta$, generate observable $F(\theta)$.
\item Determine the acceptance probability from 
\begin{equation}
\label{ProbAcc}
 \text{Prob}(\theta_n|\theta) = \text{min} \left( 1, 
 \frac{q(\theta_n|\theta) \, p(\theta|F_{obs})}{q(\theta|\theta_n) \, p(\theta_n|F_{obs})}
 \right), 
\end{equation}
%Construct  a small dimensional offline space by performing dimension reduction in the space of local snapshots.
and take 
\[
\theta_{n+1} = \left\{\begin{matrix}
\theta & \text{with probability} \quad \text{Prob}(\theta_n|\theta), \\ 
\theta_n & \text{with probability} \quad 1-\text{Prob}(\theta_n|\theta).
\end{matrix}\right.
\]
%Accept, let $\theta_{n+1}=\theta$, with probability $\text{Prob}(\theta_n|\theta).$
%$\hspace{.1in}$ 
%Reject, let $\theta_{n+1}=\theta_n,$ with probability $1-\text{Prob}(\theta_n|\theta).$
%Use small dimensional offline space to find the solution of a coarse-grid problem for any force term and/or boundary condition.
\end{enumerate}
%\comment{Note that in the above algorithm for symmetric $q(\cdot|\cdot)$, these terms cancel.} 
For symmetric $q(\cdot|\cdot)$, we have
\[
 \text{Prob}(\theta_n|\theta) 
 = \text{min} \left( 1, 
 \exp \left( -\frac{E(\theta) - E(\theta_n) }{\sigma_f^2} \right)
 \right), 
\]
with $E(\theta) = \norm{F(\theta)-F_{obs}}^2$ and $E(\theta_n) = \norm{F(\theta_n)-F_{obs}}^2$.

%\subsubsection{Two-Stage Metropolis-Hastings Algorithm}
It is well known that the above algorithm is very costly due to the expense of having to solve the fine-scale solution and the low acceptance rate of new proposals. 
A method to save computational cost is to use a coarse-scale model or/and machine learning technique as a preconditioning phase to cheaply reject proposals. Then, if the proposals is accepted on the cheap first-stage, run a fine-scale simulation and generate new acceptance-rejection probabilities. This two-stage procedure has been utilized in various other applications \cite{efendiev2005efficient, efendiev2006preconditioning, ginting2011application}.

\subsection{Two-stage Metropolis - Hastings Algorithm}

Let 
\begin{equation}
\label{likelihood}
p^*(\theta|F_{obs}) \propto 
\exp\left( 
-\frac{E^*(\theta)}{\sigma_c^2}
\right),
\end{equation}
where  $E^*(\theta) = \norm{F^*(\theta)-F_{obs}}^2$ and $\sigma_c^2$ is the precision associated with the approximate model. 

Here $F^*(\theta)$  denotes the simulated data that generated using 
\begin{itemize}
\item \textit{MS}: coarse scale simulations (multiscale solver),
\item \textit{ML}: machine learning techniques,
\end{itemize} 

\textit{The two - stage Metropolis - Hastings algorithm:}
\begin{enumerate}
\item Given $\theta_n$, generate new proposal $\theta'$, from $q(\theta'|\theta_n)$.
\item First stage:
\begin{enumerate}
\item Generate observable $F^*(\theta')$ for given $\theta'$ using forward multiscale solver (\textit{MS}) or machine learning technique (\textit{ML}).
\item Determine the acceptance probability from 
\begin{equation}
\label{ProbAccCoarse}
 g(\theta_n|\theta') = 
 \text{min}\left( 1, 
 \frac{q(\theta_n|\theta')p^*(\theta'|F_{obs})}{q(\theta'|\theta_n)p^*(\theta_n|F_{obs})}
 \right),
\end{equation}
and take  
\[
\theta =
\left\{\begin{matrix}
\theta' & \text{with probability} \quad  g(\theta_n|\theta), \\ 
\theta_n & \text{with probability} \quad 1- g(\theta_n|\theta).
\end{matrix}\right.
\]
%real proposal,  $\theta_{}=\theta'$, with probability $ g(\theta_n|\theta)$, $\theta_{}=\theta_n,$ with probability $1- g(\theta_n|\theta).$
%Use small dimensional offline space to find the solution of a coarse-grid problem for any force term and/or boundary condition.
\end{enumerate}
\item  Second stage, if accepted:
\begin{enumerate}
\item Solve fine-scale forward problem using $\theta$ and generate observable $F(\theta)$.
\item Determine the acceptance probability from 
\begin{equation}
\label{proposalprobCS}
 \rho(\theta_n|\theta) = \text{min}\left(1, 
 \frac{p^*(\theta_n|F_{obs})p(\theta|F_{obs})}{p^*(\theta|F_{obs})p(\theta_n|F_{obs})}
 \right),
\end{equation}
and take
\[
\theta_{n+1} =
\left\{\begin{matrix}
\theta & \text{with probability} \quad   \rho(\theta_n|\theta), \\ 
\theta_n & \text{with probability} \quad 1 - \rho(\theta_n|\theta).
\end{matrix}\right.
\]
%Write the effective instrumental distribution
%\[
%Q(\theta|\theta_n)=g(\theta_n|\theta)q(\theta|\theta_n)+\left(1-\int g(\theta_n|\theta)q(\theta|\theta_n)dk \right)\delta_{\theta_n}(\theta).
%\]
%Accept $\theta$ as proposal with probability 
%\begin{align}\label{proposalprob}
% \rho(\theta_n|\theta)=\text{min}\left(1, \frac{Q(\theta_n|\theta)p(\theta|F_{obs})}{Q(\theta|\theta_n)p(\theta_n|F_{obs})}\right),
%\end{align}
%that is accepted  $\theta_{n+1}=\theta$, with probability $ \rho(\theta_n|\theta)$,
% $\theta_{n+1}=\theta_n,$ with probability $1- \rho(\theta_n|\theta).$
\end{enumerate}
\end{enumerate}

%Note here that we need note compute the effective instrumental distribution as seen in refPRCMCMC, we may write 
%$$
%\frac{Q(\theta_n|\theta)}{Q(\theta|\theta_n)}= \frac{p^*(\theta_n|F_{obs})}{p^*(\theta|F_{obs})},
%$$
%and thus, may compute instead \eqref{proposalprob} as
%\begin{align}\label{proposalprobCS}
% \rho(\theta_n|\theta)=\text{min}\left(1, \frac{p^*(\theta_n|F_{obs})p(\theta|F_{obs})}{p^*(\theta|F_{obs})p(\theta_n|F_{obs})}\right),
%\end{align}
%which is much easier to compute compared to the integration of the effective instrumental distribution.

Here, we have
\[
 g(\theta_n|\theta') 
  = \text{min} \left( 1, 
 \exp \left( -\frac{E^*(\theta') - E^*(\theta_n) }{\sigma_c^2} \right)
 \right), 
\]
and 
\[
 \rho(\theta_n|\theta) 
 = \text{min}\left(1, 
 \frac{p^*(\theta_n|F_{obs})p(\theta|F_{obs})}{p^*(\theta|F_{obs})p(\theta_n|F_{obs})}
 \right)
   = \text{min} \left( 1, 
 \exp \left( 
 -\frac{ E(\theta) - E(\theta_n) }{\sigma_f^2} 
 +\frac{ E^*(\theta) - E^*(\theta_n) }{\sigma_c^2} 
 \right)
 \right). 
\]
If $E^*(\theta)$ is strongly correlated with $E(\theta)$, then the acceptance probability $\rho$ could be close to 1 for certain choice of $\sigma_c$ \cite{efendiev2006preconditioning}.

% MS
To Generate observable $E^*(\theta')$ for given $\theta'$ using forward multiscale solver (\textit{MS}), we use following steps:
\begin{itemize}
\item \textit{Offline stage}. Construction of the projection matrices $R_p$ and $R_u$.
\item \textit{Online stage}. For given $\theta'$:
\begin{enumerate}
\item generate a coarse grid system using preconstructed projection matrices,
\item solve time-depended coarse grid system \eqref{eq:coarse},
\item find observable data $u_{obs}(\theta')$ on the surface boundary from multiscale solution,
\item generate $E^*(\theta')$.
\end{enumerate}
\end{itemize}

% ML
In machine learning technique (\textit{ML}), we have following steps:
\begin{itemize}
\item \textit{Offline stage}. Dataset generation and neural networks training.
\item \textit{Online stage}. For given $\theta'$:
\begin{enumerate}
\item predict observable data $u_{obs}(\theta')$ on the surface boundary using trained neural networks,
\item generate $E^*(\theta')$.
\end{enumerate}
\end{itemize}
We construct a neural network (NN) for displacements in each direction $x, y$ and $z$. Therefore, we train two NNs for two-dimensional problems (2D) and three NNs for three-dimensional problems (3D) \cite{vasilyeva2020learning}. In \textit{ML}, we directly predict an observable data for given $\theta'$, without solution of the forward problem.

To construct a neural network for prediction of the observable data, we generate a dataset by the solution of the forward problems using a multiscale solver for some number of random coefficient $\theta_r$ ($r=1,...,N_r$)
\[
\text{ Dataset: }
\{ (X_r, Q_r), \, r = 1,...,N_r \}.
\]
where $X_r = \theta_r$ and $Q_r = u_{obs}(\theta_r)$ are the input data and output data. 
The input field is represented as a two-dimensional array for the two-dimensional problem and as a three-dimensional array for three-dimensional problems. The dataset is re-scaled to fall within the range $0$ to $1$.

We use a convolutional neural network, where several convolutional, pooling, and activation layers are stacked with several fully-connected layers with dropout. Training of the machine learning algorithms is performed using mean square loss function (MSE). Implementation of the machine learning method is based on the library Keras \cite{keras} with TensorFlow backend \cite{tensorflow} and performed on the GPU. The machine-learning algorithm learns dependence between global random coefficients ($\theta$) and observable data $u_{obs}(\theta)$.

\section{Numerical results}

In order to investigate the proposed method, we perform the following tests:
\begin{itemize}
\item In Section \ref{ss1}, we consider multiscale solver for three test cases in 2D and 3D formulations. We show the relative errors between reference (fine-grid) solution and GMsFEM solution with different numbers of the multiscale basis functions. Also, we present a solution time of the proposed multiscale method. 

\item  In Section \ref{ss2}, we consider MCMC method with GMsFEM preconditioning.
We start with tests of multiscale solver, where we present errors for pressure and displacements for 100 random realizations.  
Next, we demonstrate correlation  between $E(\theta)$ (fine grid solver) and $E^*(\theta)$ (multiscale solver with different number of basis functions). 
After that, we present results for the two-stage MCMC method with GMsFEM preconditioning for 2D and 3D test problems.  We investigate the influence of the method parameters on the numbers of the accepted fields and number of the fields that pass the first stage in the two-dimensional formulation. Some figures of accepted fields with corresponded solutions are presented.

\item  In Section \ref{ss3}, we consider the preconditioned MCMC method using a machine learning technique. 
We start with the demonstration of the neural network architectures for 2D and 3D problems, and we use them for the prediction of the observable data. 
Next, we demonstrate correlation  between $E(\theta)$ (fine grid solver) and $E^*(\theta)$ (machine learning method).  
After that, we present results for the preconditioned MCMC method, where we took an observable data for three test cases considered in Section \ref{ss1}. We shown results for both preconditioning approaches:  (\textit{MS}) multiscale solver based on GMSFEM and (\textit{ML}) machine learning technique. Finally, we discuss the computational advantage of the method. 
\end{itemize}

\subsection{Multiscale method}\label{ss1}

In this section, we present numerical examples to demonstrate the performance of the multiscale method for computing the solution of the poroelasticity problem with random heterogeneous properties.

% 2d
\begin{figure}[h!]
\centering % 1 14 7
    \begin{subfigure}[b]{0.85\textwidth}
        \includegraphics[width=1.0\linewidth]{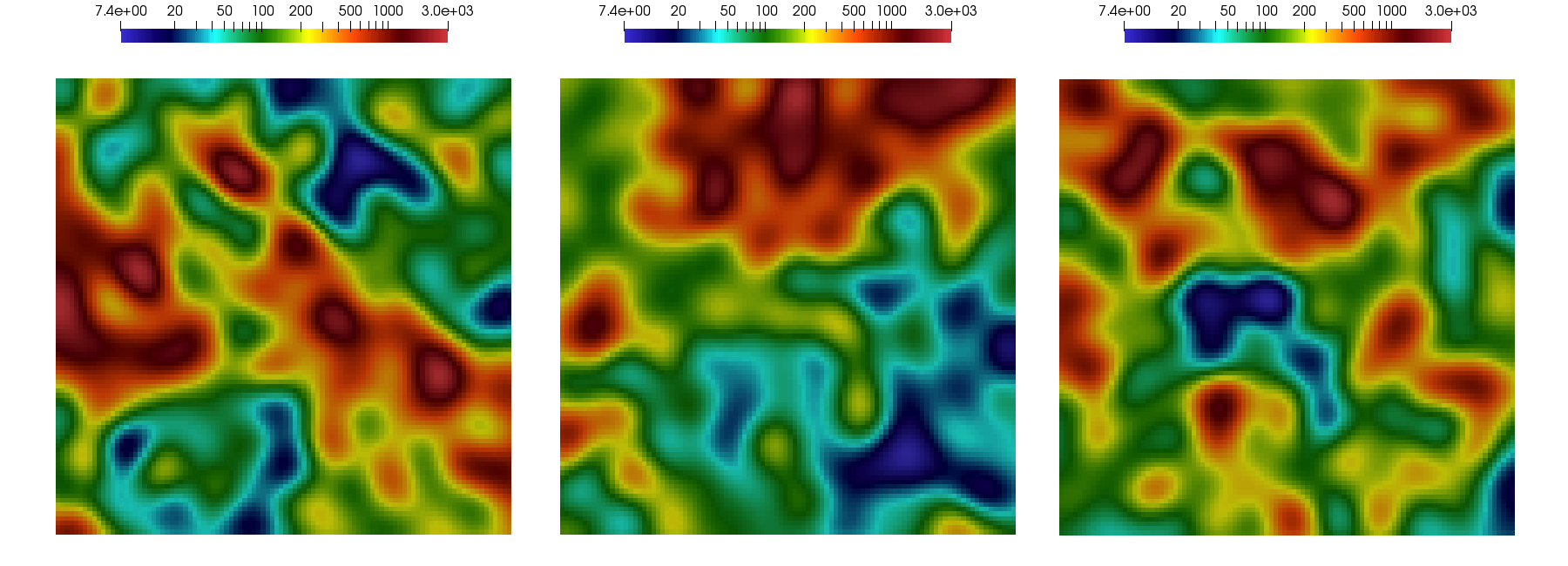}
        \caption{$k(x)$ for Case 1, 2 and 3 (from left to right). }
    \end{subfigure}\\
    \vspace{10pt}
    \begin{subfigure}[b]{0.85\textwidth}
        \includegraphics[width=1.0\linewidth]{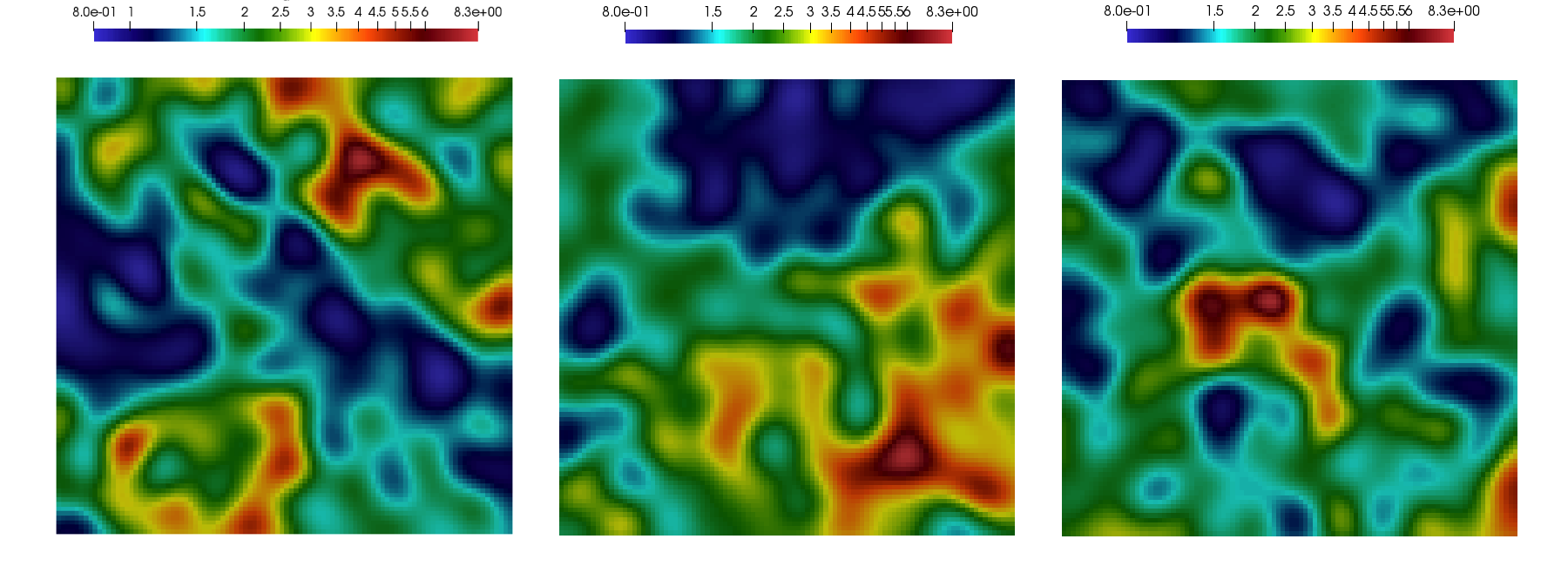}
        \caption{$E(x)$ for Case 1, 2 and 3 (from left to right). }
    \end{subfigure}
\caption{Random permeability $k$ and elastic modulus $E$. 
Two - dimensional problem (2D).  
Case 1, 2 and 3 (from left to right). 
(a) permeability, $k(x)$. 
(b) elastic modulus, $E(x)$. }
 \label{fig:2d-kE}
\end{figure} 

\begin{figure}[h!]
\centering % 1 14 7
    \begin{subfigure}[b]{0.95\textwidth}
        \includegraphics[width=1.0\linewidth]{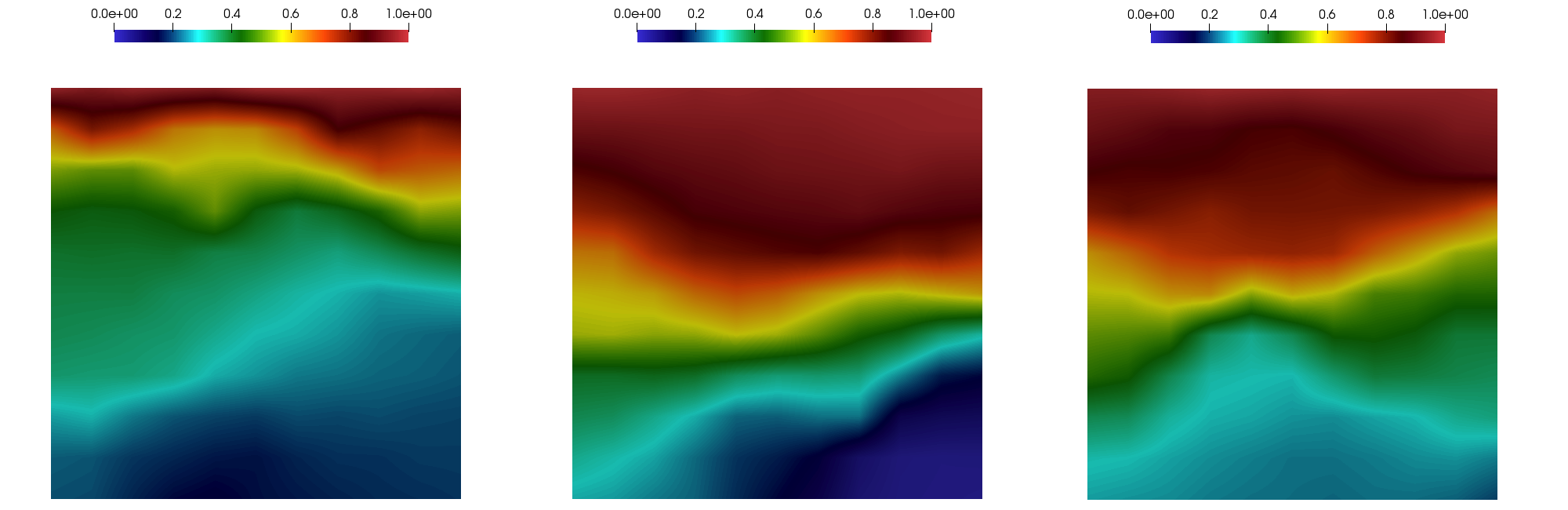}
        \caption{$p$ for Case 1, 2 and 3 (from left to right). }
    \end{subfigure}\\
    \vspace{10pt}
    \begin{subfigure}[b]{0.95\textwidth}
        \includegraphics[width=1.0\linewidth]{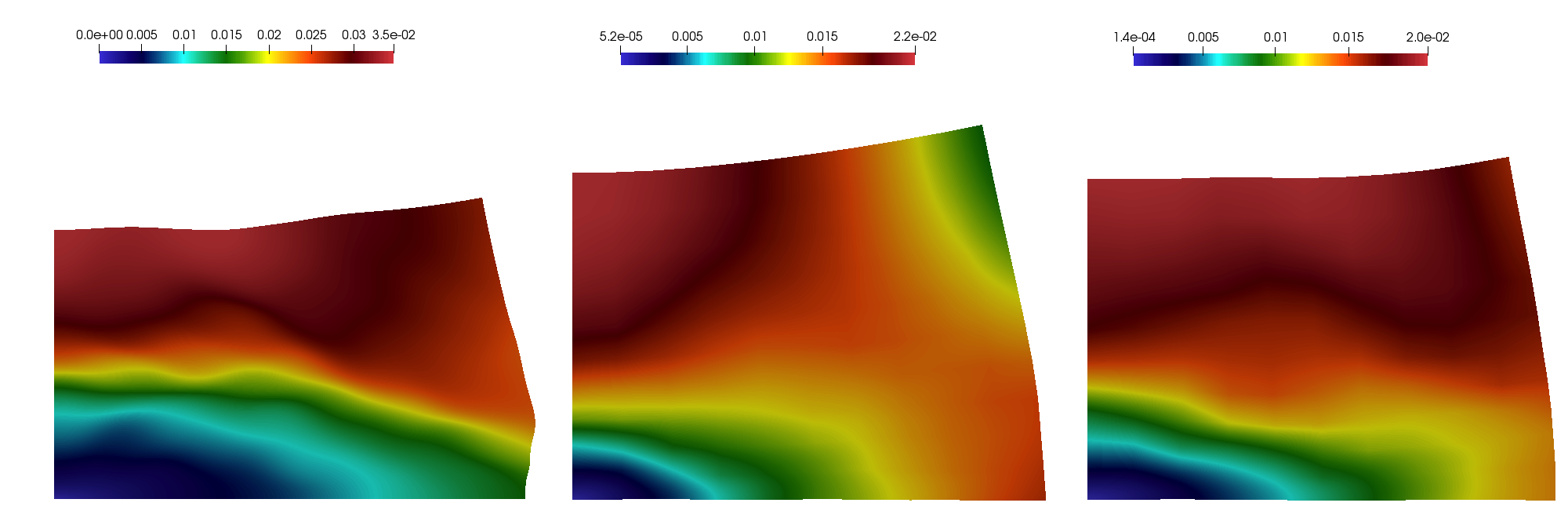}
        \caption{$u_m$ for Case 1, 2 and 3 (from left to right). }
    \end{subfigure}
\caption{Numerical results for random permeability $k$ and elastic modulus $E$. Case 1, 2 and 3 (from left to right). 
Fine grid solution of two - dimensional problem (2D). 
(a) pressure, $p$. 
(b) displacement, $u_m$. }
\label{fig:2d-up}
\end{figure} 

% 3d
\begin{figure}[h!]
\centering % 1 14 7
    \begin{subfigure}[b]{0.95\textwidth}
        \includegraphics[width=1.0\linewidth]{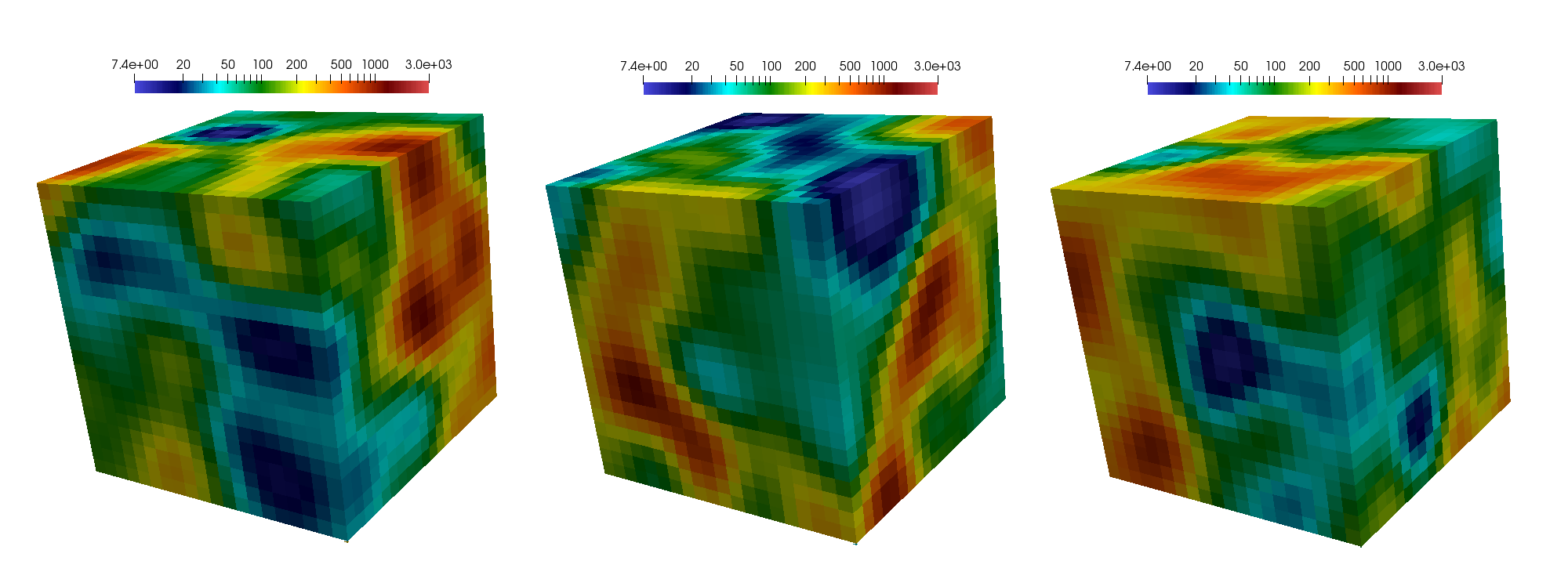}
        \caption{$k(x)$ for Case 1, 2 and 3 (from left to right). }
    \end{subfigure}\\
    \vspace{10pt}
    \begin{subfigure}[b]{0.95\textwidth}
        \includegraphics[width=1.0\linewidth]{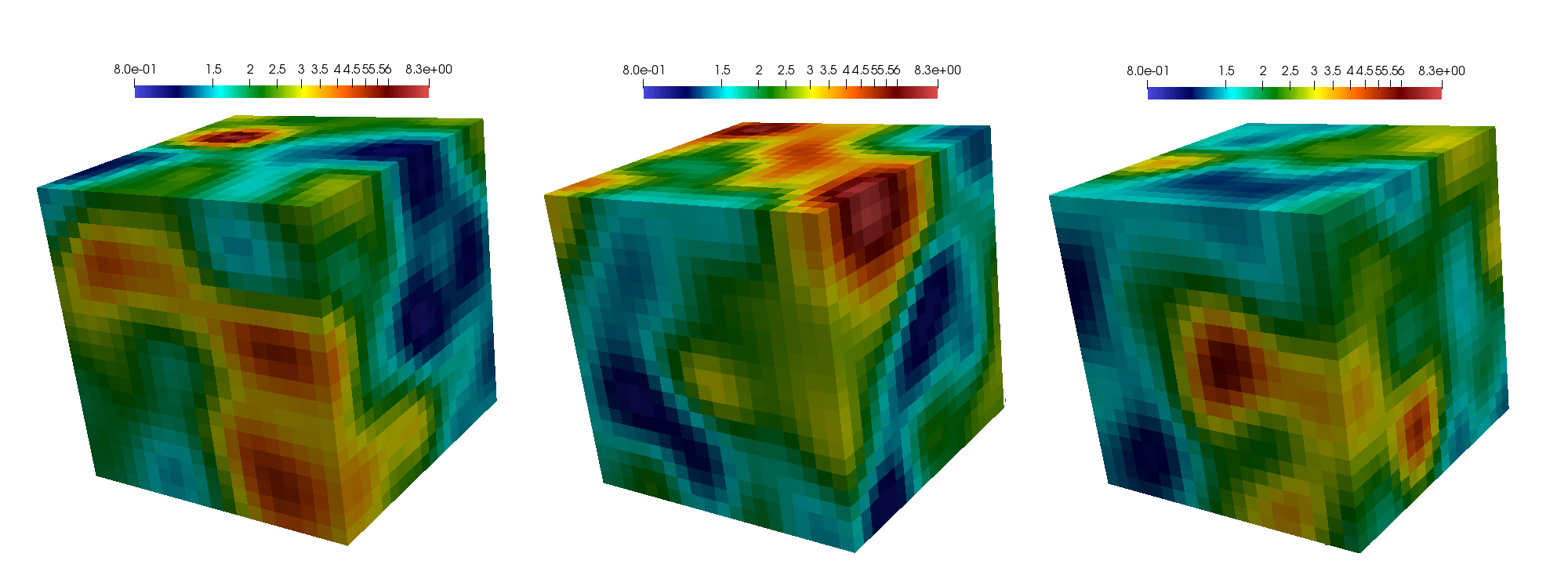}
        \caption{$E(x)$ for Case 1, 2 and 3 (from left to right). }
    \end{subfigure}
\caption{Random permeability $k$ and elastic modulus $E$. 
Three - dimensional problem (3D).  
Case 1, 2 and 3 (from left to right). 
(a) permeability, $k(x)$. 
(b) elastic modulus, $E(x)$. }
 \label{fig:3d-kE}
\end{figure} 

\begin{figure}[h!]
\centering % 1 14 7
    \begin{subfigure}[b]{0.95\textwidth}
        \includegraphics[width=1.0\linewidth]{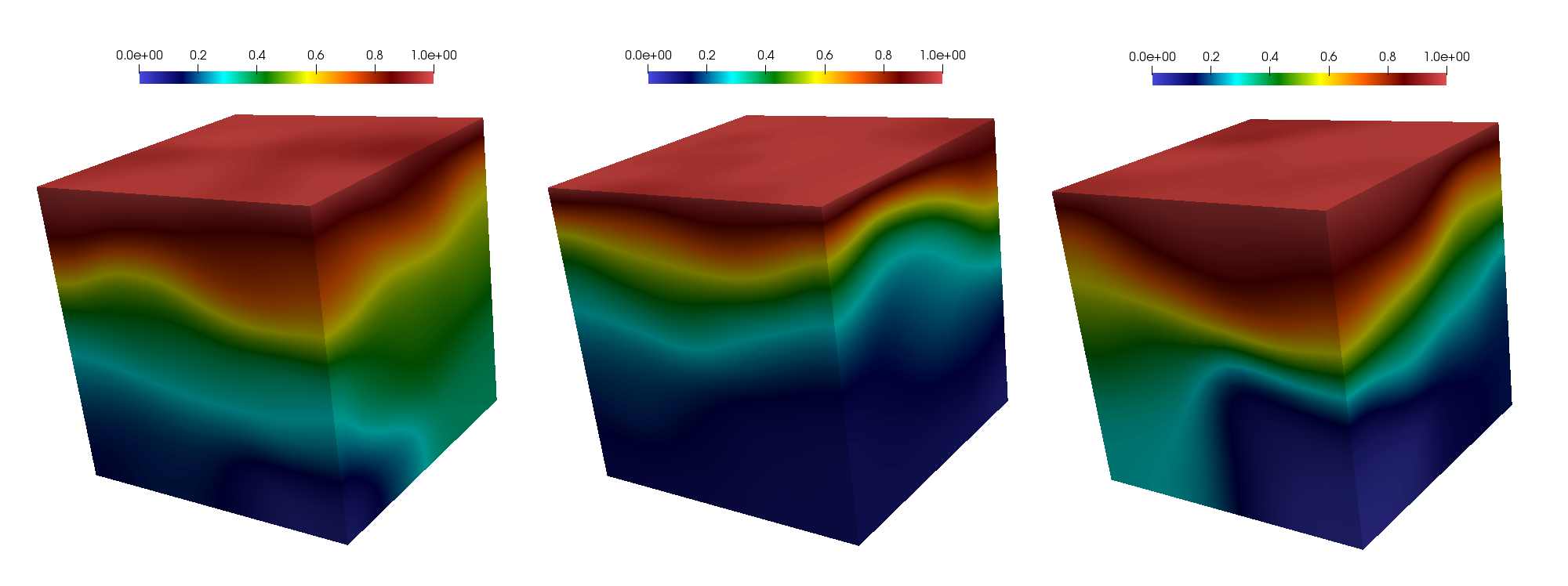}
        \caption{$p$ for Case 1, 2 and 3 (from left to right). }
    \end{subfigure}\\
    \vspace{10pt}
    \begin{subfigure}[b]{0.95\textwidth}
        \includegraphics[width=1.0\linewidth]{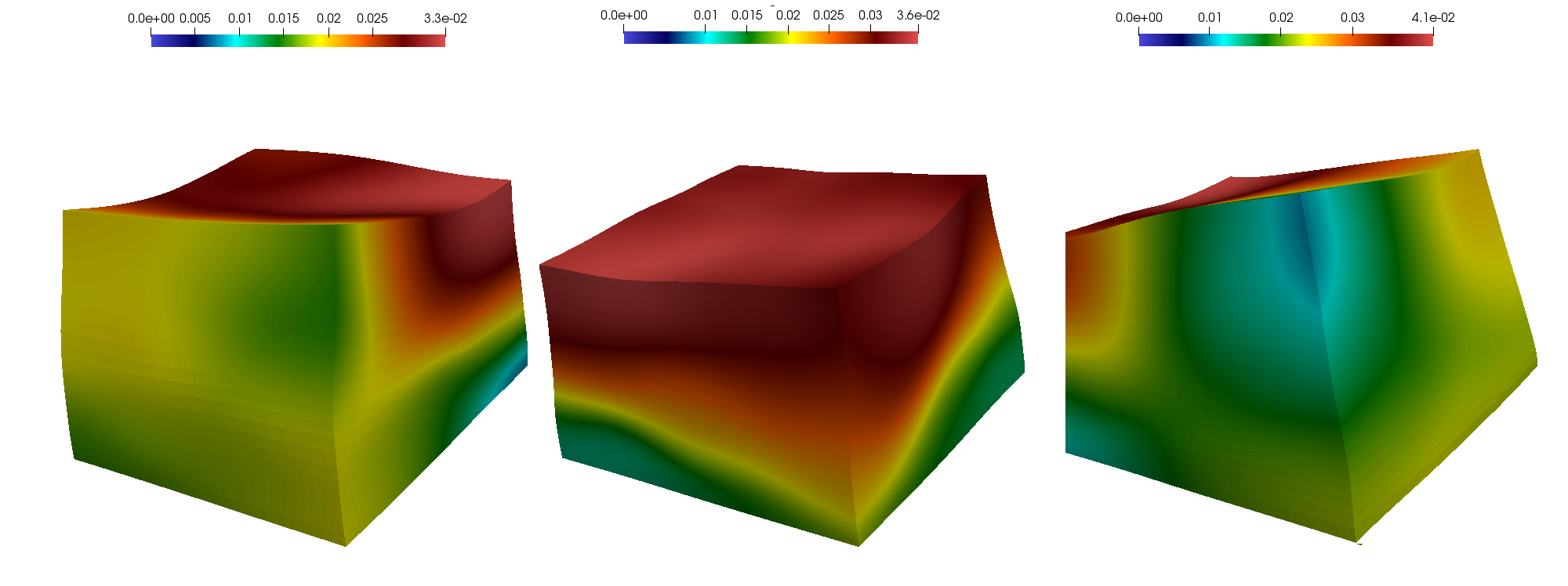}
        \caption{$u_m$ for Case 1, 2 and 3 (from left to right). }
    \end{subfigure}
\caption{Numerical results for random permeability $k$ and elastic modulus $E$. Case 1, 2 and 3 (from left to right). 
Fine grid solution of three - dimensional problem (3D). 
(a) pressure, $p$. 
(b) displacement, $u_m$. }
 \label{fig:3d-up}
\end{figure} 

\begin{table}[h!]
\begin{center}
\begin{tabular}[hp]{|c|cc|c|cc|}
\hline
\multicolumn{6}{|c|}{2D}  \\
\hline
$M_+$  & $M_p$ & $M_u$ 
& $DOF_c$ & $e_p$ (\%) & $e_u$ (\%) \\
\hline 
\multicolumn{6}{|c|}{Case 1}  \\ % 1
\hline
0	& 1	& 2 	& 363 		&  7.532  &  6.413	\\
1	& 2	& 3 	& 605 		&  4.800  &  4.048	\\
2	& 3	& 4 	& 847 		&  1.368  &  2.076	\\
3	& 4	& 5 	& 1089 	&  0.940  &  1.634	\\
4	& 5	& 6 	& 1331 	&  0.635  &  1.146	\\
6	& 7	& 8 	& 1815 	&  0.303  &  0.840	\\
8	& 9	& 10 & 2299 	&  0.182  &  0.661	\\
\hline
\multicolumn{6}{|c|}{Case 2}  \\ % 14
\hline
0	& 1	& 2 	& 363 		&  2.142  &  2.953	\\
1	& 2	& 3 	& 605 		&  1.151  &  2.127	\\
2	& 3	& 4 	& 847 		&  0.641  &  1.649	\\
3	& 4	& 5 	& 1089 	&  0.529  &  1.309	\\
4	& 5	& 6 	& 1331 	&  0.234  &  0.822	\\
6	& 7	& 8 	& 1815 	&  0.117  &  0.552	\\
8	& 9	& 10 & 2299 	&  0.068  &  0.440	\\
\hline
\multicolumn{6}{|c|}{Case 3}  \\ % 7
\hline
0	& 1	& 2 	& 363 		&  4.635  &  6.847	\\
1	& 2	& 3 	& 605 		&  2.760  &  4.303	\\
2	& 3	& 4 	& 847 		&  0.778  &  2.005	\\
3	& 4	& 5 	& 1089 	&  0.492  &  1.579	\\
4	& 5	& 6 	& 1331 	&  0.343  &  0.993	\\
6	& 7	& 8 	& 1815 	&  0.198  &  0.694	\\
8	& 9	& 10 & 2299 	&  0.124  &  0.535	\\
\hline
\end{tabular}\,\,\,\,\,\,
\begin{tabular}[hp]{|c|cc|c|cc|}
\hline
\multicolumn{6}{|c|}{3D}  \\
\hline
$M_+$  & $M_p$ & $M_u$ 
& $DOF_c$ & $e_p$ (\%) & $e_u$ (\%) \\
\hline
\multicolumn{6}{|c|}{Case 1}  \\ % 1
\hline
0	& 1	& 3 	& 864 		&  4.369  &  5.477	\\
1	& 2	& 4 	& 1296 	&  2.727  &  4.343	\\
2	& 3	& 5 	& 1728 	&  1.914  &  3.716	\\
3	& 4	& 6 	& 2160 	&  1.356  &  3.220	\\
4	& 5	& 7 	& 2592 	&  0.853  &  2.672	\\
6	& 7	& 9 	& 3456 	&  0.453  &  1.904	\\
8	& 9	& 11 & 4320 	&  0.352  &  1.619	\\
\hline
\multicolumn{6}{|c|}{Case 2}  \\ % 14
\hline
0	& 1	& 3 	& 864 		&  6.000  &  3.653	\\
1	& 2	& 4 	& 1296 	&  5.194  &  3.235	\\
2	& 3	& 5 	& 1728 	&  3.185  &  2.700	\\
3	& 4	& 6 	& 2160 	&  2.144  &  2.455	\\
4	& 5	& 7 	& 2592 	&  0.912  &  2.114	\\
6	& 7	& 9 	& 3456 	&  0.571  &  1.693	\\
8	& 9	& 11 & 4320 	&  0.452  &  1.443	\\
\hline
\multicolumn{6}{|c|}{Case 3}  \\ %7
\hline
0	& 1	& 3 	& 864 		&  3.430  &  4.520	\\
1	& 2	& 4 	& 1296 	&  2.752  &  3.997	\\
2	& 3	& 5 	& 1728 	&  1.709  &  3.318	\\
3	& 4	& 6 	& 2160 	&  1.129  &  3.071	\\
4	& 5	& 7 	& 2592 	&  0.730  &  2.841	\\
6	& 7	& 9 	& 3456 	&  0.487  &  2.331	\\
8	& 9	& 11 & 4320 	&  0.377  &  1.905	\\
\hline
\end{tabular}
\end{center}
\caption{Numerical results for random permeability $k$ and elastic modulus $E$. 
Case 1, 2 and 3 (from left to right). 
Type 1 multiscale basis functions, $M_p = 1 + M_+$, $M_u = d + M_+$ and $DOF_c = (M_p + M_u) \cdot N_c$. 
Left:  two - dimensional problem, $d = 2$ ($N_c =121$ and $DOF_f = 30603$). 
Right: three  - dimensional problem, $d = 3$ ($N_c = 216$ and $DOF_f = 37044$) }
\label{tab:err-t1}
\end{table}

\begin{table}[h!]
\begin{center}
\begin{tabular}[hp]{|c|cc|c|cc|}
\hline
\multicolumn{6}{|c|}{2D}  \\
\hline
$M_+$  & $M_p$ & $M_u$ 
& $DOF_c$ & $e_p$ (\%) & $e_u$ (\%) \\
\hline 
\multicolumn{6}{|c|}{Case 1}  \\ % 1
\hline
0	& 1	& 2 	& 363 		&  7.532  &  6.413	 \\
1	& 2	& 3 	& 605 		&  5.737  &  4.669	 \\
2	& 3	& 4 	& 847 		&  2.396  &  2.793	 \\
3	& 4	& 5 	& 1089 	&  1.548  &  1.796	 \\
4	& 5	& 6 	& 1331 	&  1.102  &  1.300	 \\
6	& 7	& 8 	& 1815 	&  0.639  &  0.481	 \\
8	& 9	& 10 & 2299 	&  0.364  &  0.272	 \\
\hline
\multicolumn{6}{|c|}{Case 2}  \\ % 14
\hline
0	& 1	& 2 	& 363 		&  2.142  &  2.953	 \\
1	& 2	& 3 	& 605 		&  1.219  &  1.847	 \\
2	& 3	& 4 	& 847 		&  0.703  &  1.430	 \\
3	& 4	& 5 	& 1089 	&  0.596  &  1.085	 \\
4	& 5	& 6 	& 1331 	&  0.380  &  0.708	 \\
6	& 7	& 8 	& 1815 	&  0.180  &  0.220	 \\
8	& 9	& 10 & 2299 	&  0.101  &  0.146	 \\
\hline
\multicolumn{6}{|c|}{Case 3}  \\ % 7
\hline
0	& 1	& 2 	& 363 		&  4.635  &  6.847 	\\
1	& 2	& 3 	& 605 		&  2.397  &  4.294	 \\
2	& 3	& 4 	& 847 		&  0.842  &  2.596	 \\
3	& 4	& 5 	& 1089 	&  0.570  &  1.641	 \\
4	& 5	& 6 	& 1331 	&  0.379  &  1.109	 \\
6	& 7	& 8 	& 1815 	&  0.218  &  0.336	 \\
8	& 9	& 10 & 2299 	&  0.149  &  0.224	 \\
\hline
\end{tabular}\,\,\,\,\,\,
\begin{tabular}[hp]{|c|cc|c|cc|}
\hline
\multicolumn{6}{|c|}{3D}  \\
\hline
$M_+$  & $M_p$ & $M_u$ 
& $DOF_c$ & $e_p$ (\%) & $e_u$ (\%) \\
\hline
\multicolumn{6}{|c|}{Case 1}  \\ % 1
\hline
0	& 1	& 3 	& 864 		&  4.369  &  5.477 \\
1	& 2	& 4 	& 1296 	&  2.750  &  4.045 \\
2	& 3	& 5 	& 1728 	&  1.707  &  3.530 \\
3	& 4	& 6 	& 2160 	&  1.348  &  3.052 \\
4	& 5	& 7 	& 2592 	&  0.906  &  2.691 \\
6	& 7	& 9 	& 3456 	&  0.467  &  2.329 \\
8	& 9	& 11 & 4320 	&  0.358  &  1.962 \\
\hline
\multicolumn{6}{|c|}{Case 2}  \\ % 14
\hline
0	& 1	& 3 	& 864 		&  6.000  &  3.653 \\
1	& 2	& 4 	& 1296 	&  4.940  &  3.224 \\
2	& 3	& 5 	& 1728 	&  3.585  &  2.731 \\
3	& 4	& 6 	& 2160 	&  2.514  &  2.496 \\
4	& 5	& 7 	& 2592 	&  1.805  &  2.327 \\
6	& 7	& 9 	& 3456 	&  1.042  &  2.052 \\
8	& 9	& 11 & 4320 	&  0.830  &  1.694 \\
\hline
\multicolumn{6}{|c|}{Case 3}  \\ %7
\hline
0	& 1	& 3 	& 864 		&  3.430  &  4.520 \\
1	& 2	& 4 	& 1296 	&  2.705  &  3.991 \\
2	& 3	& 5 	& 1728 	&  1.613  &  3.197 \\
3	& 4	& 6 	& 2160 	&  1.142  &  2.659 \\
4	& 5	& 7 	& 2592 	&  0.812  &  2.544 \\
6	& 7	& 9 	& 3456 	&  0.554  &  2.356 \\
8	& 9	& 11 & 4320 	&  0.454  &  2.062 \\
\hline
\end{tabular}
\end{center}
\caption{Numerical results for random permeability $k$ and elastic modulus $E$. 
Case 1, 2 and 3 (from left to right). 
Type 2 multiscale basis functions, $M_p = 1 + M_+$, $M_u = d + M_+$ and $DOF_c = (M_p + M_u) \cdot N_c$. 
Left:  two - dimensional problem, $d = 2$ ($N_c =121$ and $DOF_f = 30603$). 
Right: three  - dimensional problem, $d = 3$ ($N_c = 216$ and $DOF_f = 37044$) }
\label{tab:err-t2}
\end{table}

%\begin{tabular}[hp]{|c|cc|c|ccc|}
%\hline
% & & & & \multicolumn{3}{|c|}{Time (sec)} \\
%$M_+$  & $M_p$  & $M_u$  & $DOF_c$ & $R$ & $A_c$ & solve \\
%\hline 
%\multicolumn{3}{|c|}{fine grid}		& 30603	&  -	&  -	& 6.211 \\
%\hline 
%0	& 1	& 2 	& 363 		&  1.576	&  0.077 	& 0.350  \\
%1	& 2	& 3 	& 605 		&  2.183	&  0.170	& 0.587 \\
%2	& 3	& 4 	& 847 		&  2.824	&  0.278 	& 0.833 \\
%3	& 4	& 5 	& 1089 	&  3.340	&  0.412 	& 1.171 \\
%4	& 5	& 6 	& 1331 	&  4.197	&  0.561 	& 1.611 \\
%6	& 7	& 8 	& 1815 	&  5.478 	& 0.918  	& 2.767 \\
%8	& 9	& 10 & 2299 	&  6.448	& 1.371  	& 4.569 \\
%\hline
%\end{tabular}\,\,
%\begin{tabular}[hp]{|c|cc|c|ccc|}
%\hline
% & & & & \multicolumn{3}{|c|}{Time (sec)} \\
%$M_+$  & $M_p$  & $M_u$  & $DOF_c$ & $R$ & $A_c$ & solve \\
%\hline 
%\multicolumn{3}{|c|}{fine grid}	 & 37044	&  -	 &  -	 & 158.745 \\
%\hline 
%0	& 1	& 3 	& 864 		&  5.797	&  1.140 	& 2.525  \\
%1	& 2	& 4 	& 1296 	&  7.378	&  2.364	& 4.940 \\
%2	& 3	& 5 	& 1728 	&  9.159	&  3.587 	& 9.633 \\
%3	& 4	& 6 	& 2160 	&  10.907	&  5.282 	& 16.285 \\
%4	& 5	& 7 	& 2592 	&  12.761	&  7.338 	& 20.495 \\
%6	& 7	& 9 	& 3456 	&  16.154	& 12.217  & 38.951 \\
%8	& 9	& 11 & 4320 	&  19.713	& 18.442  & 69.322 \\
%\hline
%\end{tabular}

\begin{table}[h!]
\begin{center}
\begin{tabular}[hp]{|c|cc|c|c|}
\hline
\multicolumn{5}{|c|}{2D}  \\
\hline
$M_+$  & $M_p$  & $M_u$  & $DOF_c$ & Time (sec) \\
\hline 
\multicolumn{3}{|c|}{fine grid}		& 30603	& 6.211 \\
\hline 
0	& 1	& 2 	& 363 		& 0.350  \\
1	& 2	& 3 	& 605 		& 0.587 \\
2	& 3	& 4 	& 847 		& 0.833 \\
3	& 4	& 5 	& 1089 	& 1.171 \\
4	& 5	& 6 	& 1331 	& 1.611 \\
6	& 7	& 8 	& 1815 	& 2.767 \\
8	& 9	& 10 & 2299 	& 4.569 \\
\hline
\end{tabular}\,\,\,\,
\begin{tabular}[hp]{|c|cc|c|c|}
\hline
\multicolumn{5}{|c|}{3D}  \\
\hline
$M_+$  & $M_p$  & $M_u$  & $DOF_c$ & Time (sec) \\
\hline 
\multicolumn{3}{|c|}{fine grid}	 & 37044	& 158.745 \\
\hline 
0	& 1	& 3 	& 864 		& 2.525  \\
1	& 2	& 4 	& 1296 	& 4.940 \\
2	& 3	& 5 	& 1728 	& 9.633 \\
3	& 4	& 6 	& 2160 	& 16.285 \\
4	& 5	& 7 	& 2592 	& 20.495 \\
6	& 7	& 9 	& 3456 	& 38.951 \\
8	& 9	& 11 & 4320 	& 69.322 \\
\hline
\end{tabular}
\end{center}
\caption{Solution time for Case 3. Type 2 multiscale basis functions. 
Left:  two - dimensional problem, $d = 2$. 
Right: three  - dimensional problem, $d = 3$
}
\label{tab:err-time}
\end{table}

We consider two and three-dimensional model problems in domain $\Omega = [0,1]^d$ with  $d = 2,3$:
\begin{itemize}
\item Two - dimensional problem (2D). Coarse grid contains 121 nodes ($10 \times 10$) and the fine grid has 10201 nodes. %($100 \times 100$).  
\item Three - dimensional problem (3D). Coarse grid contains of 216 nodes ($5 \times 5 \times 5$) and fine grid has 9261 nodes. %($20 \times 20 \times 20$).  
\end{itemize}
The coarse grid and fine grid are presented in Figure \ref{schematic}.  
%We note that the MCMC method is independent of fine mesh size. %, and we used that can handle.

We perform simulations for $T_{max} = 0.001$ with 20 time steps. As an initial condition, we set $p_0 = 0$ and impose the following boundary conditions:
\[
u_x = 0, \quad \sigma_y = 0 , \quad \sigma_z = 0,
\quad x = 0, 
\]\[
\sigma_x = 0, \quad u_y = 0, \quad \sigma_z = 0,
\quad y = 0,
\]\[
\sigma_x = 0, \quad \sigma_y = 0, \quad u_z = 0,
\quad z = 0,
\]\[
\sigma \cdot n = 0, \quad x,y,z=1,
\]
and
\[
q \cdot n = \gamma (p - p_1), \quad y = 1,
\]\[ 
q \cdot n = 0, \quad y \neq 1,
\]
with $\gamma = 10^4$ and  $p_1 = 1.0$.

The high dimensional heterogeneity is represented using Karhunen-Lo{\'e}ve expansion. Random heterogeneous porosity fields are generated using $L = 200$ basis functions ($Y_L(x,\theta) $) with exponential covariance ($l_x = l_y = l_z =0.2$ and $\sigma_R^2 = 2$)
\[
\phi(x,\theta) = \phi(Y(x,\theta)),
\]
with normalization, that give $\phi \in [0.05, 0.2]$.

Heterogeneous permeability and elasticity modulus fields are given by
\[
k(x,\theta) = \exp\left(
a \, \phi(x,\theta) 
\right),  
\]\[
E(x,\theta) = 
b \left( 
\frac{1.0 - \phi(x,\theta) }{\phi(x,\theta) }
\right)^m,
\]
with $a = 40$, $b = 0.1$ and $m = 1.5$  \cite{yang2013simulator}. 
The Biot modulus is $M = 1.0$, fluid viscosity is $\nu = 1$,  fluid-solid coupling constant is $\alpha = 0.1$, the Poisson's ratio is $\eta = 0.3$. Permeability field and elastic modulus are shown in Figure \ref{fig:2d-kE} for 2D problems and in Figure \ref{fig:3d-kE} for 3D problems. We consider three test cases (Case 1, 2 and 3 are depicted from left to right). All permeabilities and elastic modulus are depicted in the log scale.

The reference solution computed by a standard finite element method with linear basis functions for pressure and displacements on the fine grid. The reference pressure and the displacement fields at final time are presented in  Figures \ref{fig:2d-up} and  \ref{fig:3d-up}.

To compare a multiscale solution,  we calculate relative errors in $L^2$ norm in \%
\[
e_p =  \sqrt{ \frac{\int_{\Omega}  (p^{ms} - p, p^{ms} - p) \, dx}{\int_{\Omega} (p, p) \,dx} } \cdot 100 \%, 
\]\[
e_u = \sqrt{ \frac{\int_{\Omega}  (u^{ms} - u, u^{ms} - u) \,dx}{\int_{\Omega} (u, u) \,dx} }  \cdot 100 \%, 
\]
where $p^{ms}$ and $u^{ms}$ are multiscale solutions, $p$ and $u$ are reference solutions. 

For multiscale basis construction on the offline stage, we use a 10 random permeability and elastic fields. We present numerical results in Tables \ref{tab:err-t1} and \ref{tab:err-t2} for Type 1 and 2 multiscale basis functions.
$DOF_c$ and $DOF_f$ are degrees of freedom for multiscale and reference (fine grid) solutions. $M_p$ and $M_u$ are the number of the multiscale basis functions for pressure and displacements, respectively. Varying the basis functions in both pressure and displacement multiscale spaces we presented the errors at the final times.
We note that the first basis for pressure is standard linear basis function because the first eigenvalue of the local spectral problem is constant.
Similarly for displacements, the first $d$ basis functions are  standard linear basis functions. $M_+$ is used to denote the number of the additional spectral basis functions calculated using algorithm presented in Section 3.
Therefore, $M_p = 1 + M_+$, $M_u = d + M_+$ and $DOF_c = (M_p + M_u) \cdot N_c$, where $N_c$ is the number of coarse grid nodes. For two - dimensional problem, we have $N_c =121$ and $DOF_f = 30603$. For three - dimensional problem, we have $N_c = 216$ and $DOF_f = 37044$.
We can obtain good multiscale solution when we take sufficient number of multiscale basis functions for pressure and for displacements.
For $M_+ = 2$, we have near $1-3\%$ of errors in two - dimensional and three-dimensional problems for all cases.
In Case 3 of heterogeneous permeability field, we have $0.7$ \% of pressure error and $2.0$ \% of displacement error in two-dimensional problems with 3 multiscale basis functions for pressure and 4 multiscale basis functions for displacements $M_+=2$. In this case, we reduce size of the system from $DOF_f = 30603$ to $DOF_c = 847$. For three - dimensional problem, we have $1.7$ \% of pressure error and $3.3$ \% of displacement error with 3 basis functions for pressure and 5 basis functions for displacements ($M_+=2$). We reduce size of the system from $DOF_f = 37044$ to $DOF_c = 1728$.

In Table \ref{tab:err-time}, we present solution time for the coarse grid and fine grid solvers for Case 3. Solution time of the fine grid solver is $6.2$ second for two - dimensional problem  ($DOF_f = 30603$) and $158.7$ second for the three-dimensional problem ($DOF_f = 37044$). When we use a multiscale method for the solution with $M_+ = 2$, we solve two - dimensional problem by $0.8$ seconds ($DOF_c = 847$) and three-dimensional problem by $9.6$ seconds ($DOF_c = 1728$). We see that a smaller number of basis functions give a coarse grid system with a smaller size and therefore solution time is faster. Here for the solution of the coarse and fine grid systems, we used a direct solver (default solver in FEniCS \cite{logg2012automated}). Note that, the solution time doesn't include a time of multiscale basis construction because they are constructed on the offline stage as precalculations.

\subsection{Preconditioned MCMC using GMsFEM}\label{ss2}

In this section, we present results for the Two-stage MCMC method.
MCMC simulations on the fine grid (single-stage) are generally very computationally expensive because each proposal requires solving a forward coupled poroelasticity problem over a large time interval. 
In presented preconditioned MCMC simulations, we use inexpensive computations in the first stage. 
Preconditioning procedure is performed using 
\begin{enumerate}
\item \textit{MS} multiscale solver based on GMSFEM,
\item \textit{ML} machine learning technique.
\end{enumerate}
Note that, the Preconditioning procedure should be inexpensive, but not necessarily very accurate. 
We start with a multiscale solver with a small number of multiscale basis functions. As we showed above, the size of the coarse grid system ($DOF_c$) depends on a number of multiscale basis functions. The multiscale basis functions are constructed only once on the offline stage, and we use them for all proposed permeability and elastic modulus without online recalculations. 

% t2
\begin{figure}[h!]
\centering
    \begin{subfigure}[b]{0.49\textwidth}
        \includegraphics[width=1.0\linewidth]{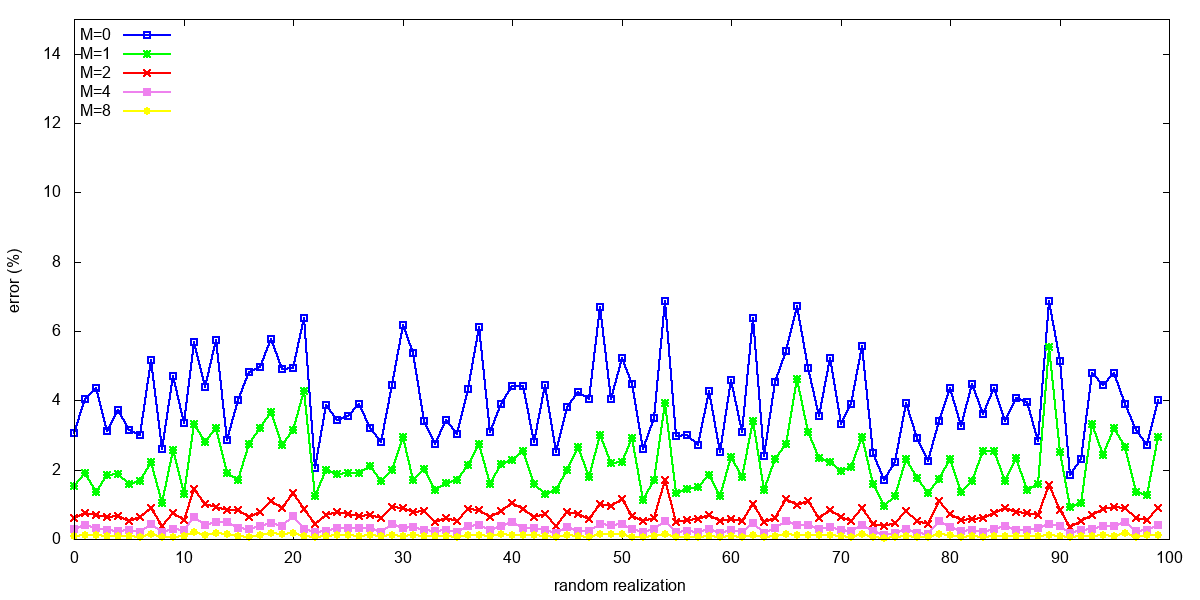}
        \caption{2D, $e_p$ (\%)}
    \end{subfigure}
    \begin{subfigure}[b]{0.49\textwidth}
        \includegraphics[width=1.0\linewidth]{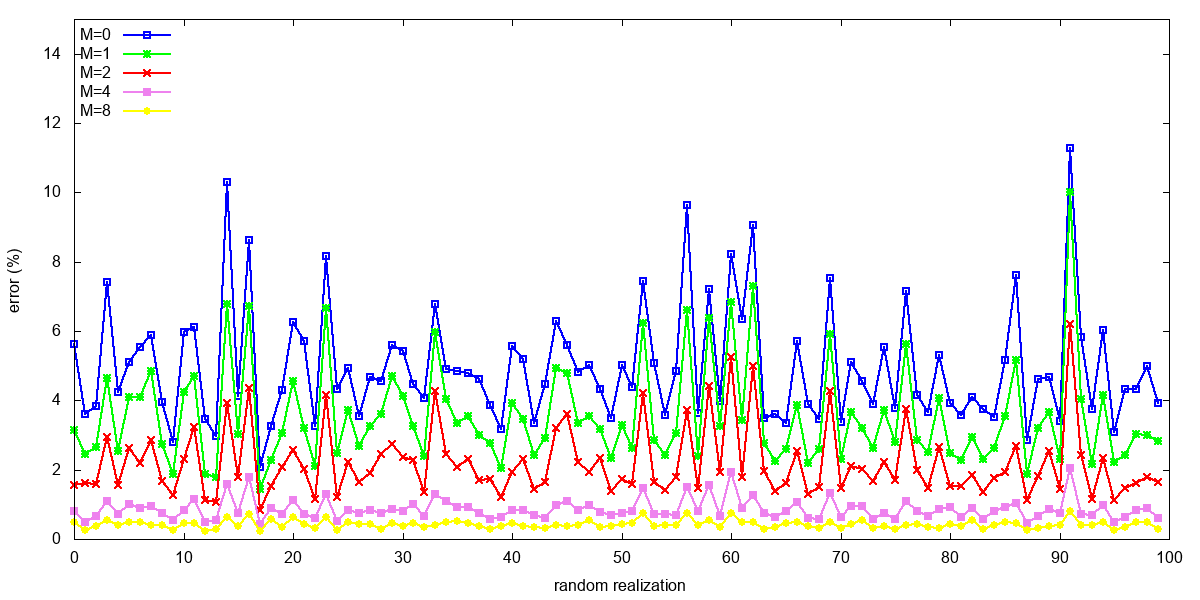}
        \caption{3D, $e_p$ (\%)}
    \end{subfigure}
    \\
    \vspace{10pt}
   \begin{subfigure}[b]{0.49\textwidth}
        \includegraphics[width=1.0\linewidth]{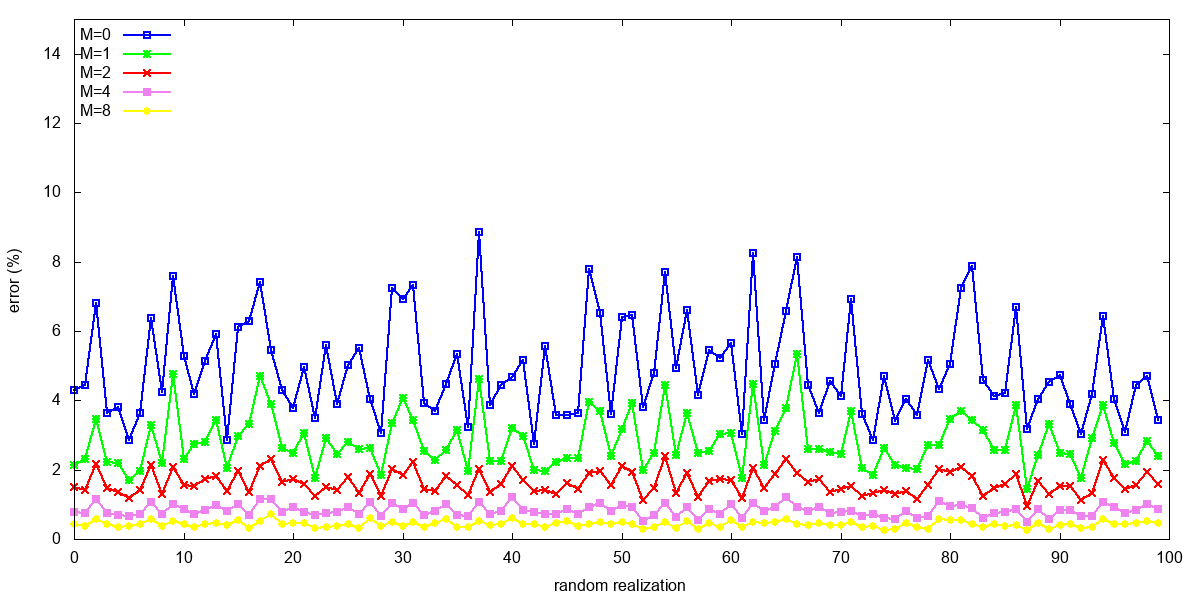}
        \caption{2D, $e_u$ (\%)}
    \end{subfigure}
    \begin{subfigure}[b]{0.49\textwidth}
        \includegraphics[width=1.0\linewidth]{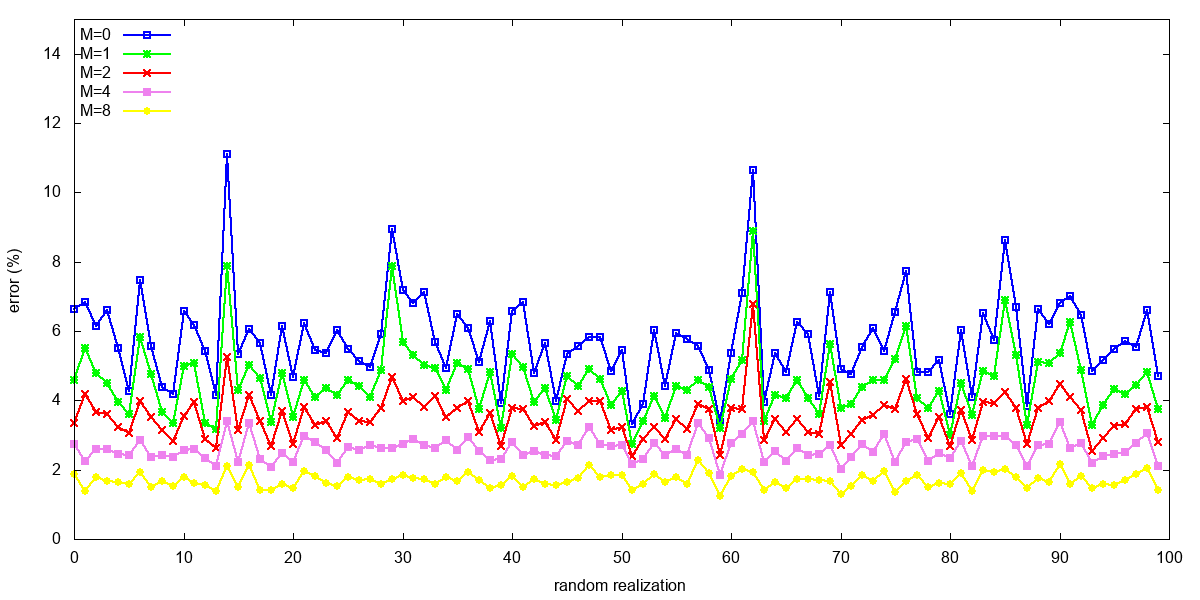}
        \caption{3D, $e_u$ (\%)}
    \end{subfigure}
\caption{Multiscale method errors for 100 random realizations of heterogeneous filed for $M_+=0, 1, 2, 4$ and $8$.
First row: pressure relative errors in \%.
Second row: displacement relative errors in \%.
Left: two - dimensional problem.
Right: three - dimensional problem. }
\label{fig:err-100}
\end{figure}

\begin{figure}[h!]
\centering
    \begin{subfigure}[b]{0.4\textwidth}
        \includegraphics[width=1.0\linewidth]{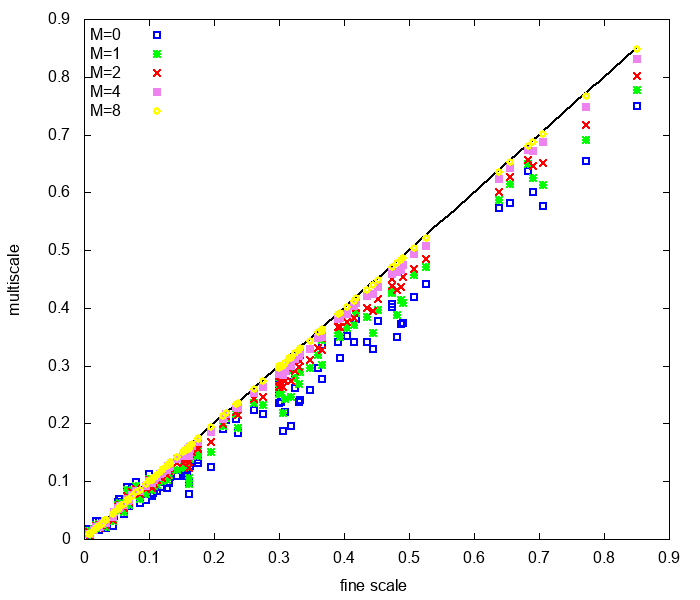}
        \caption{2D}
    \end{subfigure}
    \begin{subfigure}[b]{0.4\textwidth}
        \includegraphics[width=1.0\linewidth]{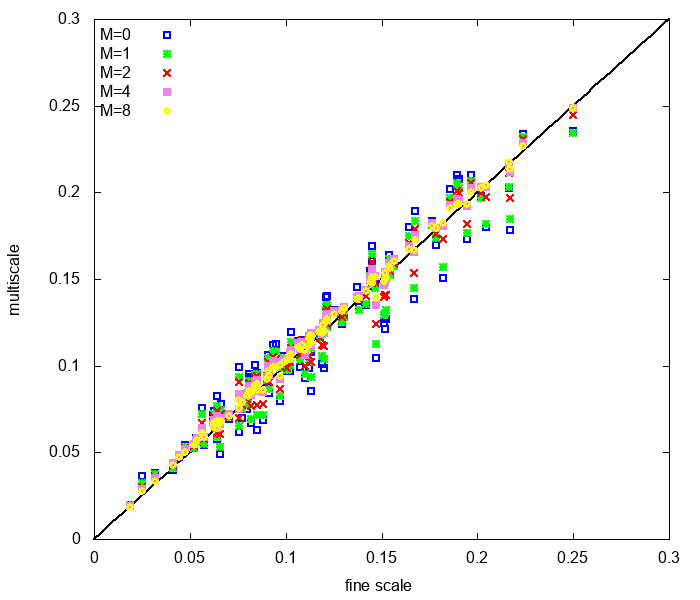}
        \caption{3D}
    \end{subfigure}
\caption{Cross-plot  between ${E(\theta)}$ (fine scale) and ${E^*(\theta)}$ (multiscale solver, \textit{MS}) for $M_+=0, 1, 2, 4$ and $8$.  For given observation data. 
(a) two - dimensional problem. 
(b) three  - dimensional problem. }
\label{fig:F-100}
\end{figure}

\begin{figure}[h!]
\centering
\includegraphics[width=0.19\textwidth]{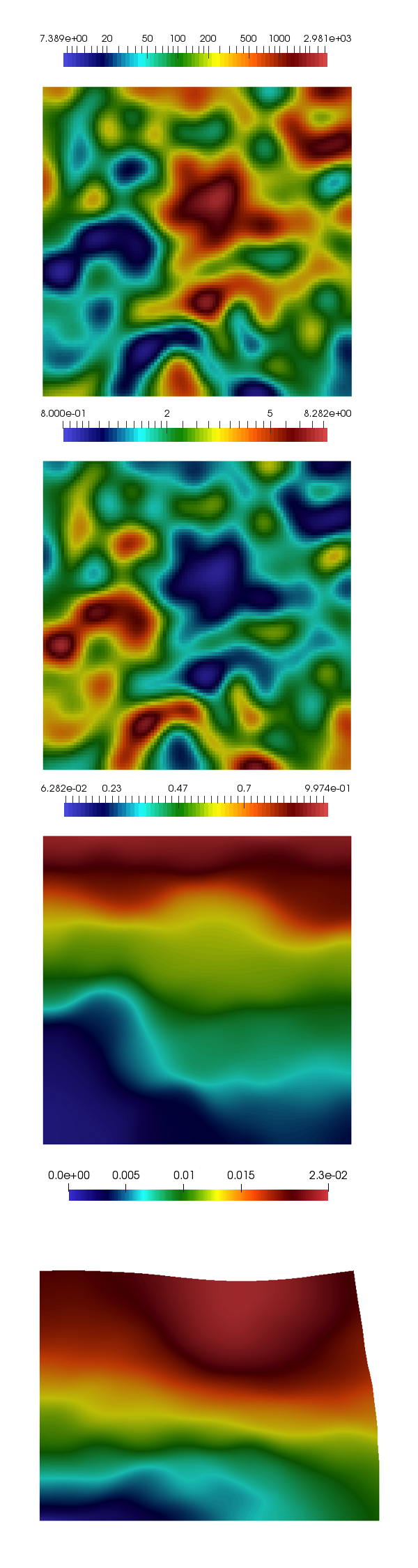} 
\includegraphics[width=0.19\textwidth]{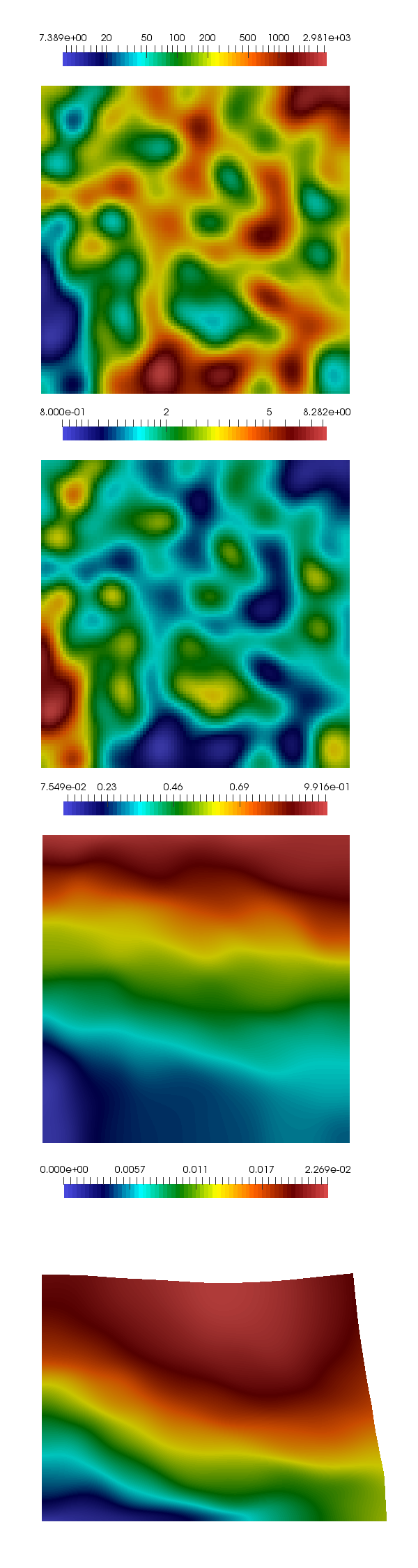} 
\includegraphics[width=0.19\textwidth]{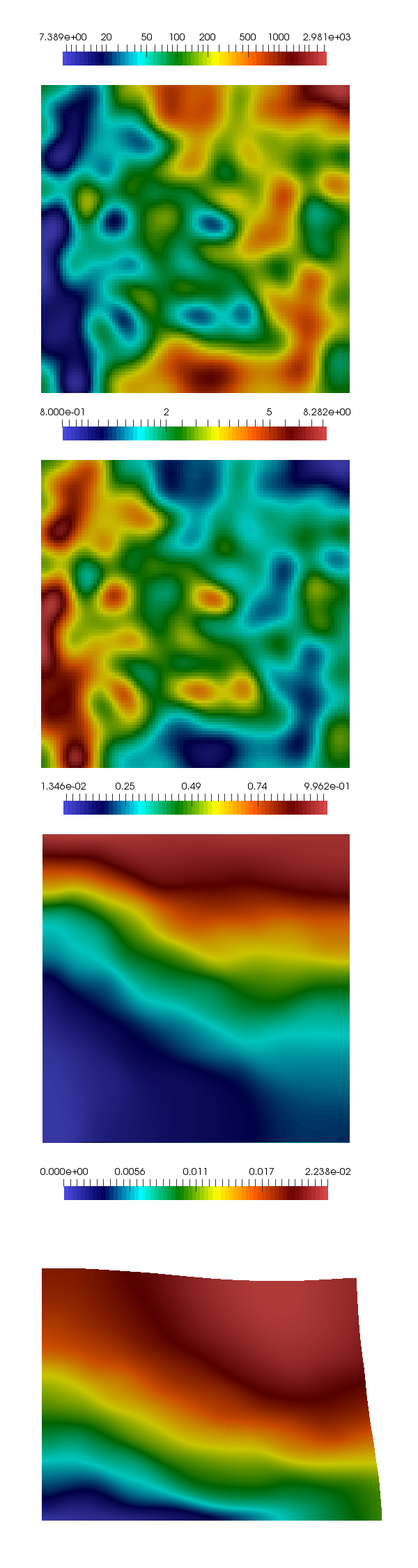} 
\includegraphics[width=0.19\textwidth]{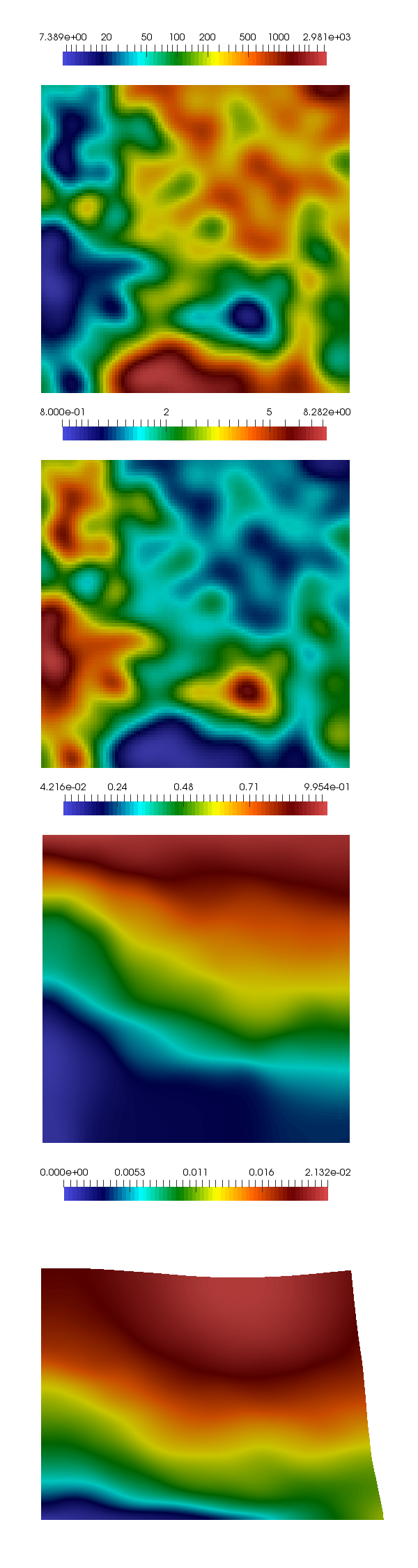} 
\includegraphics[width=0.19\textwidth]{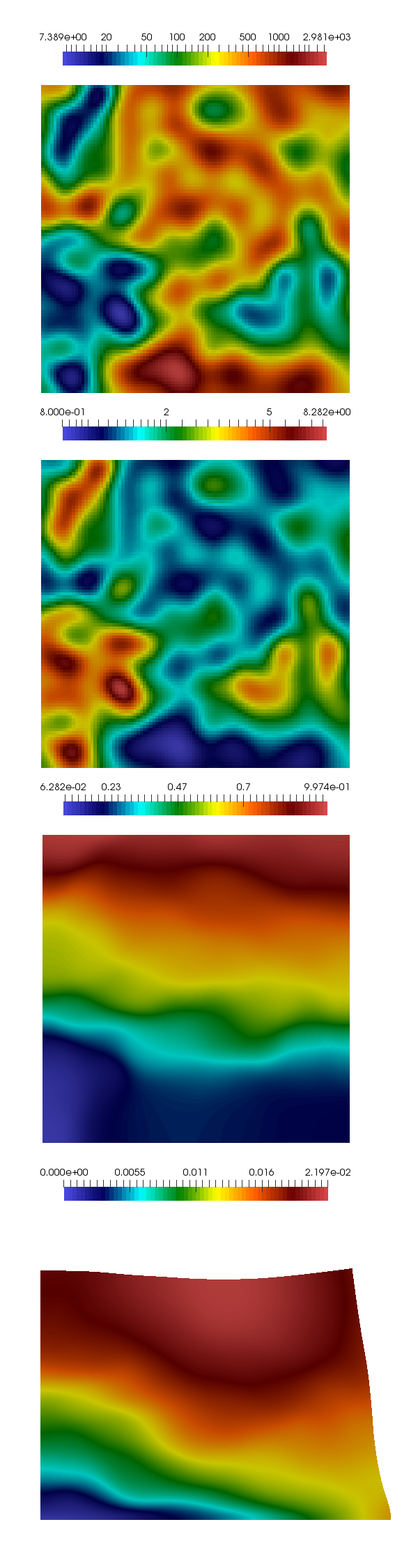} 
\caption{Reference solution and solutions for accepted random fields for two - dimensional problem (from left to right).
Two-stage MCMC with $\sigma_f = 0.02$, $\beta = 2$ and $M_+ = 2$.
Random walk samples with $\delta = 0.5$.
31 accepted and 69 passed the first stage.
First row: permeability, $k$.
Second row: elastic modulus, $E$.
Third row: pressure, $p$.
Fourth row: displacement, $u_m$.
 }
\label{fig:passed-2d}
\end{figure}

\begin{figure}[h!]
\centering
\includegraphics[width=0.19\textwidth]{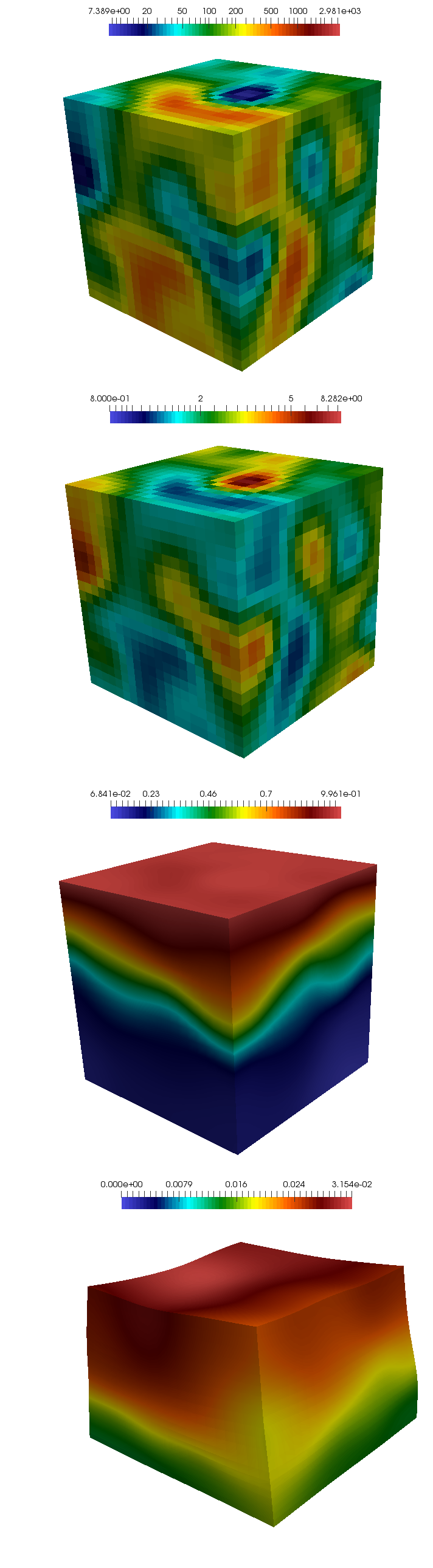} 
\includegraphics[width=0.19\textwidth]{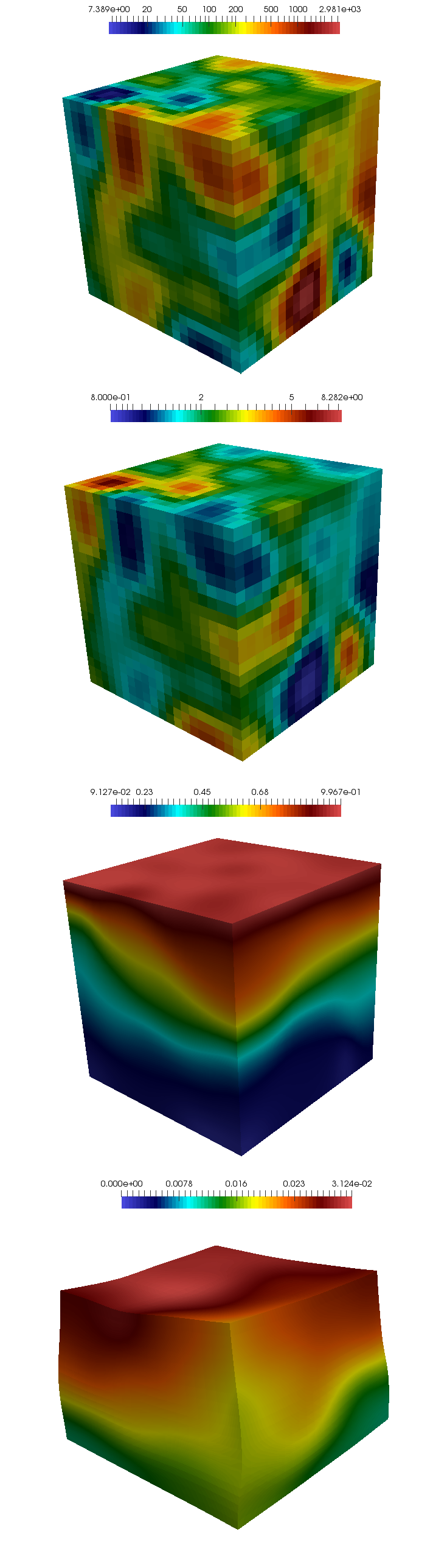} 
\includegraphics[width=0.19\textwidth]{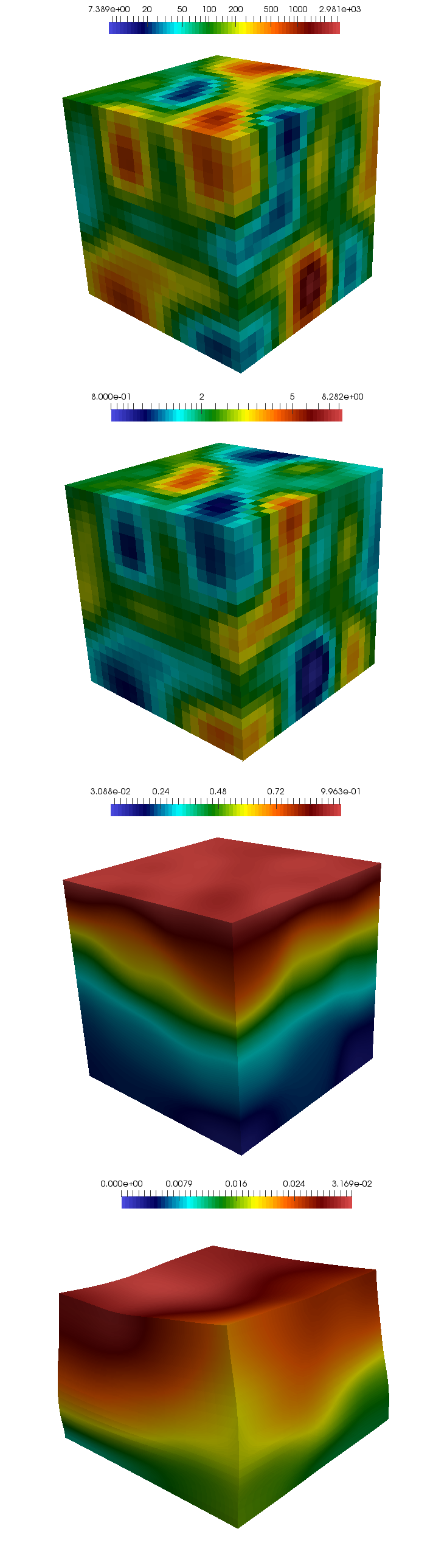}
\includegraphics[width=0.19\textwidth]{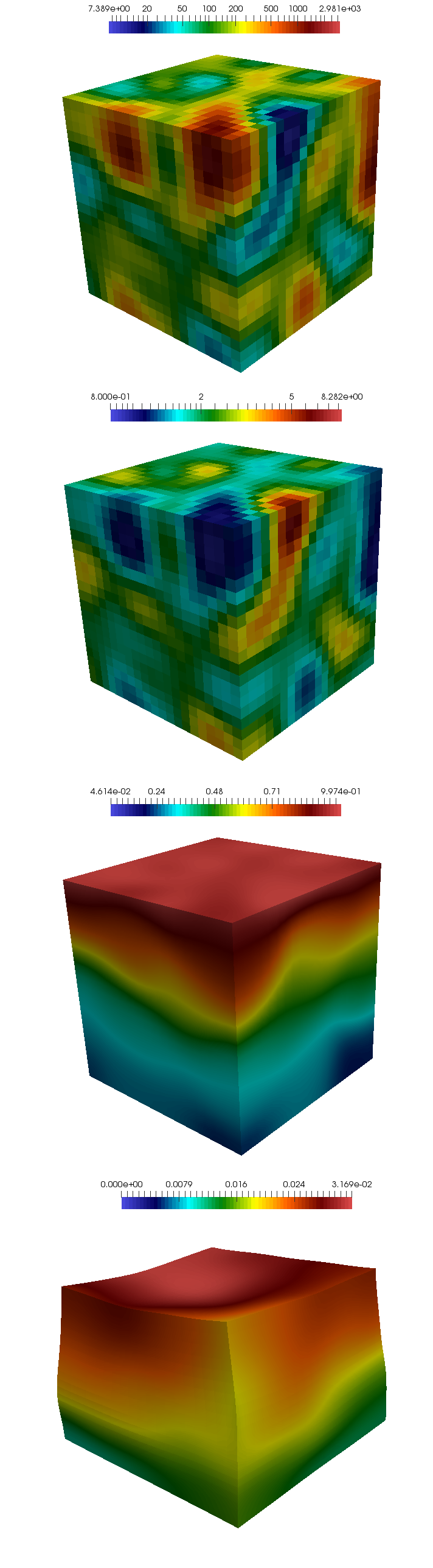} 
\includegraphics[width=0.19\textwidth]{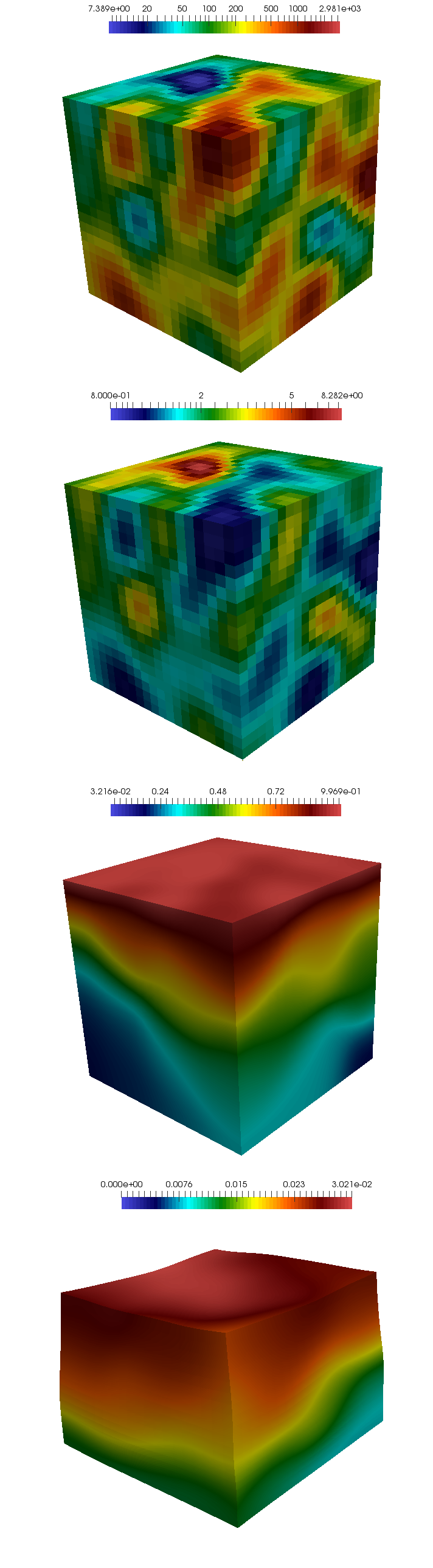} 
\caption{Reference solution and solutions for accepted random fields for three - dimensional problem (from left to right).  
Two-stage MCMC with $\sigma_f = 0.02$,  $\beta = 2$ and $M_+ = 2$. 
Random walk samples with $\delta = 0.5$. 
29 accepted and 141 passed the first stage.
First row: permeability, $k$. 
Second row: elastic modulus, $E$. 
Third row: pressure, $p$. 
Fourth row: displacement, $u_m$. 
 }
\label{fig:passed-3d}
\end{figure}

We present results for the proposed multiscale solver (\textit{MS}) for 100 random realizations of the heterogeneous fields, $\theta_r$.  
Relative errors for pressure and displacements are presented in Figure \ref{fig:err-100}.   
In Figure \ref{fig:F-100}, we plot $E(\theta)$ against $E^*(\theta)$ to demonstrate correlation between ${E^*(\theta)}$ and ${E(\theta)}$ for 100 realizations of random fields.  
Here we use a relative difference between observation data and proposed solution for given  $\theta$
\begin{equation}
\label{eq:E}
E^*(\theta) = \frac{||F^*(\theta) - F_{obs}||^2}{ ||F_{obs}||^2}, \quad 
E(\theta) = \frac{||F(\theta) - F_{obs}||^2}{ ||F_{obs}||^2}, 
\end{equation}
where $F =u(x, t)$, $F^* = u^{ms}(x, t)$ where $t = T_{max}$ (final time)  and $x \in  \partial \Omega_{surf}$ ($ \partial \Omega_{surf}$ is the top boundary of the domain).  As observation data, we take solution that presented in Figure \ref{fig:passed-2d} and \ref{fig:passed-3d}. 
We observe that if we take a sufficient number of multiscale basis, we can obtain a good correlation between these two characteristics. 
From this perspective we will use in our calculations $M_+ = 2$ in MCMC algorithm which have $DOF_c = 847$ for 2D problem ($M_p = 3$ and $M_u = 4$) and $DOF_c = 1728$ for 3D problem ($M_p = 3$ and $M_u = 5$).  Fine grid system size is $DOF_f = 30603$ for  two-dimensional problem and $DOF_f = 37044$ for  three-dimensional problem.
The random heterogeneity is represented using Karhunen-Lo{\'e}ve expansion. In MCMC iterations, the heterogeneous porosity field is generated using $L = 200$ basis functions ($Y_L(x,\theta) $) with exponential covariance ($l_x = l_y = l_z =0.2$ and $\sigma_R^2 = 2$).  Form of $\phi$, $k$ and $E$ are presented in \ref{porosity}, \ref{logperm} and \ref{E.stoch} with $a = 40$, $b = 0.1$.  
In numerical simulation, we  pick a some realization of the random field $\theta_{ref}$ and use corresponded solution as observation data (see first columns of Figures \ref{fig:passed-2d} and \ref{fig:passed-3d}).

% t34
\begin{figure}[h!]
\centering
 	\begin{subfigure}[b]{0.32\textwidth}
        \includegraphics[width=1.0\linewidth]{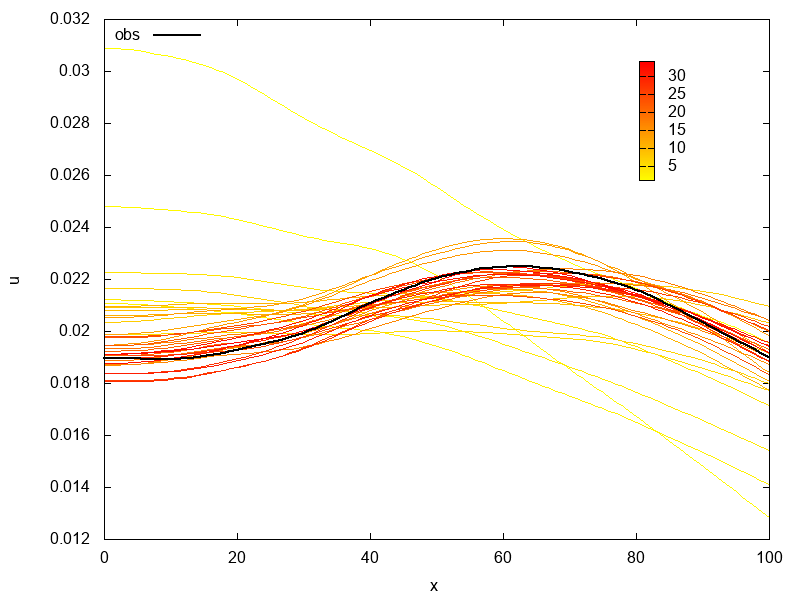}\\
 		\includegraphics[width=1.0\linewidth]{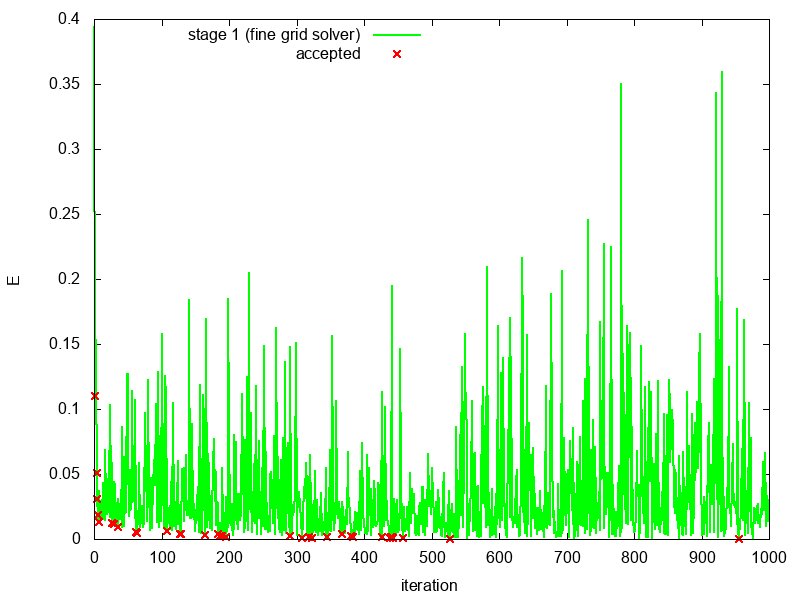}
        \caption{Single-stage MCMC}
    \end{subfigure} 
    \begin{subfigure}[b]{0.32\textwidth}
        \includegraphics[width=1.0\linewidth]{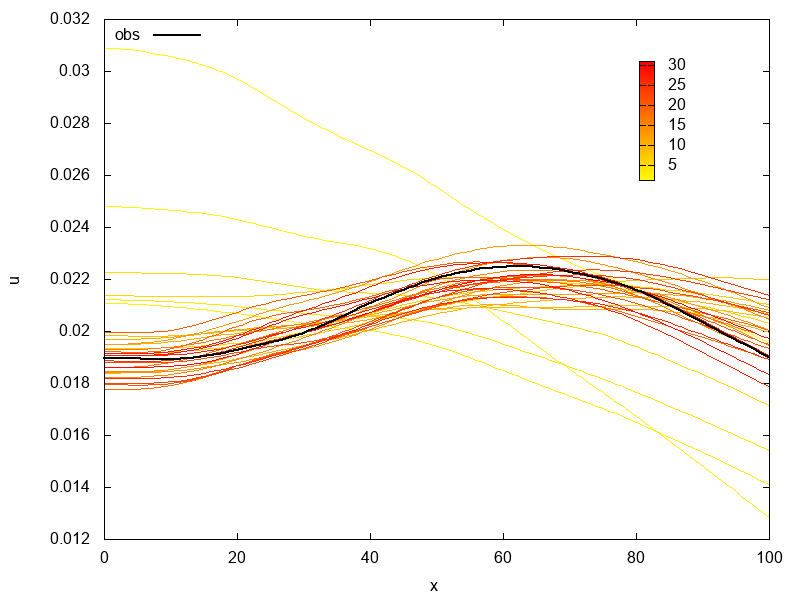}\\
        \includegraphics[width=1.0\linewidth]{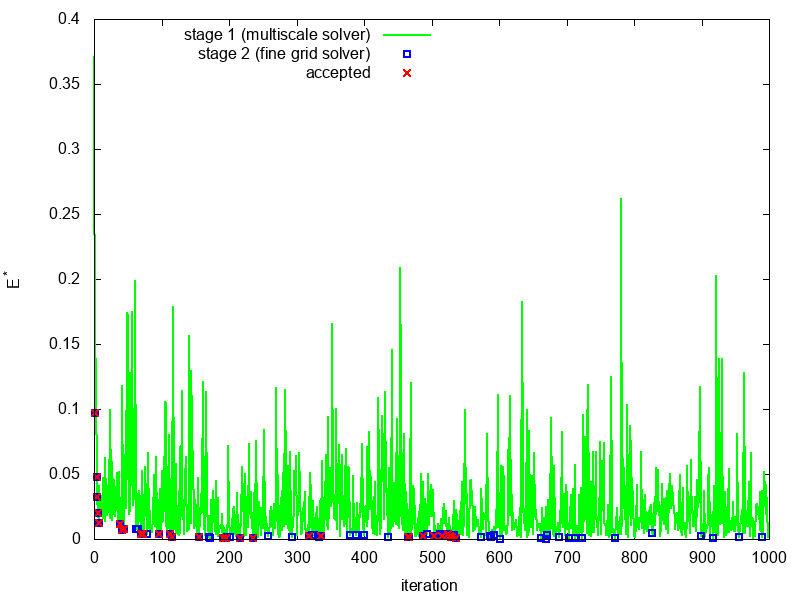}
        \caption{Two-stage MCMC, $M_+ = 2$}
    \end{subfigure}
    \begin{subfigure}[b]{0.32\textwidth}
        \includegraphics[width=1.0\linewidth]{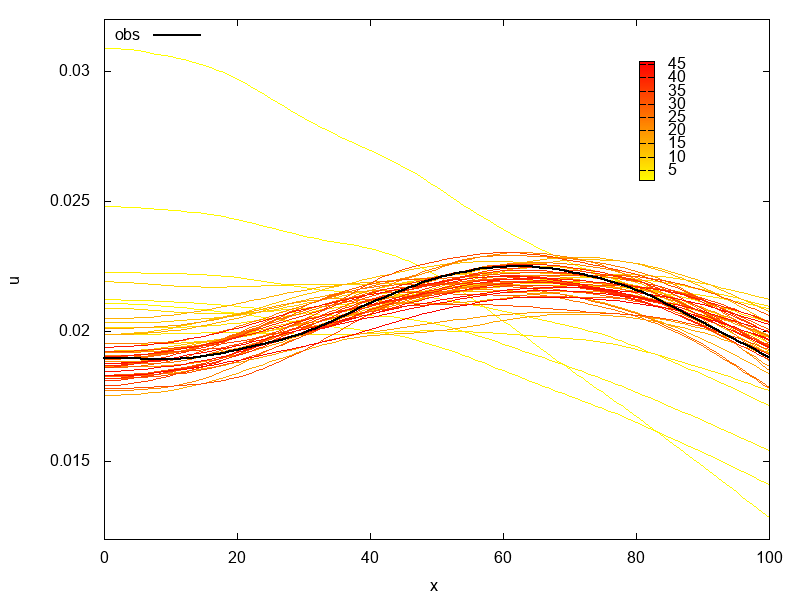}\\
        \includegraphics[width=1.0\linewidth]{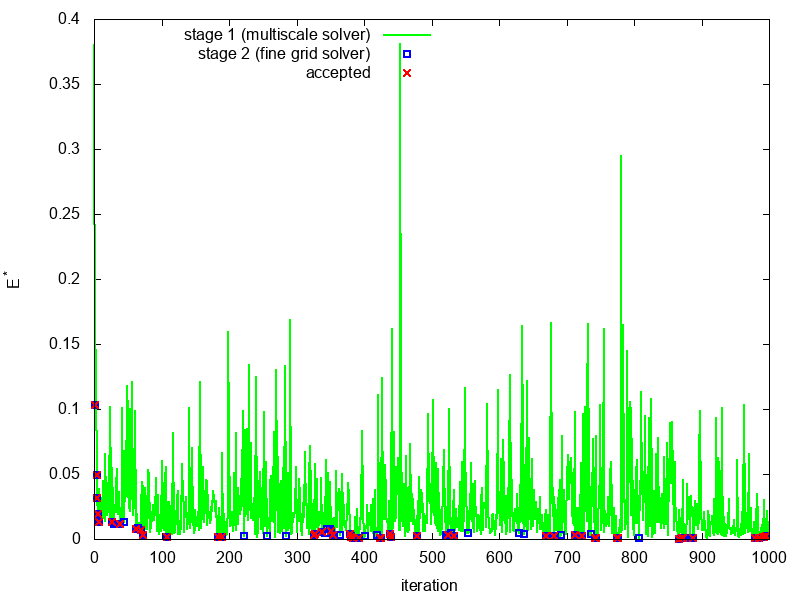}
        \caption{Two-stage MCMC, $M_+ = 4$}
    \end{subfigure}
\caption{Two - dimensional problem with $\sigma_f = 0.02$.
Random walk samples with $\delta = 0.5$.
First row: accepted $u_{obs}$ on surface boundary. 
Second row:  $E^*(\theta)$ (green color) and $E(\theta)$(blue color) in each MCMC iteration. 
(a) Single-stage MCMC (34 accepted).
(b) Two-stage MCMC with $\beta = 2$ and $M_+ = 2$ (31 accepted and 69 passed first stage).
(b) Two-stage MCMC with $\beta = 2$ and $M_+ = 4$ (46 accepted and 68 passed first stage). }
\label{fig:acc-mcmc-M}
\end{figure} 

% 3d
\begin{figure}[h!]
\centering
\includegraphics[width=0.49\linewidth]{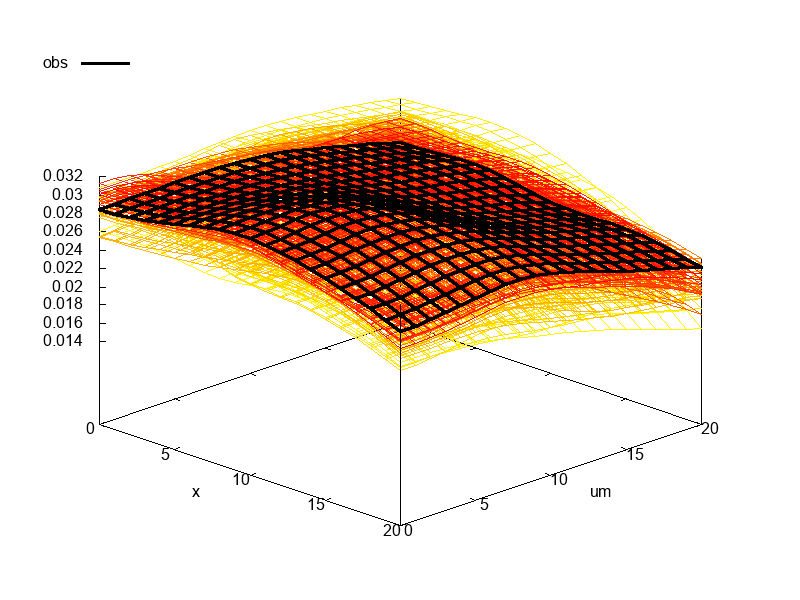}
\includegraphics[width=0.49\linewidth]{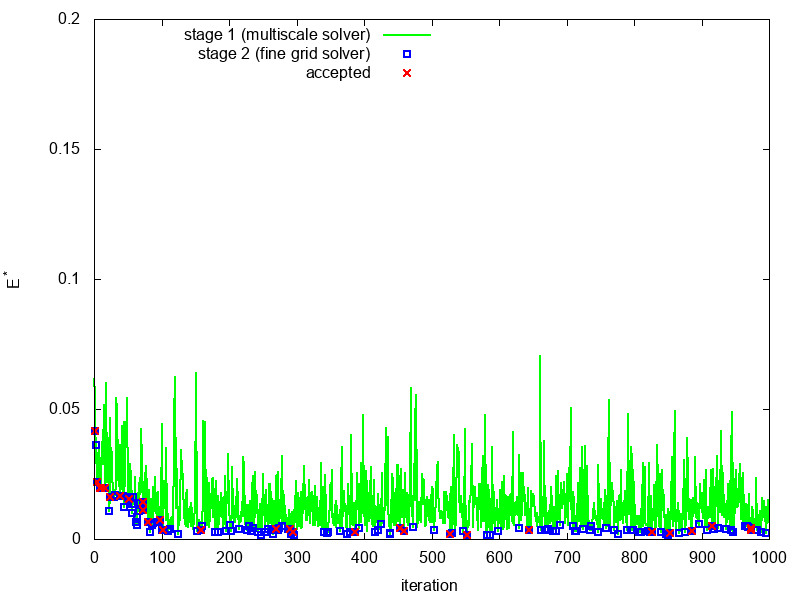}
\caption{Three - dimensional problem with $M_+ = 2$.
Two-stage MCMC with $\sigma_f = 0.02$ and $\beta = 2$. 
Random walk samples with $\delta = 0.5$. 
First row: accepted $u_{obs}$ on surface boundary. 
Second row:  $E^*(\theta)$ (green color) and $E(\theta)$(blue color) in each MCMC iteration. 
29 accepted and 141 passed the first stage.
}
\label{fig:acc-mcmc-3d}
\end{figure}

% 2d param
\begin{figure}[h!]
\centering
 	\begin{subfigure}[b]{0.32\textwidth}
        \includegraphics[width=1.0\linewidth]{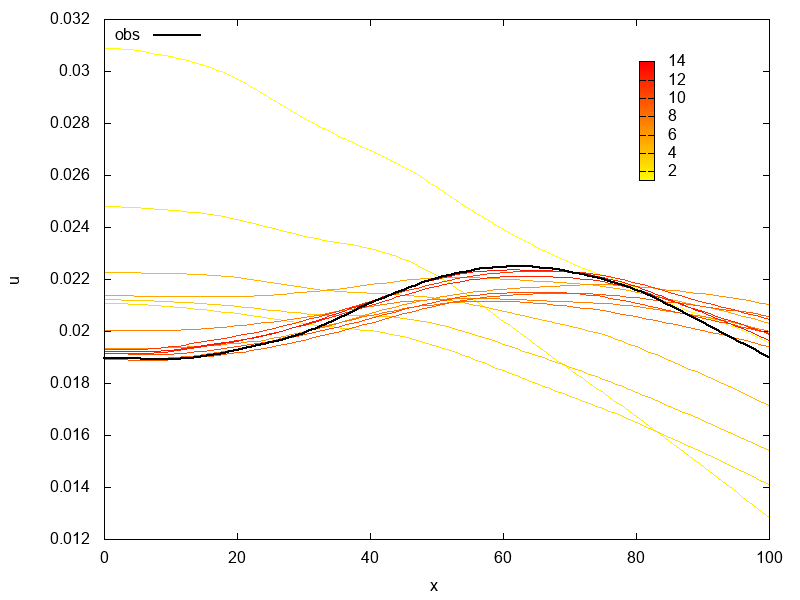}\\
 		\includegraphics[width=1.0\linewidth]{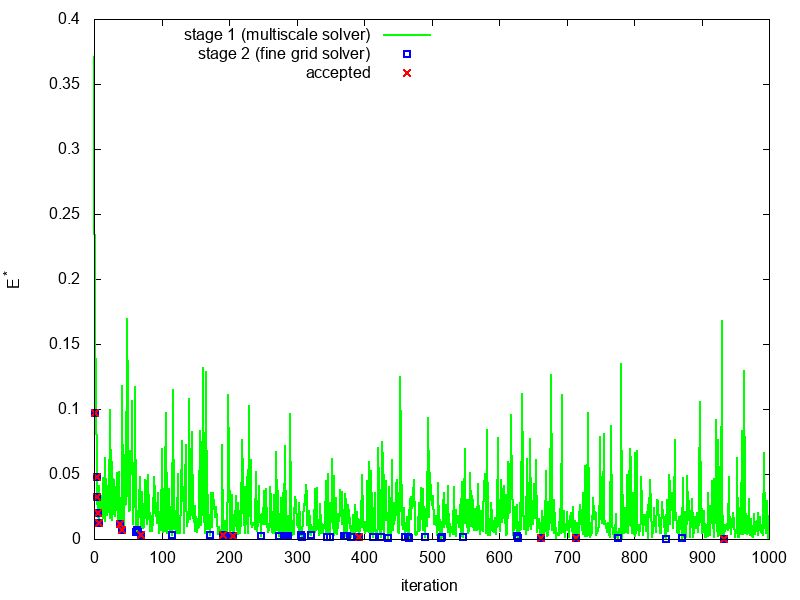}
        \caption{Two-stage MCMC, $\sigma_f = 0.01$}
    \end{subfigure} 
    \begin{subfigure}[b]{0.32\textwidth}
        \includegraphics[width=1.0\linewidth]{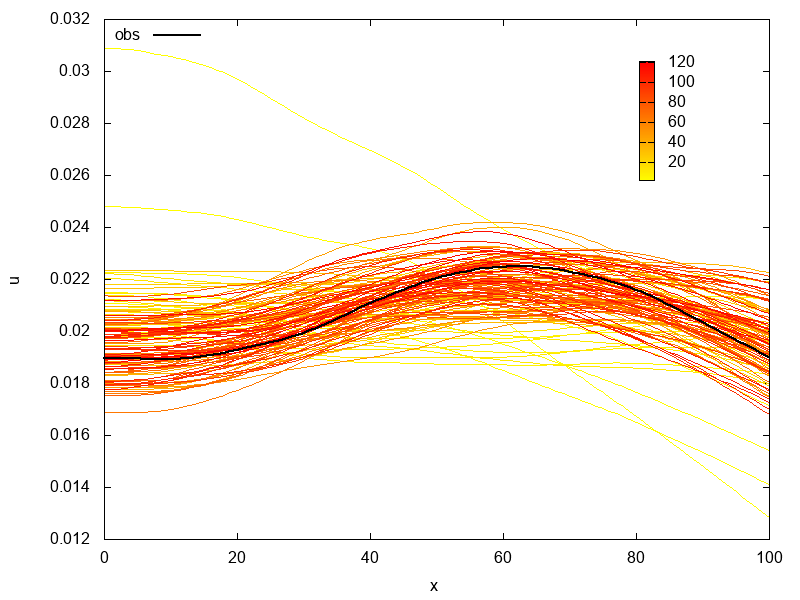}\\
 		\includegraphics[width=1.0\linewidth]{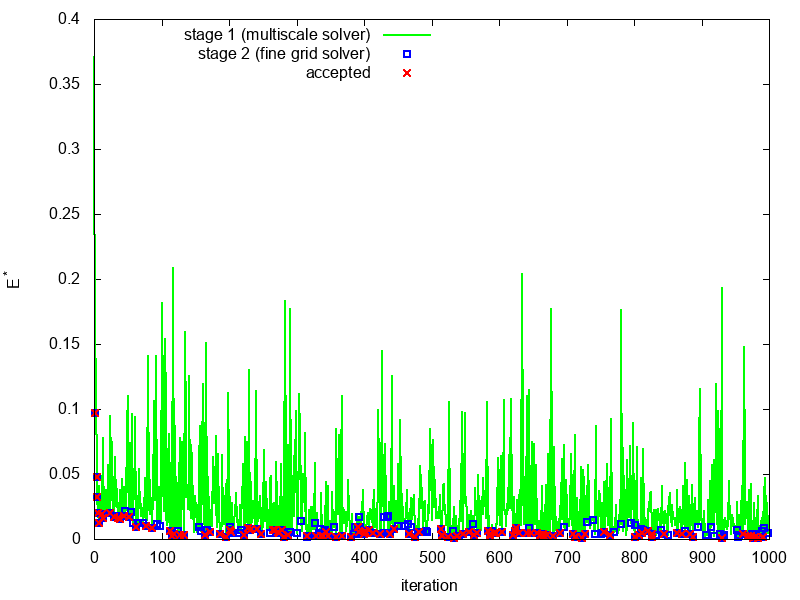}
        \caption{Two-stage MCMC, $\sigma_f = 0.04$}
    \end{subfigure} 
    \begin{subfigure}[b]{0.32\textwidth}
        \includegraphics[width=1.0\linewidth]{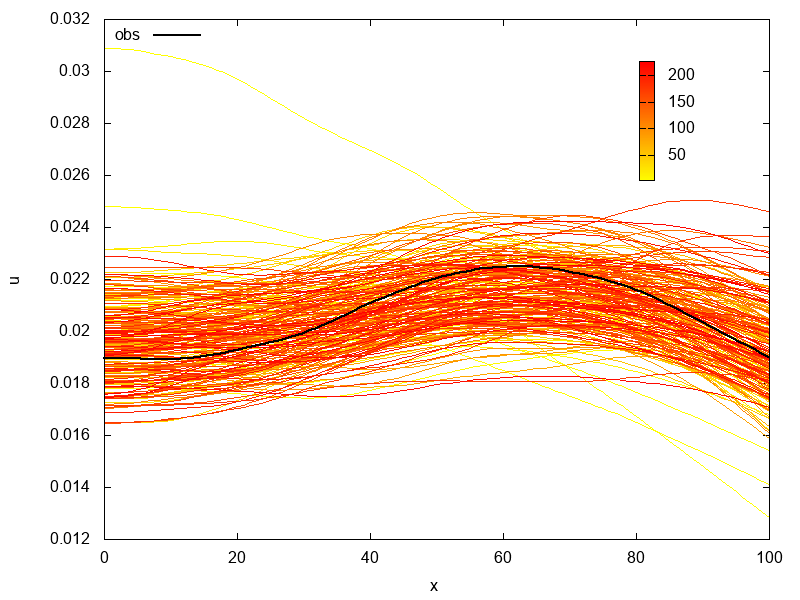}\\
        \includegraphics[width=1.0\linewidth]{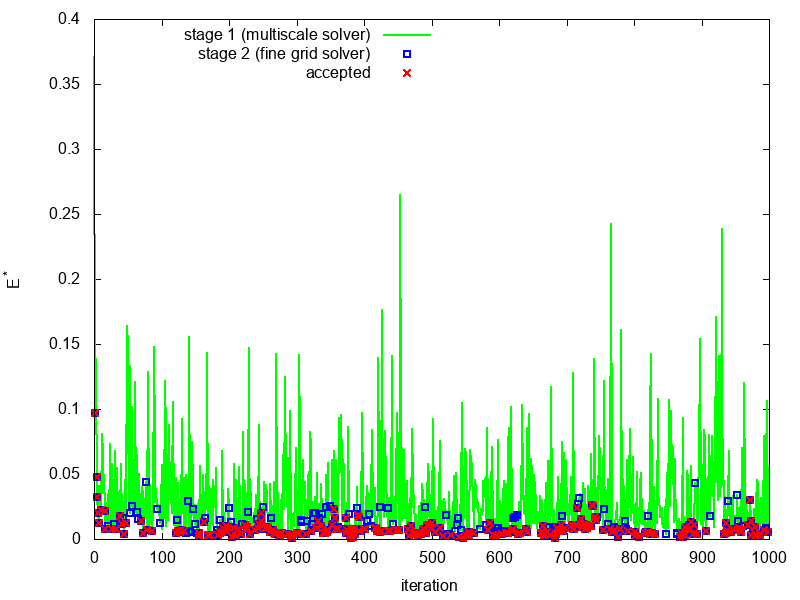}
        \caption{Two-stage MCMC, $\sigma_f = 0.06$}
    \end{subfigure}
\caption{Two - dimensional problem with $M_+ = 2$.
Two-stage MCMC with $\beta = 2$ for different $\sigma_f = 0.01, 0.04, 0.06$.
Random walk samples with $\delta = 0.5$. 
First row: accepted $u_{obs}$ on surface boundary. 
Second row:  $E^*(\theta)$ (green color) and $E(\theta)$(blue color) in each MCMC iteration. 
(a) $\sigma_f = 0.01$ (14 accepted and 50 passed first stage).
(b) $\sigma_f = 0.04$ (121 accepted and 212 passed first stage).
(c) $\sigma_f = 0.06$ (226 accepted and 336 passed first stage). }
\label{fig:acc-mcmc-sigma}
\end{figure}

% 2d hist
\begin{figure}[h!]
\centering
	\begin{subfigure}[b]{0.32\textwidth}
        \includegraphics[width=1.0\linewidth]{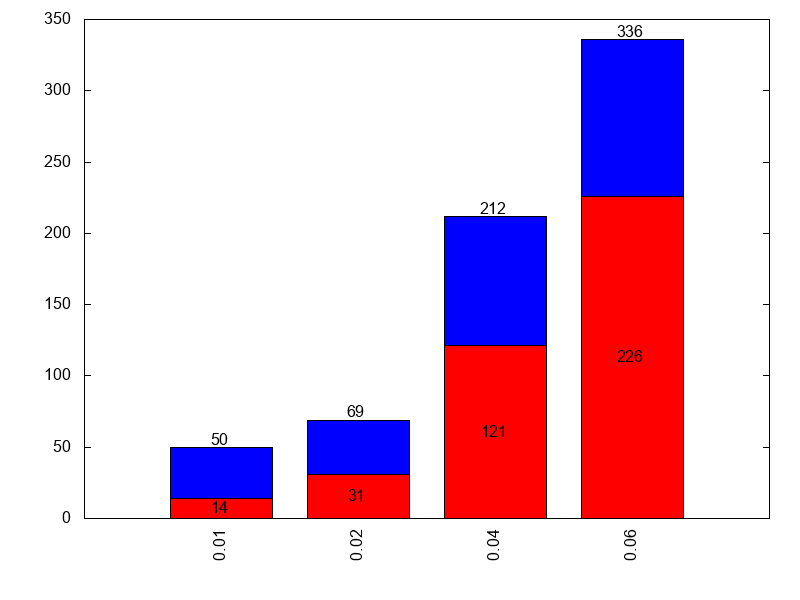}
        \caption{$\sigma_f$}
    \end{subfigure}
	\begin{subfigure}[b]{0.32\textwidth}
        \includegraphics[width=1.0\linewidth]{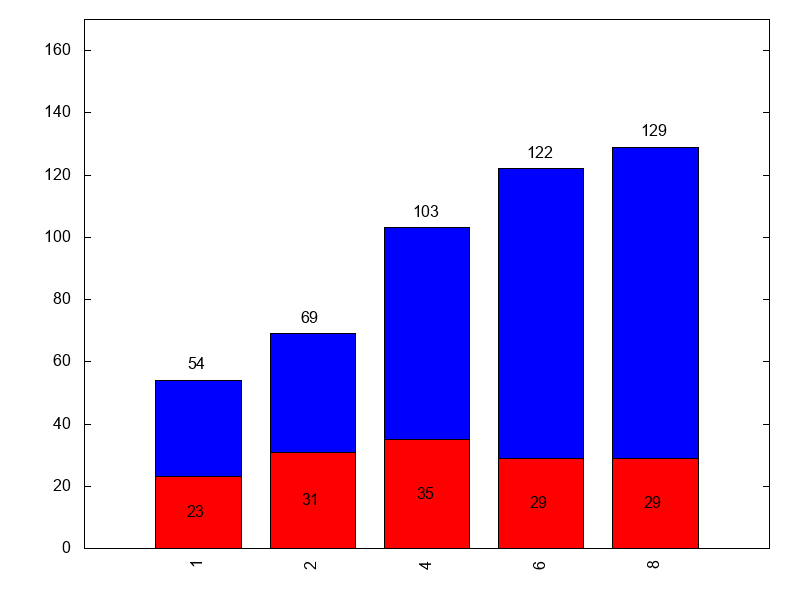}
        \caption{$\beta$}
    \end{subfigure}
	\begin{subfigure}[b]{0.32\textwidth}
        \includegraphics[width=1.0\linewidth]{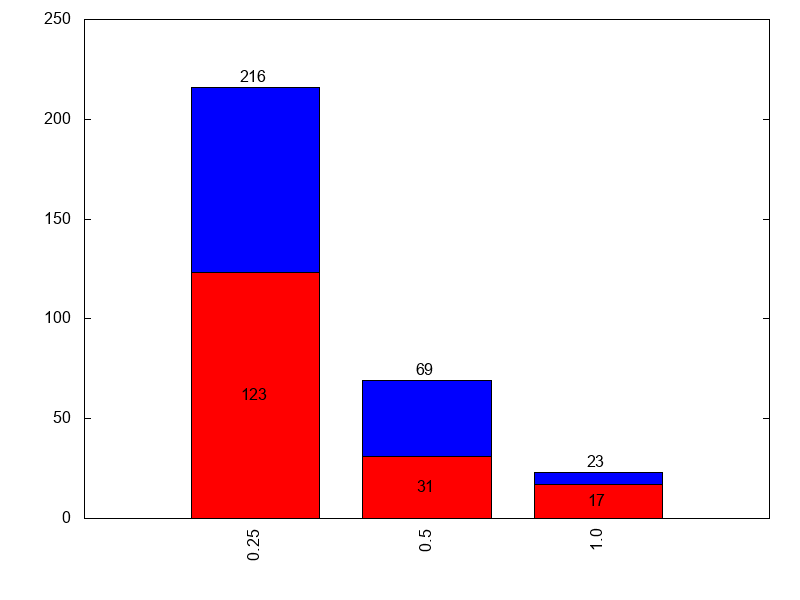}
        \caption{$\delta$}
    \end{subfigure}
\caption{Number of accepted  and number of passed first stage. 
Two - dimensional problem. 
Two-stage MCMC with $M_+ = 2$. 
(a) $\delta = 0.5$, $\sigma_f = 0.02$, $\beta =1, 2, 4, 6$ and $8$.  
(b) $\delta = 0.5$, $\beta = 2$, $\sigma_f = 0.01, 0.02, 0.04$ and $0.06$.  
(c) $\beta = 2$, $\sigma_f =0.02$, $\delta = 0,25, 0.5$ and $1.0$. 
}
\label{fig:acc-mcmc-hist}
\end{figure}

% 2d pic

% 3d pic

Next, we present results for MCMC algorithms with $\sigma_c^2 = \beta \cdot \, \sigma_f^2$. We test  1000 random fields proposals and use a random walk sampler for transitional probability distribution 
\[
q(\theta|\theta_n) = \theta_n + \delta \cdot r,
\] 
where $r$ is a Gaussian random variable with zero mean and variance 1. 
We use relative difference between observation data and proposed solution form \eqref{eq:E} with  $F =u(x, t)$, $F^* = u^{ms}(x, t)$ where $t = T_{max}$ (final time)  and $x \in  \partial \Omega_{surf}$ ($ \partial \Omega_{surf}$ is the top boundary of the domain).  

We consider the influence of the following parameters on the results and perform the following tests:
\begin{itemize}
\item Single-stage and two-stage MCMC method with $M_+ = 2$ and $4$ for 2D problem.
\item Two-stage MCMC for three-dimensional test problem with $M_+ = 2$.
\item Influence of the $\beta$, $\sigma_f$ and $\delta$ to the results of the two-stage MCMC method.
\end{itemize}

% 2d
In Figure \ref{fig:acc-mcmc-M}, we present results for single-stage and two-stage MCMC method with $M_+ = 2$ and $4$. In the first row, we depict an accepted $u_{obs}$ on the surface boundary, and on the second row, we present acceptance errors. Results are presented for two - dimensional problem. In MCMC algorithm, we use a random walk sampler with $\delta = 0.5$ and set $\sigma_f = 0.02$. We have 34 accepted fields in a single-stage method. In the two-stage MCMC with $M_+ = 2$ we have 31 accepted and 69 passed first stage. For $M_+ = 4$, we obtain 46 accepted and 68 passed first stage fields. 
% Fig9
In the first row of Figure \ref{fig:acc-mcmc-M}, we depicted an accepted $u_{obs}$ on the surface boundary, where with gradient coloring from yellow to red related to the number of accepted fields. We observe a convergence of the method from  second row, where we depicted relative differences between observation data and the proposed solution. Here $E^*(\theta)$ is shown in green color for multiscale solver (fist stage) and $E(\theta)$is shown in blue color for $\theta$ that passed the first stage and calculated using fine grid solver.  

% 2d
In Figure \ref{fig:passed-2d}, we present examples of accepted permeabilities with a reference solution that we used to calculate observation data (two-dimensional problem). 2D results are shown for two-stage MCMC algorithm with $\sigma_f = 0.02$, $\beta = 2$, $M_+ = 2$ and random walk sampler with $\delta = 0.5$. On the first and second rows, we depict permeability and elastic modulus, $k$ and $E$. Pressure and magnitude of the displacements are depicted in the third and fourth rows.
% 3d
In Figures \ref{fig:acc-mcmc-3d} and \ref{fig:passed-3d}, we present results for three-dimensional case. Results are presented for the two-stage MCMC algorithm with $\sigma_f = 0.02$, $\beta = 2$, $M_+ = 2$ and random walk sampler with $\delta = 0.5$. 
% Figs 10, 11 and 13, 14
For 2D and 3D results with the same parameters, we obtain:
\begin{itemize}
\item 2D: we have 31 accepted and 69 passed the first stage. 
\item 3D: we have 29 accepted and 141 passed the first stage. 
\end{itemize}

% 2d params
Numerical results for two-stage MCMC with different $\sigma_f = 0.01, 0.04, 0.06$ are presented in Figure \ref{fig:acc-mcmc-sigma}. We perform calculations for $\beta = 2$ and random walk sampler with $\delta = 0.5$.
%We obtain 14 accepted  (50 passed first stage) for $\sigma_f = 0.01$; 121 accepted (212 passed first stage) for $\sigma_f = 0.04$ and  226 accepted (336 passed first stage) for $\sigma_f = 0.06$.  
In Figure \ref{fig:acc-mcmc-hist}, we show a histogram with number of accepted and number of passed first stage for two - dimensional problem with $M_+ = 2$ and $\delta = 0.5$.
In the first picture in Figure \ref{fig:acc-mcmc-hist}, we consider $\delta = 0.5$, $\beta = 2$, $\sigma_f = 0.01, 0.02, 0.04$ and $0.06$. We have 
\begin{itemize}
\item $\sigma_f = 0.01$: 14 accepted and 50 passed first stage.
\item $\sigma_f = 0.02$: 31 accepted and 69 passed first stage.
\item $\sigma_f = 0.04$: 121 accepted and 212 passed first stage.
\item $\sigma_f = 0.06$: 226 accepted and 336 passed first stage.
\end{itemize}
In the second picture, we  show results for
\[
\sigma_c = \beta \cdot \sigma_f,
\]
with  $\sigma_f = 0.02$, $\beta =1, 2, 4, 6$ and $8$. We obtain following results
\begin{itemize}
\item $\beta =1$: 23 accepted and 54 passed first stage.
\item $\beta =2$: 31 accepted and 69 passed first stage.
\item $\beta =4$: 35 accepted and 103 passed first stage.
\item $\beta =6$: 29 accepted and 122 passed first stage.
\item $\beta =8$: 29 accepted and 129 passed first stage.
\end{itemize}
 In third picture in Figure \ref{fig:acc-mcmc-hist}, we present results for  $\beta = 2$, $\sigma_f =0.02$, $\delta = 0,25, 0.5$ and $1.0$. We obtain
 \begin{itemize}
\item $\delta = 0.25$: 123 accepted and 216 passed first stage.
\item $\delta = 0.5$: 31 accepted and 69 passed first stage.
\item $\delta = 1.0$: 17 accepted and 23 passed first stage.
\end{itemize}
We see that larger $\beta$ leads to larger acceptance on the first stage. We observe that $\beta = 2$ is better to obtain an optimal number of the accepted on the first stage vs the number of the accepted on the fine grid. 
By increasing the parameter $\sigma_f$, we increase the number of the accepted fields by reducing requirements between proposed and observation data differences. 
The smaller parameter $\delta$ from random walk leads to the larger number of accepted fields, but larger $\delta$ leads to the larger jump in random field generator.

\subsection{Preconditioned MCMC using Machine Learning}\label{ss3}

Finally, we consider the preconditioned MCMC method using a machine learning technique (\textit{ML}). We construct the neural networks for the prediction of the observable data (displacements on the surface boundary). Neural networks are constructed for each direction of displacements separately. 

% 2d NN
\begin{figure}[h!]
\centering
   \begin{subfigure}[b]{1.0\textwidth}
        \includegraphics[width=1.0\linewidth]{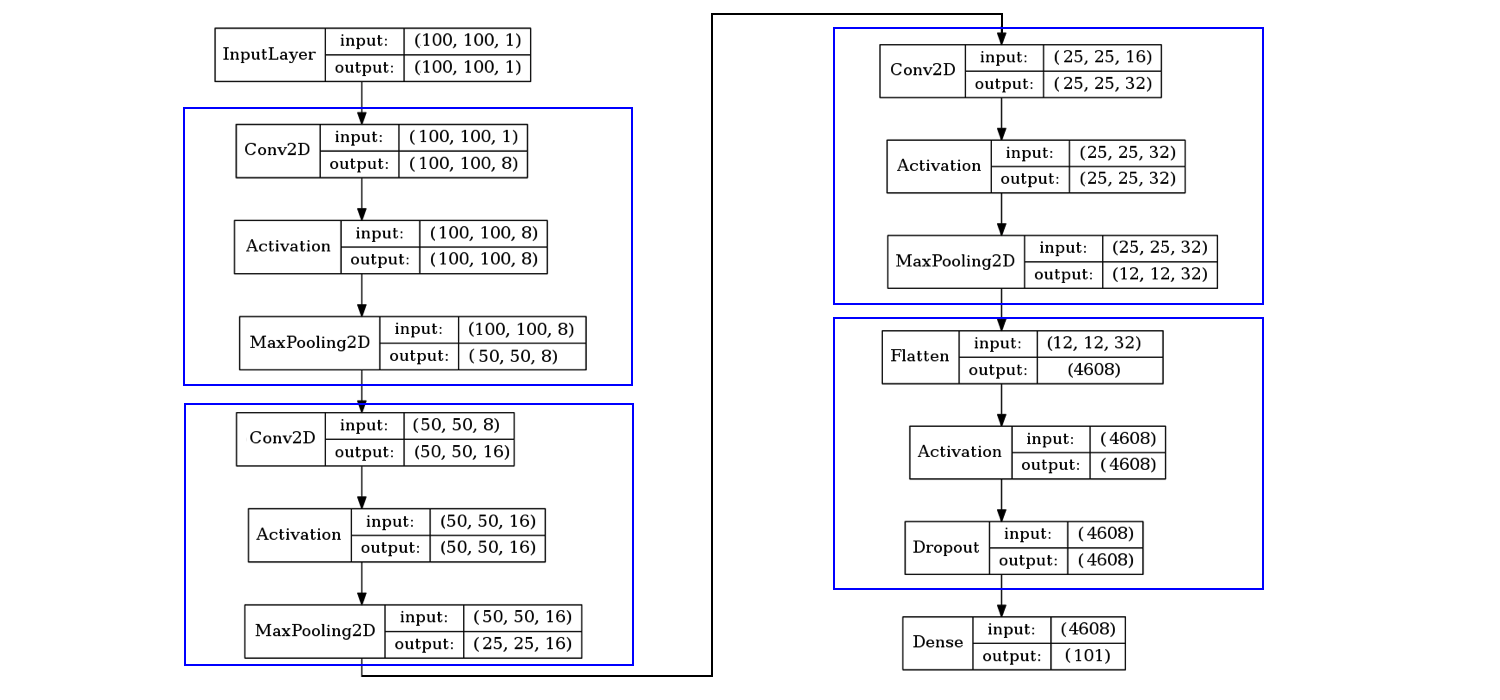}
        \caption{architecture of NN}
    \end{subfigure}
    \begin{subfigure}[b]{1.0\textwidth}
        \includegraphics[width=0.49\linewidth]{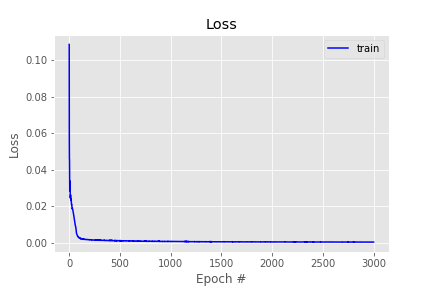}
        \includegraphics[width=0.49\linewidth]{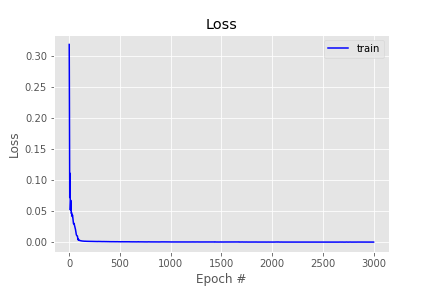}
        \caption{loss function vs epoch}
    \end{subfigure}
\caption{Two neural networks for the 2D problem.
(a) the architecture of neural networks.
(b) learning performance, loss function vs epoch.
$NN_x$ for $x$-displacement with RMSE = 1.887 \%, and  
$NN_y$ for $y$-displacement with RMSE = 1.386 \% (from left to right). }
\label{fig:nn-2d}
\end{figure} 

% 3d NN
\begin{figure}[h!]
\centering
   \begin{subfigure}[b]{1.0\textwidth}
        \includegraphics[width=1.0\linewidth]{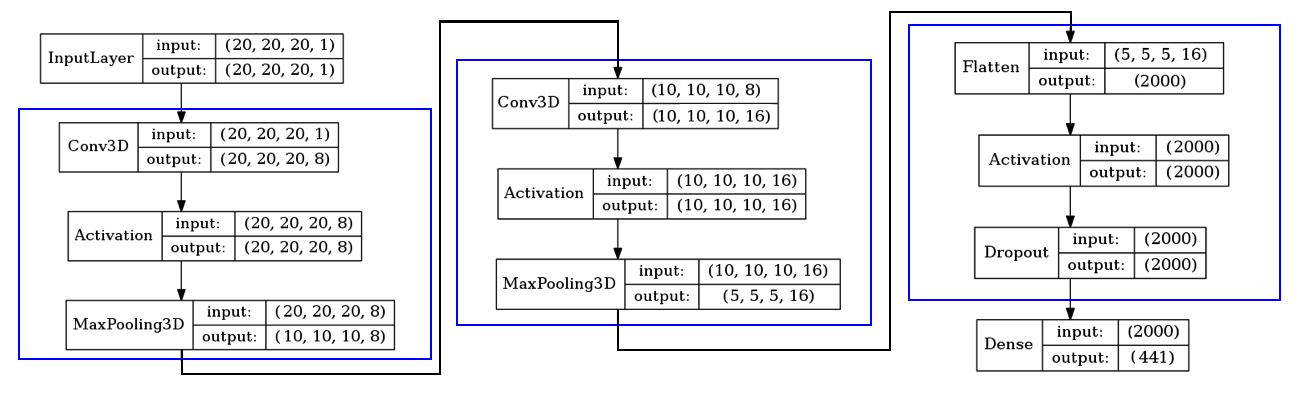}
        \caption{architecture of NN}
    \end{subfigure}
    \begin{subfigure}[b]{1.0\textwidth}
        \includegraphics[width=0.32\linewidth]{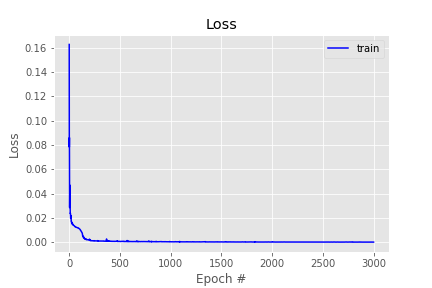}
        \includegraphics[width=0.32\linewidth]{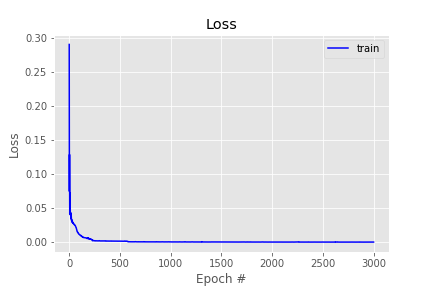}
        \includegraphics[width=0.32\linewidth]{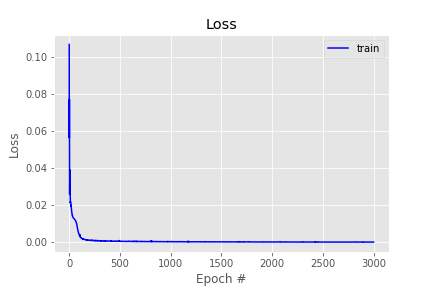}        
        \caption{loss function vs epoch}
    \end{subfigure}
\caption{Three neural networks for the 3D problem. 
(a) the architecture of neural networks. 
(b) learning performance, loss function vs epoch. 
$NN_x$ for $x$-displacement with RMSE = 1.104 \%, 
$NN_y$ for $y$-displacement with  RMSE = 1.746 \%, and  
$NN_z$ for $z$-displacement with  RMSE = 1.734 \% (from left to right). }
\label{fig:nn-3d}
\end{figure}

\begin{table}[h!]
\begin{center}
\begin{tabular}{ | c | c c c | }
\hline
& MSE & RMSE  (\%)  & MAE (\%) \\
\hline
\multicolumn{4}{|c|}{2D}  \\
\hline
$NN_x$ 	& 0.0007 & 1.8878 & 1.6719 \\
$NN_y$ 	& 0.0002 & 1.3864 & 1.2660 \\
\hline
\multicolumn{4}{|c|}{3D}  \\
\hline
$NN_x$ 	& 0.0001 & 1.1049 & 0.8978 \\
$NN_y$ 	& 0.0003 & 1.7468 & 1.5461 \\
$NN_z$ 	& 0.0003 & 1.7343 & 1.4264 \\
\hline
\end{tabular}
\end{center}
\caption{Learning performance of a machine learning algorithms for 2D and 3D formulations. }
\label{tab:ml}
\end{table}

% errors 100
\begin{figure}[htb!]
\centering
   \begin{subfigure}[b]{0.45\textwidth}
        \includegraphics[width=1.0\linewidth]{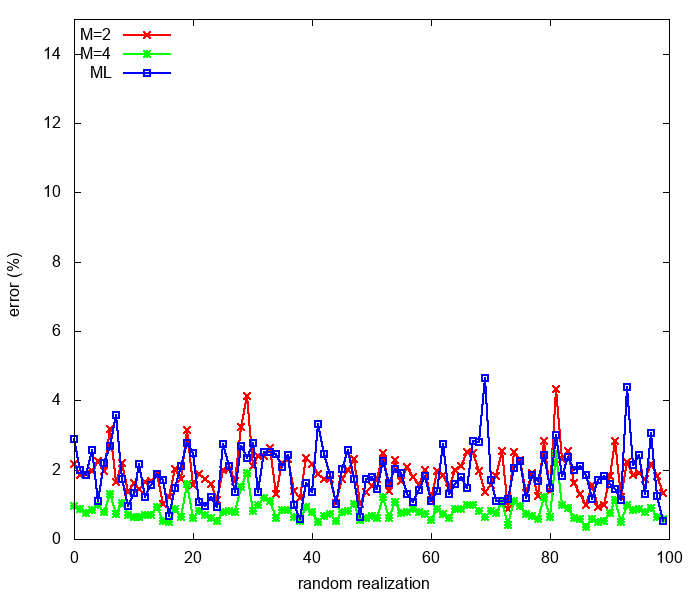}
        \caption{2D}
    \end{subfigure}
    \begin{subfigure}[b]{0.45\textwidth}
        \includegraphics[width=1.0\linewidth]{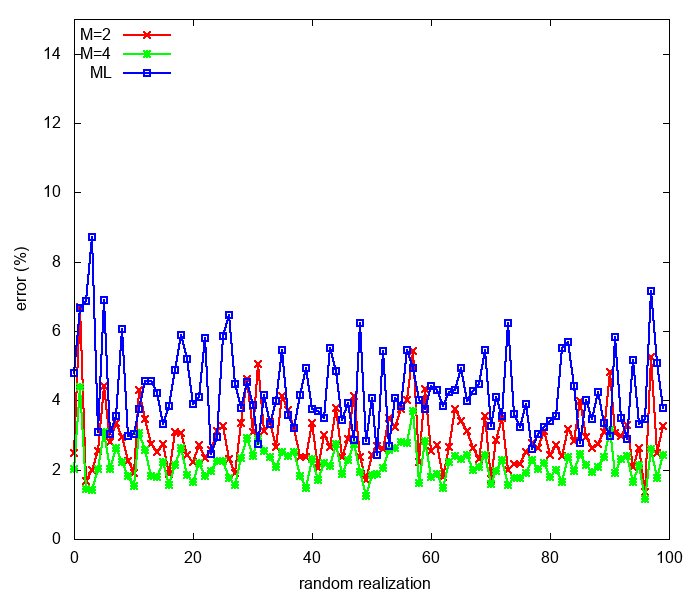}
        \caption{3D}
    \end{subfigure}
\caption{Multiscale method errors for 100 random realizations of heterogeneous filed. Prediction of the machine learning algorithm and multiscale solver with  $M_+=2, 4$.  Relative errors in \% for displacements on the top boundary.  }
\label{fig:err-ml-100}
\end{figure} 

% cp 100
\begin{figure}[h!]
\centering
    \begin{subfigure}[b]{0.32\textwidth}
        \includegraphics[width=1.0\linewidth]{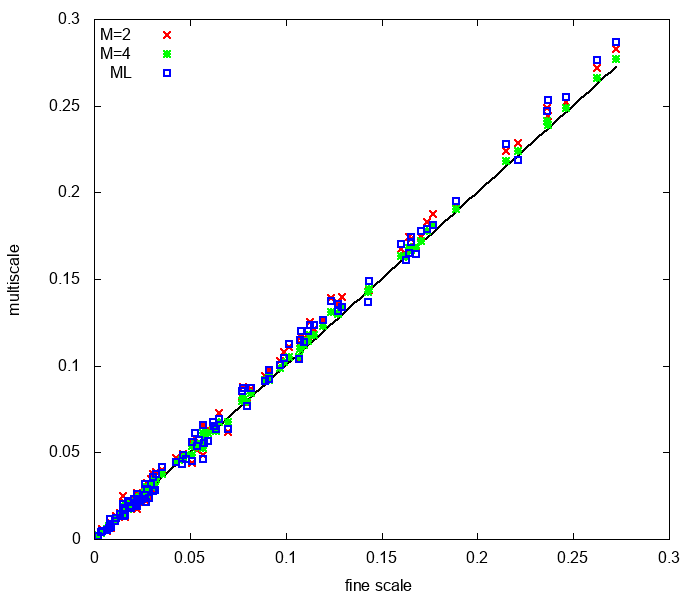}
        \caption{Case 1 for 2D}
    \end{subfigure}
    \begin{subfigure}[b]{0.32\textwidth}
        \includegraphics[width=1.0\linewidth]{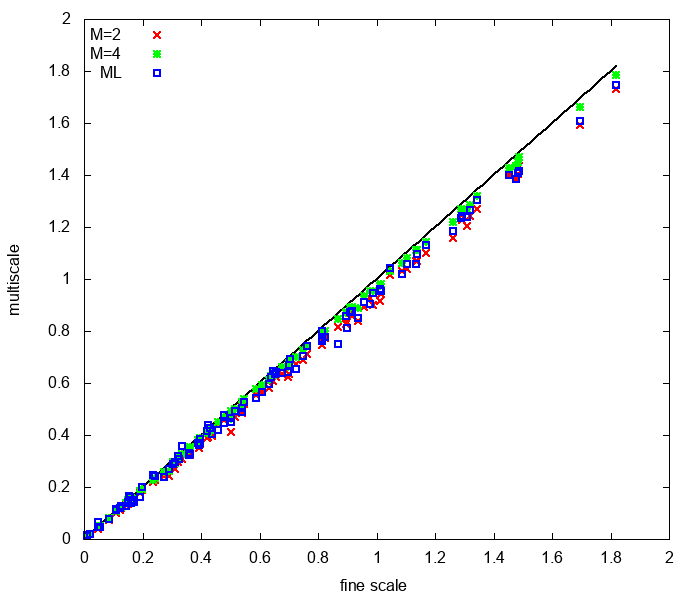}
        \caption{Case 2 for 2D}
    \end{subfigure}
    \begin{subfigure}[b]{0.32\textwidth}
        \includegraphics[width=1.0\linewidth]{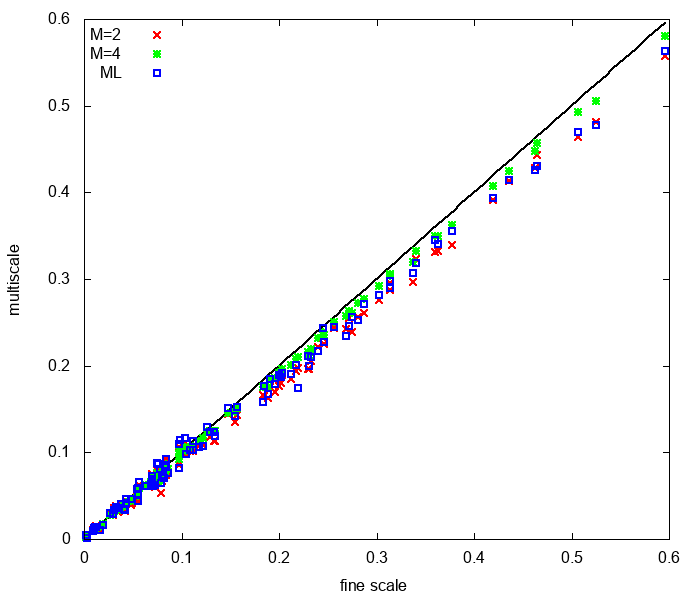}
        \caption{Case 3 for 2D}
    \end{subfigure} \\
    \vspace{10pt}
    \begin{subfigure}[b]{0.32\textwidth}
        \includegraphics[width=1.0\linewidth]{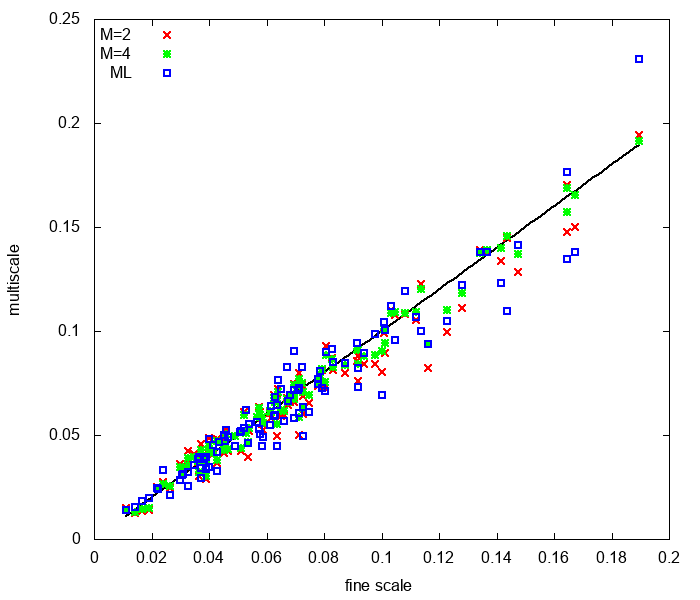}
        \caption{Case 1 for 3D}
    \end{subfigure}
    \begin{subfigure}[b]{0.32\textwidth}
        \includegraphics[width=1.0\linewidth]{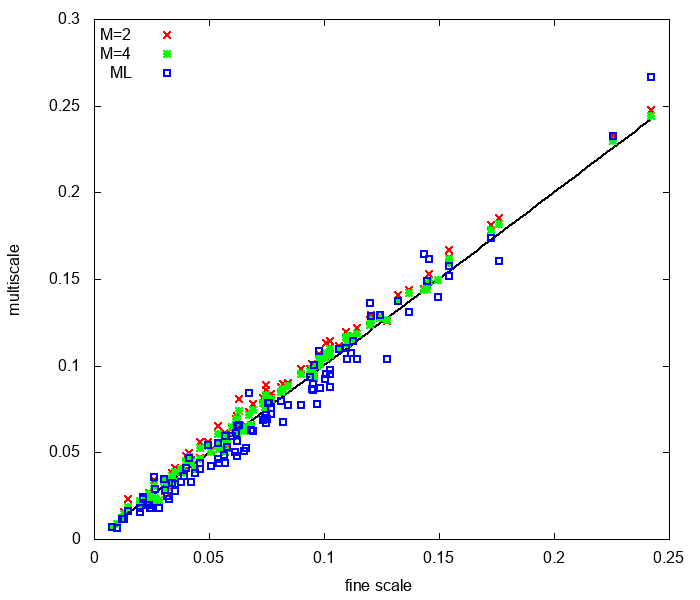}
        \caption{Case 2 for 3D}
    \end{subfigure}
    \begin{subfigure}[b]{0.32\textwidth}
        \includegraphics[width=1.0\linewidth]{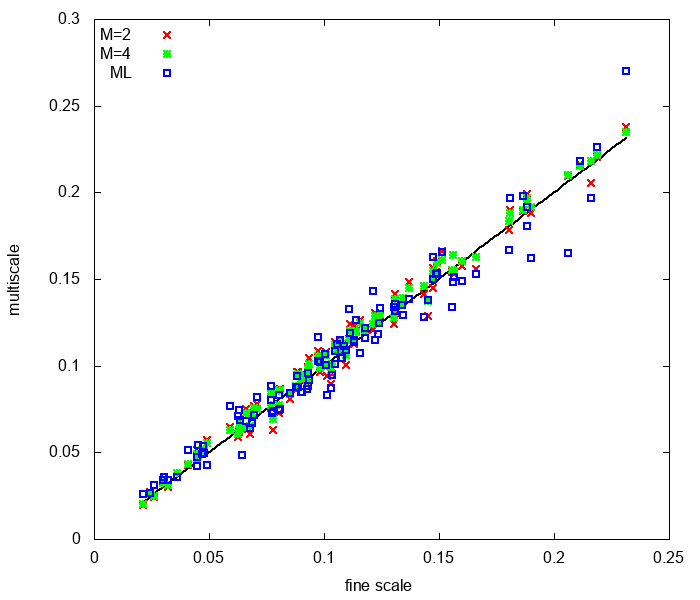}
        \caption{Case 3 for 3D}
    \end{subfigure}
\caption{Cross-plot  between ${E(\theta)}$ (fine scale) and ${E^*(\theta)}$. 
Prediction of the machine learning algorithm and multiscale solver with  $M_+=2, 4$. Case 1, 2 and 3 observation data (from left to right). 
First row: two-dimensional problem. 
Second row: three-dimensional problem. }
\label{fig:ml-cp}
\end{figure} 

% ms 2d
\begin{figure}[h!]
\centering
	\begin{subfigure}[b]{0.32\textwidth}
        \includegraphics[width=1.0\linewidth]{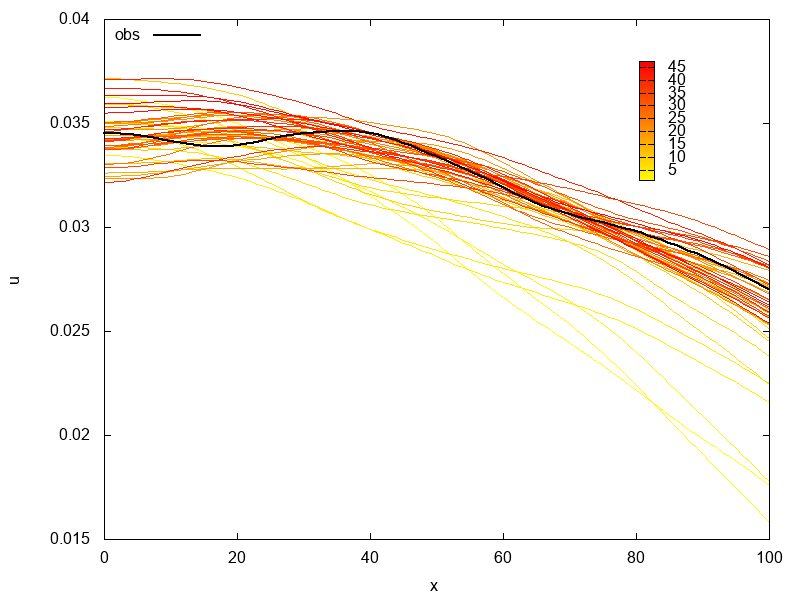}\\
 		\includegraphics[width=1.0\linewidth]{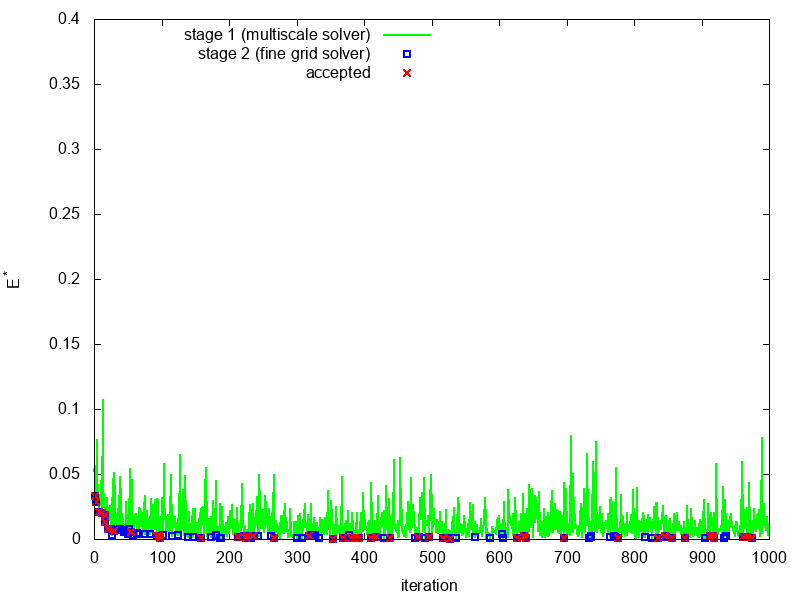}
        \caption{Case 1 for 2D  (\textit{MS})} 
    \end{subfigure}
    \begin{subfigure}[b]{0.32\textwidth}
        \includegraphics[width=1.0\linewidth]{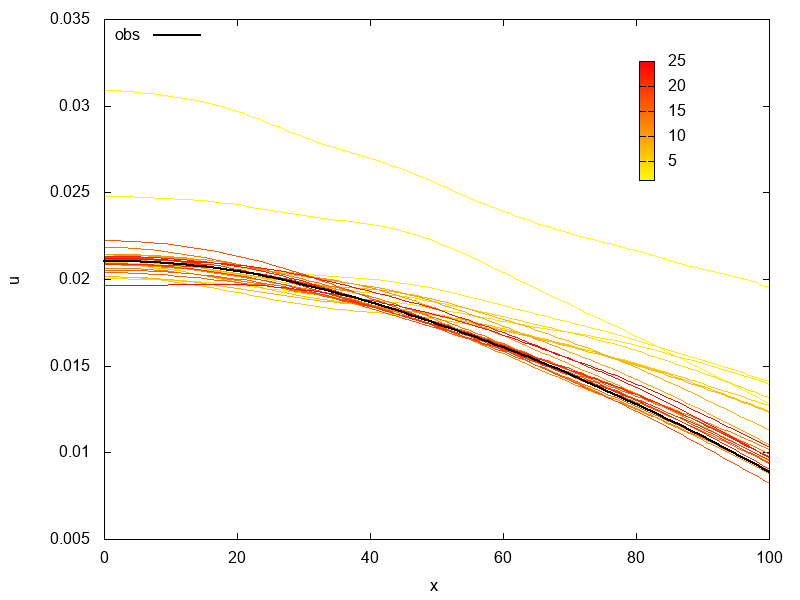}\\
 		\includegraphics[width=1.0\linewidth]{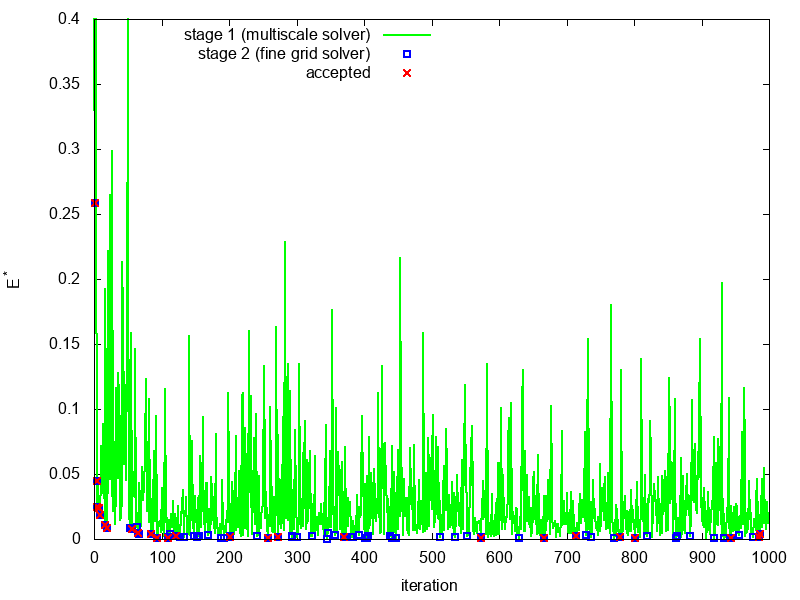}
        \caption{Case 3 for 2D  (\textit{MS})} 
    \end{subfigure}
    \begin{subfigure}[b]{0.32\textwidth}
        \includegraphics[width=1.0\linewidth]{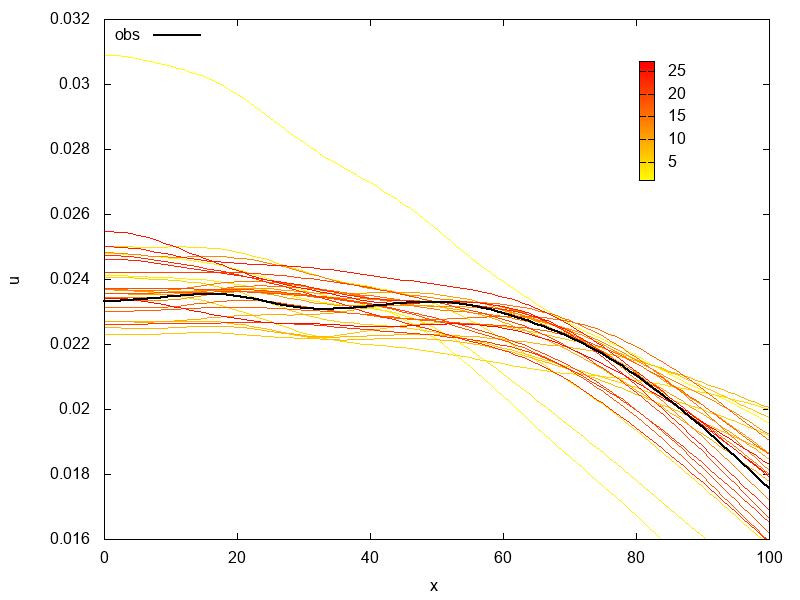}\\
 		\includegraphics[width=1.0\linewidth]{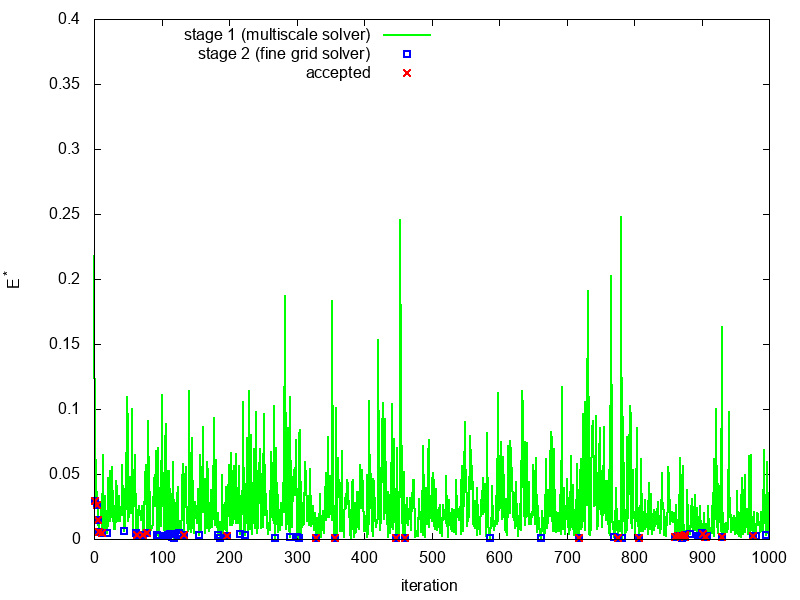}
        \caption{Case 2 for 2D  (\textit{MS})} 
    \end{subfigure}
\caption{Two - dimensional problem with \textit{MS} preconditioning. 
Two-stage MCMC with $\sigma_f = 0.02$ and  $\beta = 2$. 
First row: accepted $u_{obs}$ on surface boundary. 
Second row:  $E^*(\theta)$ (green color) and $E(\theta)$(blue color) in each MCMC iteration. 
(a) Case 1. (b) Case 2. (c) Case 3. }
\label{fig:ms-sol-2d}
\end{figure} 

% ml 2d
\begin{figure}[h!]
\centering
 \centering
 	\begin{subfigure}[b]{0.32\textwidth}
        \includegraphics[width=1.0\linewidth]{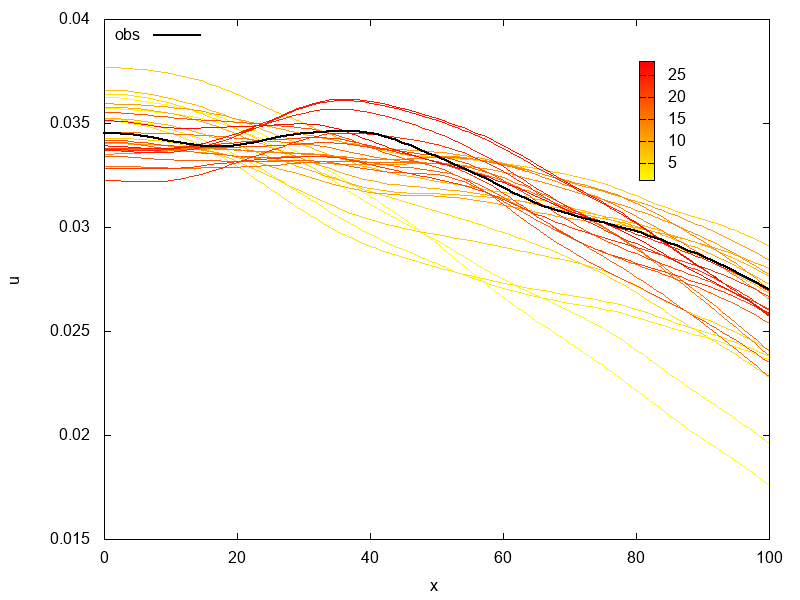}\\
 		\includegraphics[width=1.0\linewidth]{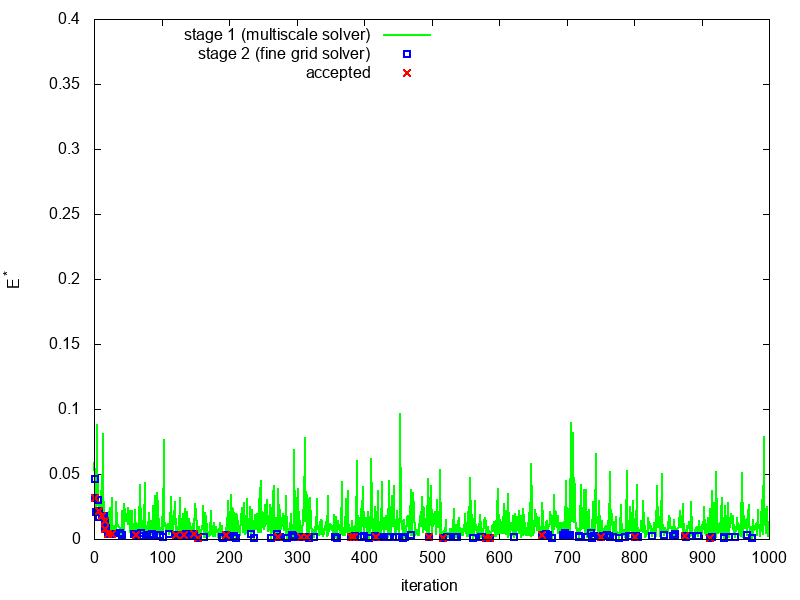}
        \caption{Case 1 for 2D  (\textit{ML})} 
    \end{subfigure}
    \begin{subfigure}[b]{0.32\textwidth}
        \includegraphics[width=1.0\linewidth]{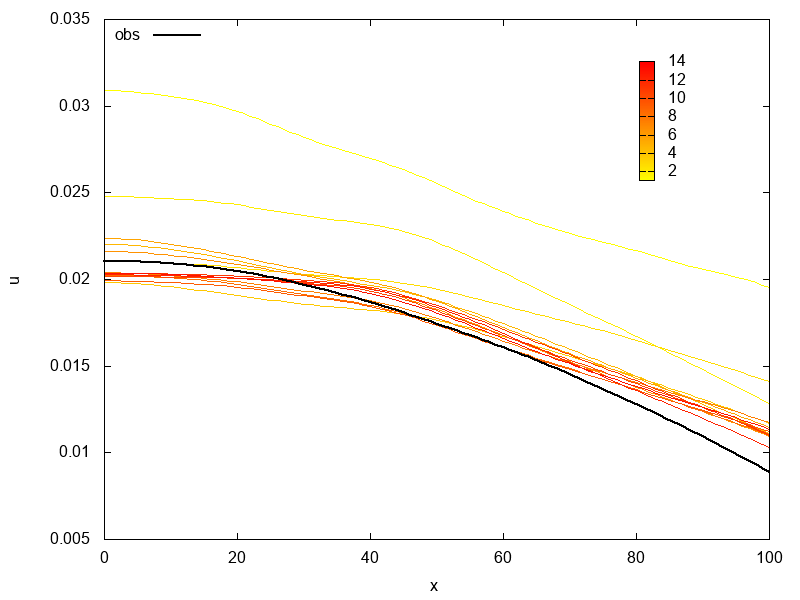}\\
 		\includegraphics[width=1.0\linewidth]{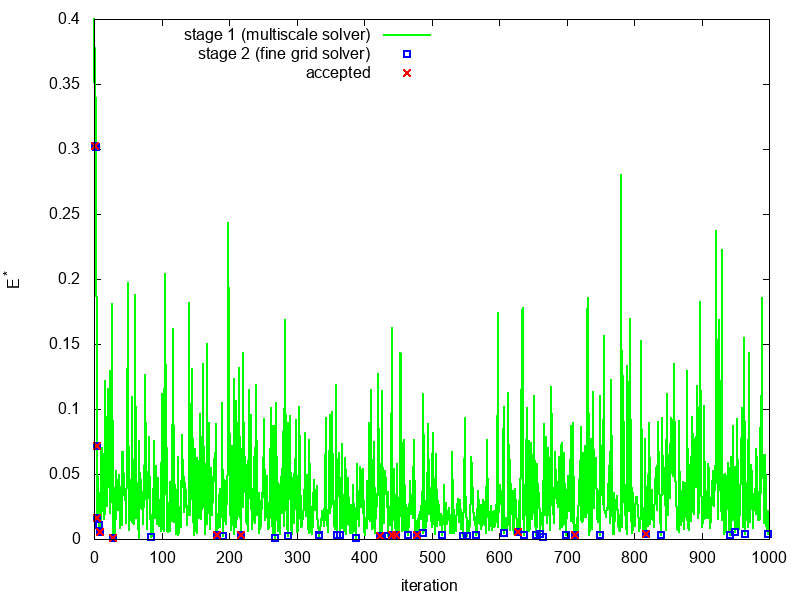}
        \caption{Case 3 for 2D  (\textit{ML})} 
    \end{subfigure}
    \begin{subfigure}[b]{0.32\textwidth}
        \includegraphics[width=1.0\linewidth]{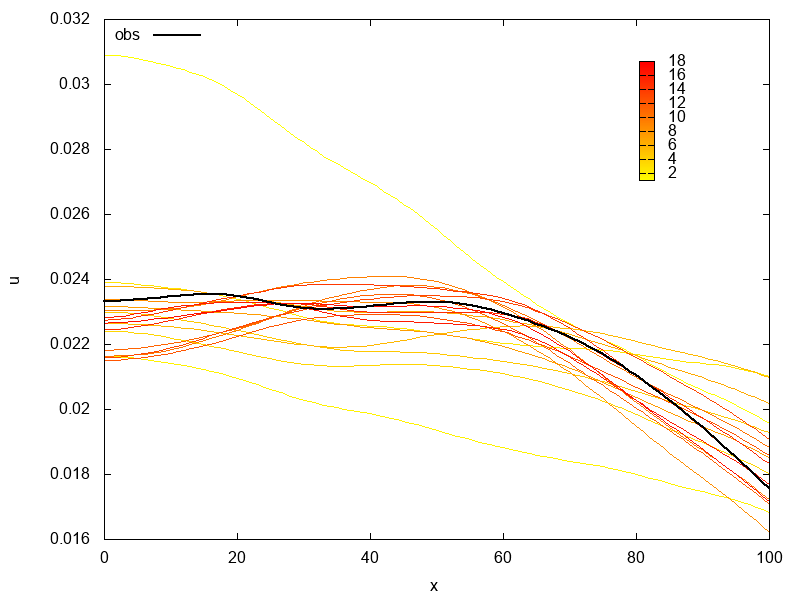}\\
 		\includegraphics[width=1.0\linewidth]{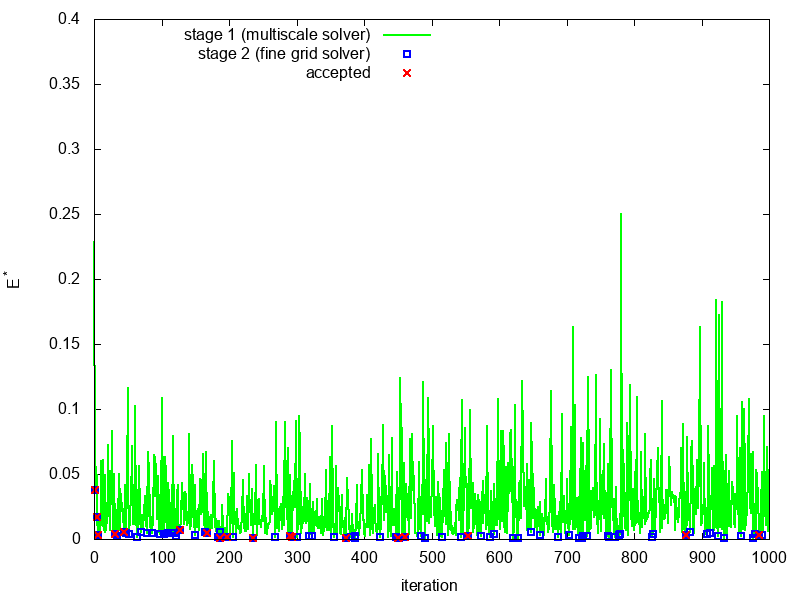}
        \caption{Case 2 for 2D  (\textit{ML})} 
    \end{subfigure}
\caption{Two - dimensional problem with \textit{ML} preconditioning. 
Two-stage MCMC with $\sigma_f = 0.02$ and  $\beta = 2$. 
First row: accepted $u_{obs}$ on surface boundary. 
Second row:  $E^*(\theta)$ (green color) and $E(\theta)$(blue color) in each MCMC iteration. 
(a) Case 1. (b) Case 2. (c) Case 3.}
\label{fig:ml-sol-2d}
\end{figure}

% 2d pic
\begin{figure}[h!]
\centering
\begin{subfigure}[b]{1.0\textwidth}
\includegraphics[width=1.0\textwidth]{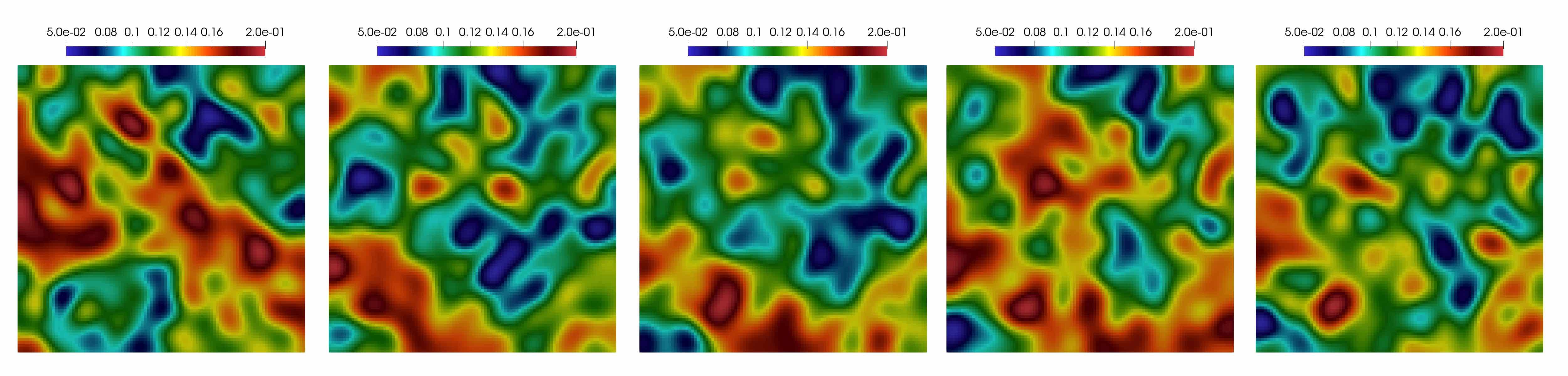}
\caption{Case 1 for 2D  (\textit{ML})} 
\end{subfigure}\\
\vspace{10pt}
\begin{subfigure}[b]{1.0\textwidth}
\includegraphics[width=1.0\textwidth]{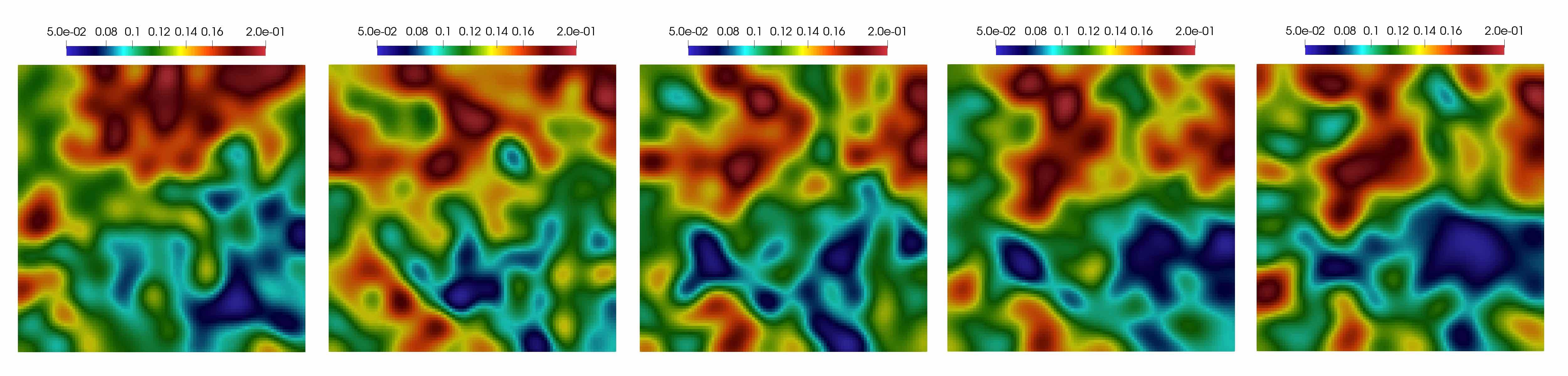}
\caption{Case 2 for 2D  (\textit{ML})} 
\end{subfigure}\\
\vspace{10pt}
\begin{subfigure}[b]{1.0\textwidth}
\includegraphics[width=1.0\textwidth]{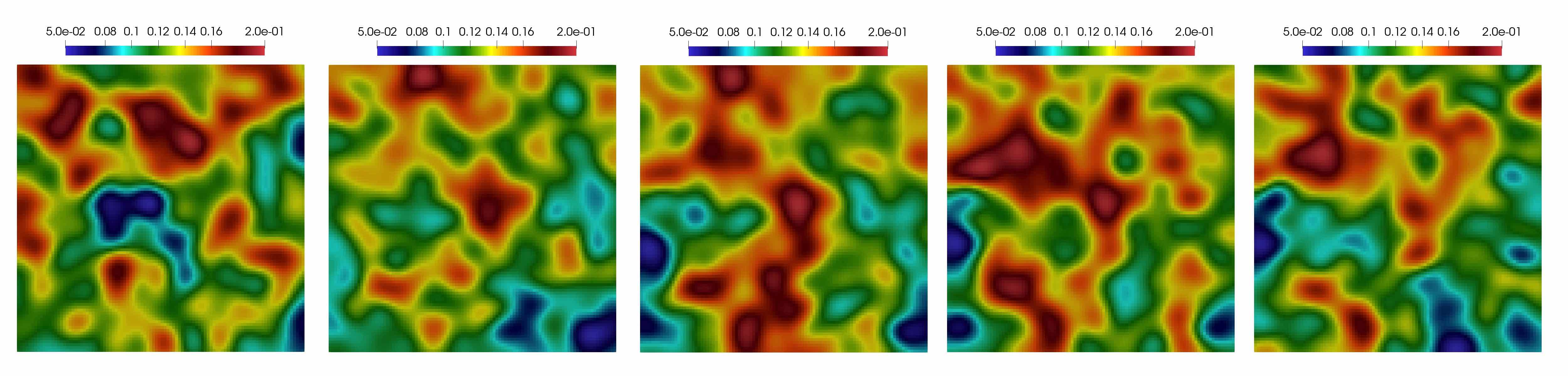}
\caption{Case 3 for 2D  (\textit{ML})} 
\end{subfigure}
\caption{Reference field and accepted random fields for three - dimensional problem (from left to right).  
Two - dimensional problem with \textit{ML} preconditioning. 
Two-stage MCMC with $\sigma_f = 0.02$ and  $\beta = 2$. 
(a) $\phi$ for Case 1. (b) $\phi$ for Case 2. (c) $\phi$ for Case 3.
 }
\label{fig:passed-2d-ml}
\end{figure}

% ms 3d
\begin{figure}[h!]
 \centering
 	\begin{subfigure}[b]{0.32\textwidth}
        \includegraphics[width=1.0\linewidth]{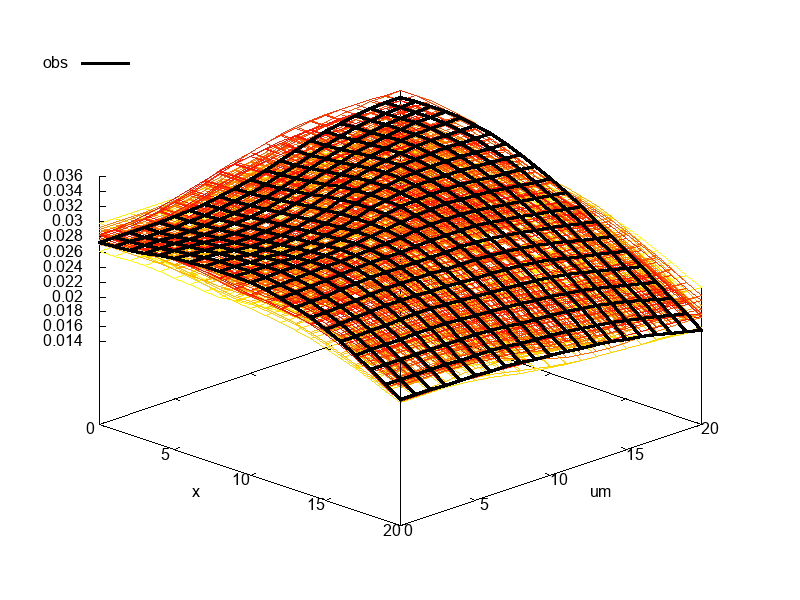}\\
 		\includegraphics[width=1.0\linewidth]{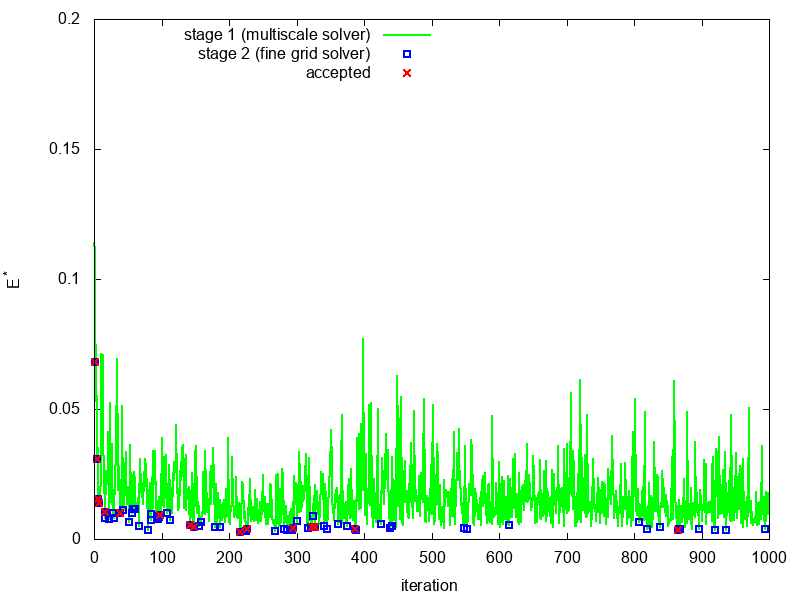}
        \caption{Case 1 for 3D (\textit{MS})} 
    \end{subfigure}
    \begin{subfigure}[b]{0.32\textwidth}
        \includegraphics[width=1.0\linewidth]{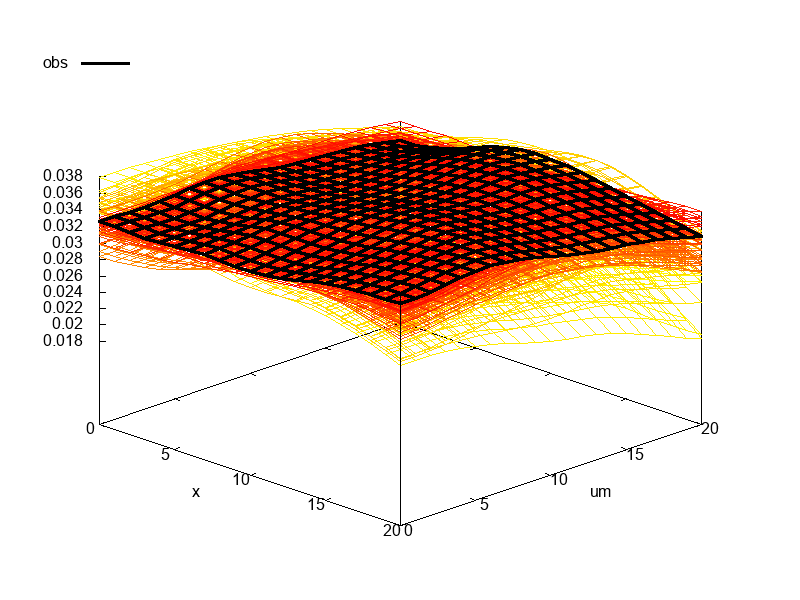}\\
 		\includegraphics[width=1.0\linewidth]{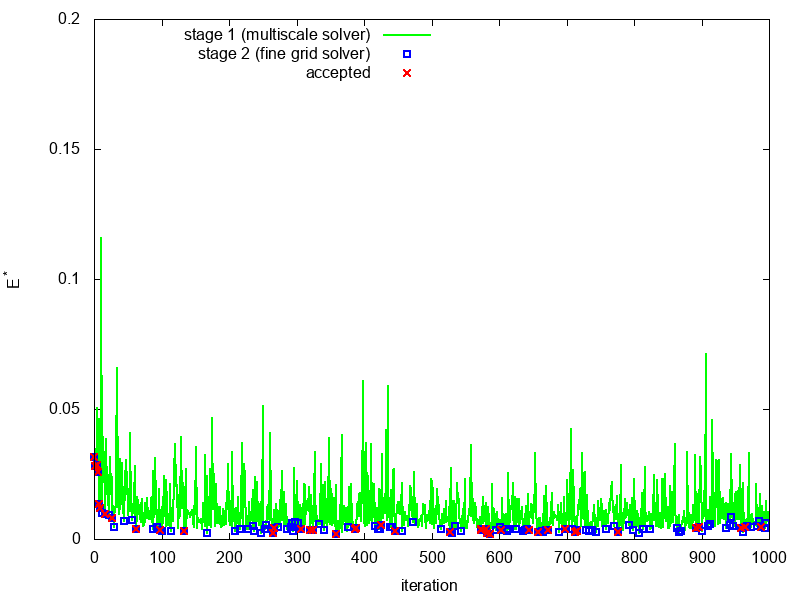}
        \caption{Case 3 for 3D  (\textit{MS})} 
    \end{subfigure}
    \begin{subfigure}[b]{0.32\textwidth}
        \includegraphics[width=1.0\linewidth]{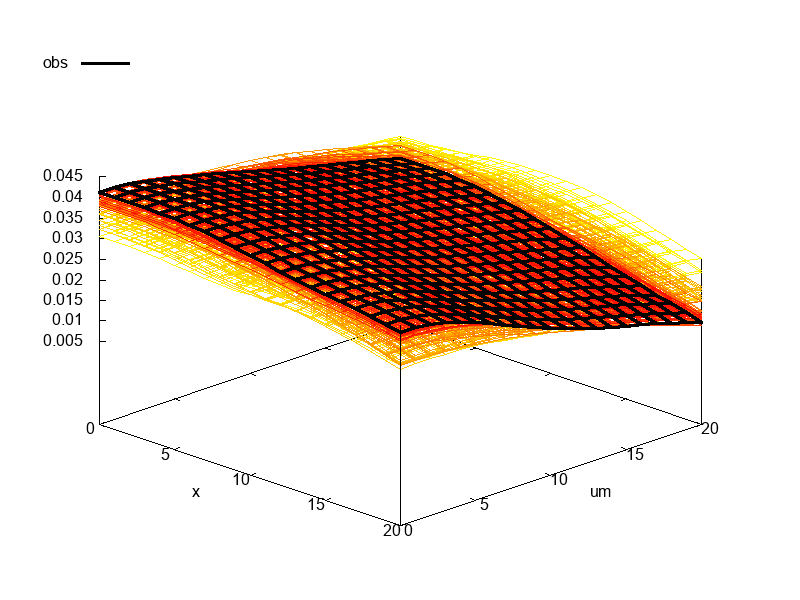}\\
 		\includegraphics[width=1.0\linewidth]{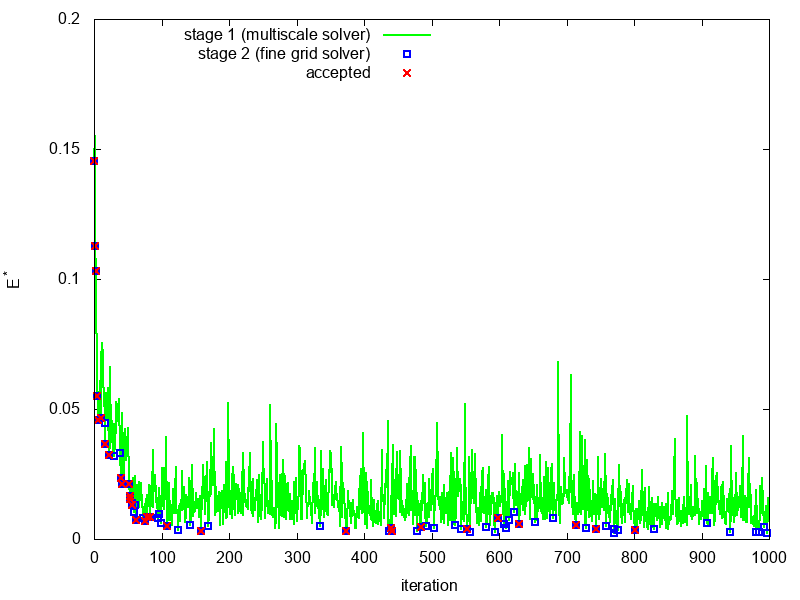}
        \caption{Case 2 for 3D  (\textit{MS})} 
    \end{subfigure}
\caption{Three - dimensional problem with \textit{MS} preconditioning. 
Two-stage MCMC with $\sigma_f = 0.02$ and  $\beta = 2$. 
First row: accepted $u_{obs}$ on surface boundary. 
Second row:  $E^*(\theta)$ (green color) and $E(\theta)$(blue color) in each MCMC iteration. 
(a) Case 1. (b) Case 2. (c) Case 3.}
\label{fig:ms-sol-3d}
\end{figure} 

% ml 3d
\begin{figure}[h!]
 \centering
 	\begin{subfigure}[b]{0.32\textwidth}
        \includegraphics[width=1.0\linewidth]{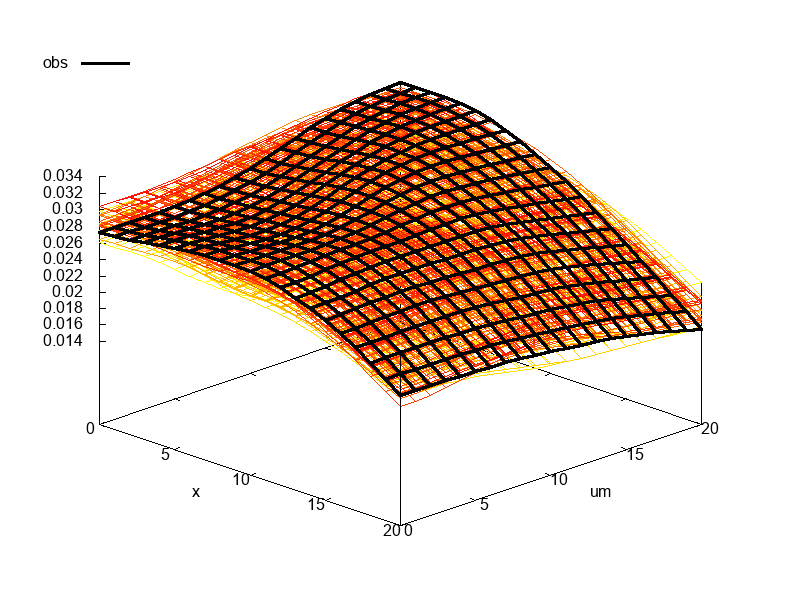}\\
        \includegraphics[width=1.0\linewidth]{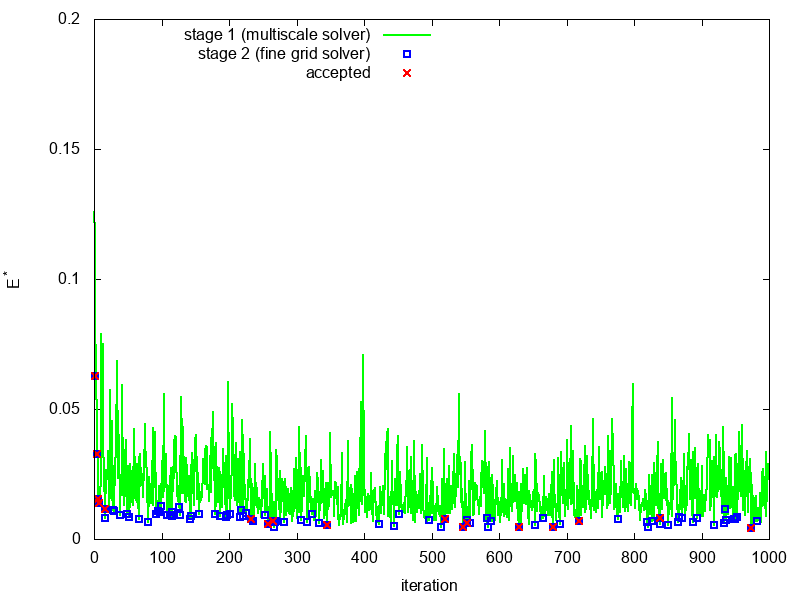}
        \caption{Case 1 for 3D  (\textit{ML})} 
    \end{subfigure}
    \begin{subfigure}[b]{0.32\textwidth}
        \includegraphics[width=1.0\linewidth]{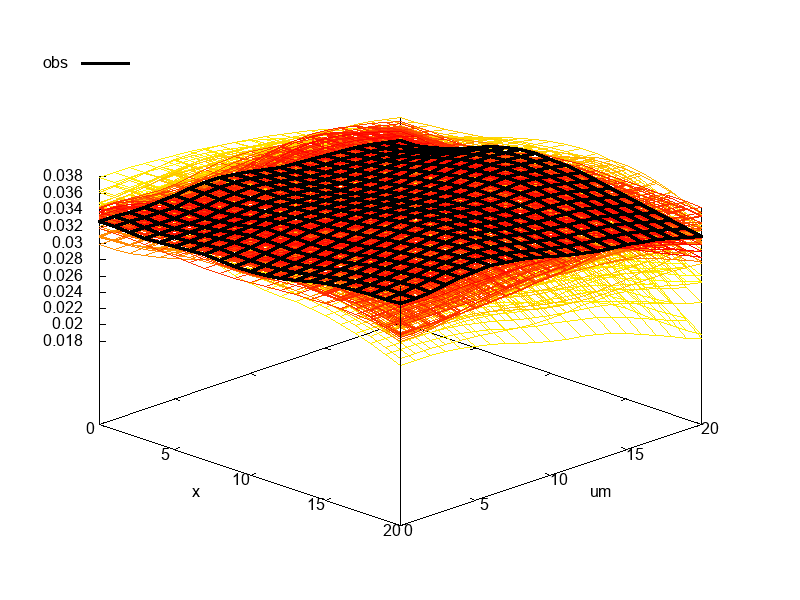}\\
        \includegraphics[width=1.0\linewidth]{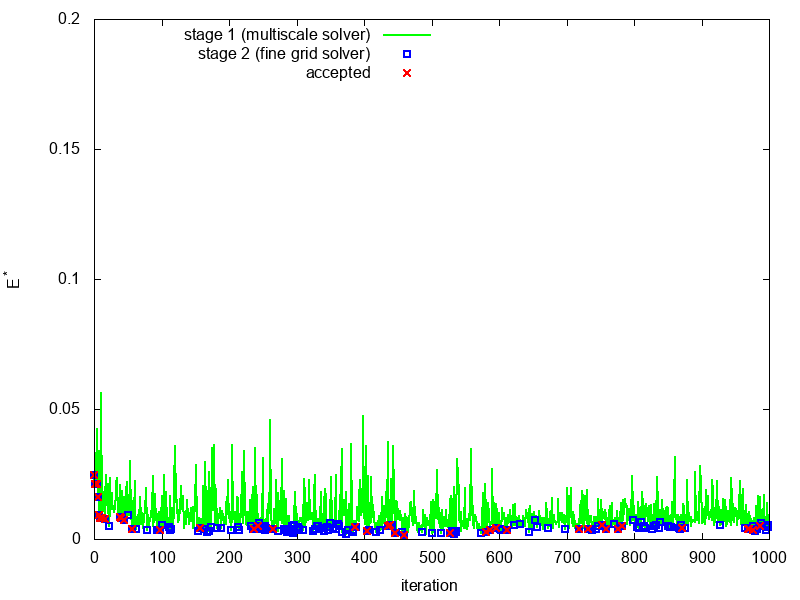}
        \caption{Case 3 for 3D  (\textit{ML})} 
    \end{subfigure}
    \begin{subfigure}[b]{0.32\textwidth}
        \includegraphics[width=1.0\linewidth]{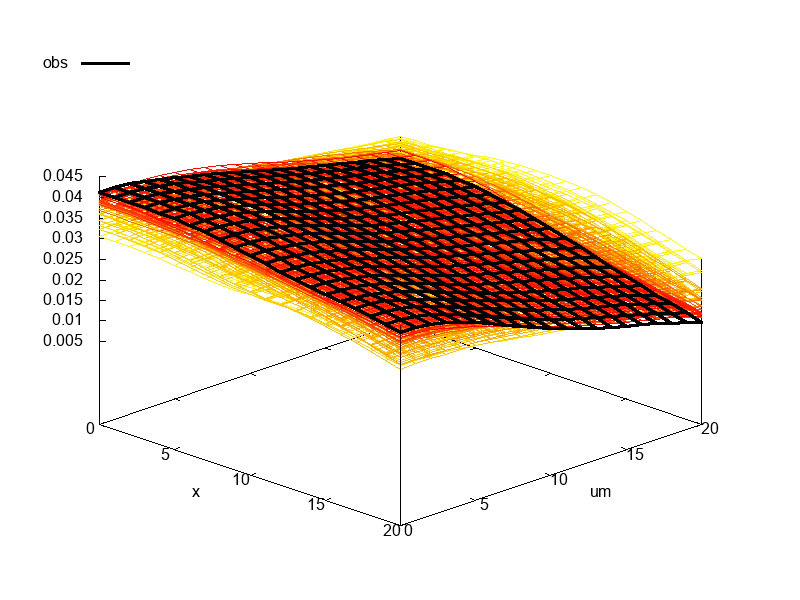}\\
        \includegraphics[width=1.0\linewidth]{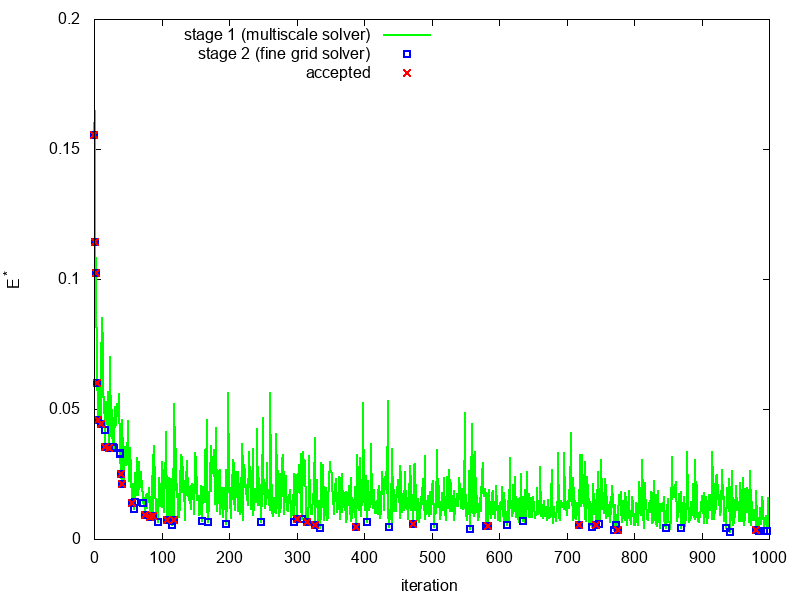}
        \caption{Case 2 for 3D  (\textit{ML})} 
    \end{subfigure}
\caption{Three - dimensional problem with \textit{ML} preconditioning. 
Two-stage MCMC with $\sigma_f = 0.02$ and  $\beta = 2$. 
First row: accepted $u_{obs}$ on surface boundary. 
Second row:  $E^*(\theta)$ (green color) and $E(\theta)$(blue color) in each MCMC iteration. 
(a) Case 1. (b) Case 2. (c) Case 3.}
\label{fig:ml-sol-3d}
\end{figure}

% 3d pic
\begin{figure}[h!]
\centering
\begin{subfigure}[b]{1.0\textwidth}
\includegraphics[width=1.0\textwidth]{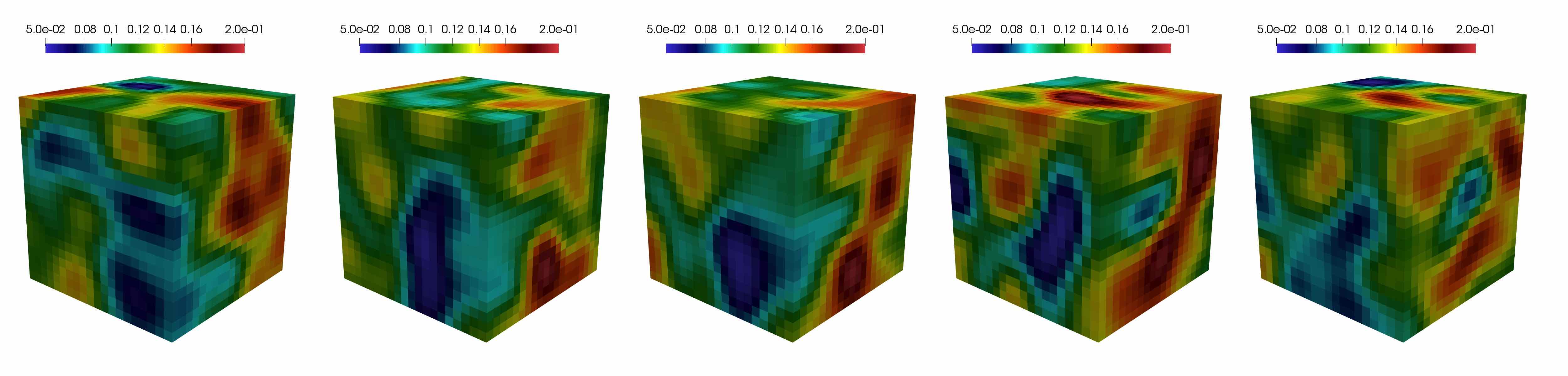}
\caption{Case 1 for 3D  (\textit{ML})} 
\end{subfigure}\\
\vspace{10pt}
\begin{subfigure}[b]{1.0\textwidth}
\includegraphics[width=1.0\textwidth]{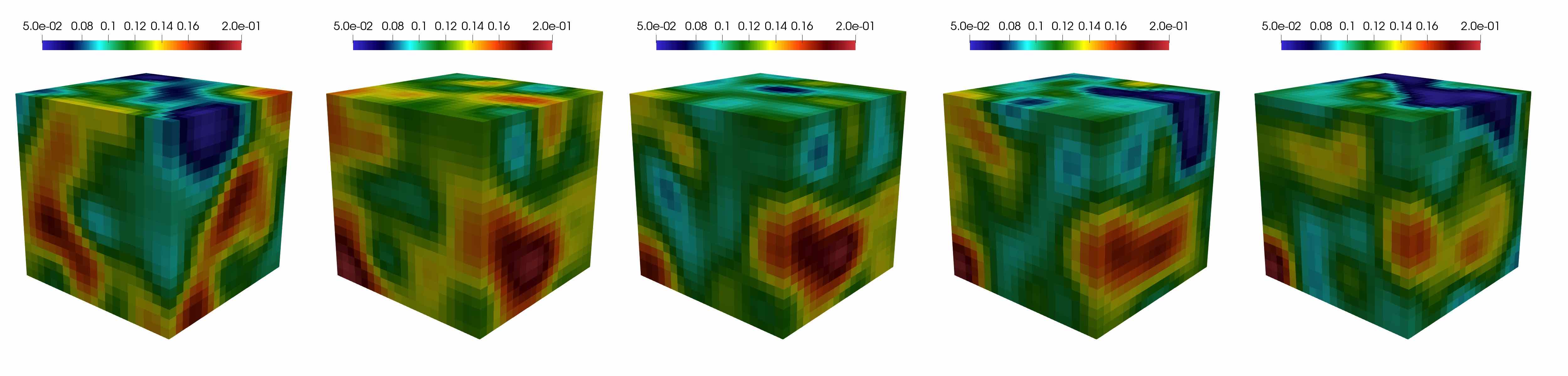}
\caption{Case 2 for 3D  (\textit{ML})} 
\end{subfigure}\\
\vspace{10pt}
\begin{subfigure}[b]{1.0\textwidth}
\includegraphics[width=1.0\textwidth]{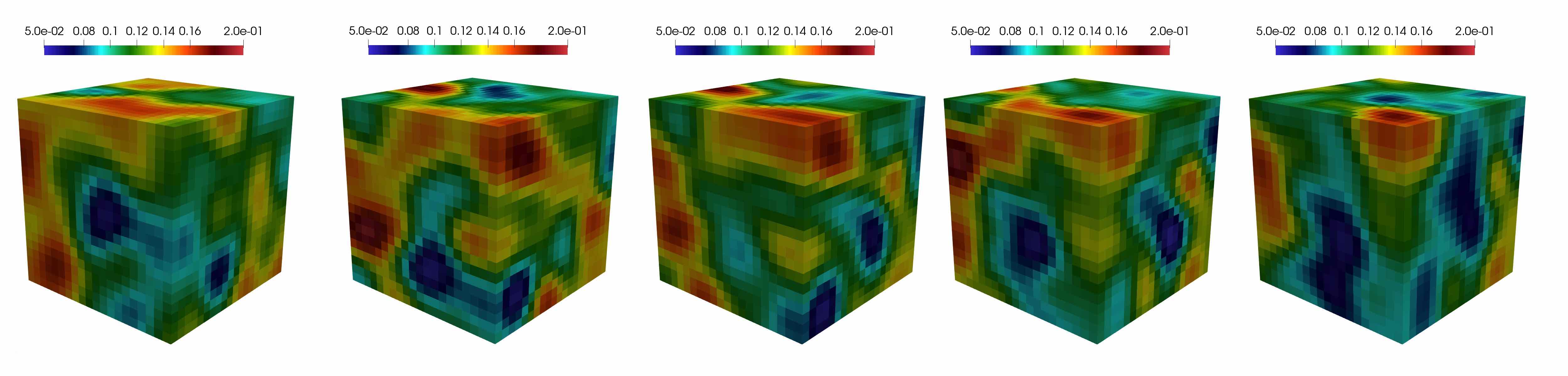}
\caption{Case 3 for 3D  (\textit{ML})} 
\end{subfigure}
\caption{Reference field and accepted random fields for three - dimensional problem (from left to right).  
Three - dimensional problem with \textit{ML} preconditioning. 
Two-stage MCMC with $\sigma_f = 0.02$ and  $\beta = 2$. 
(a) $\phi$ for Case 1. (b) $\phi$ for Case 2. (c) $\phi$ for Case 3.
 }
\label{fig:passed-3d-ml}
\end{figure}

The architectures of the neural networks are presented in Figures \ref{fig:nn-2d} and \ref{fig:nn-3d} for 2D and 3D problems, respectively. In training process, we perform 3000 epochs using Adam optimizer with learning rate $\epsilon = 0.001$. 
For accelerating the training process of the CNN, we use GPU (GTX 1800 Ti). 
We use $3^d$ convolutions and $2^d$ maxpooling layers with RELU activation for $d = 2,3$.  We have several layers of convolutions with several final fully connected layers.  We use dropout with rate 10 \%  in order to prevent over-fitting.
As a loss function, we use a mean square error (MSE). 
Convergence of the loss function presented in Figures \ref{fig:nn-2d} and \ref{fig:nn-3d}, where we plot the MSE loss function vs epoch number.

For error calculation, we use mean square errors, relative mean absolute and relative root mean square errors
\[
MSE = \sum_i |Q_i - \tilde{Q}_i|^2,
\quad
RMSE = \sqrt{ \frac{\sum_i |Q_i - \tilde{Q}_i|^2 }{\sum_i |Q_i|^2 } },
\quad
MAE = \frac{\sum_i |Q_i - \tilde{Q}_i|  }{\sum_i |Q_i|},
\]
where $Q_i$  and $\tilde{Q}_i$ denotes reference and predicted values for sample $X_i$. 
Learning performance for neural networks are presented in Tables \ref{tab:ml} for two and three - dimensional problems. We observe good convergence of the relative errors with $\approx 1 \%$ of RMSE.

In Figure \ref{fig:ml-cp}, we present parity plots comparing reference values against predicted using trained neural networks for Case 1, 2 and 3. Cross plot demonstrate correlation between ${E^*(\theta)}$ and ${E(\theta)}$ for 100 realizations of random fields for three cases in 2D and 3D formulations.  
In Figure \ref{fig:err-ml-100}, we present a relative errors for displacements on the top boundary.  We compare prediction using reference fine grid values of the displacements.  We observe sufficiently good errors for predicted values, where we have less than $5 \%$ of errors for 2D and less than $10 \%$ of errors for 3D. 
Later, we will discuss the computational efficiency of the machine learning algorithms compared with direct forward calculations and preconditioning using multiscale solver. 

Next, we consider numerical results for the two-stage  MCMC method with machine learning technique (\textit{ML}) and multiscale solver (\textit{MS}). 
In Figures \ref{fig:ms-sol-2d} and \ref{fig:ml-sol-2d}, we present results two-dimensional problem for Case 1, 2 and 3 (see Section \ref{ss1}). In the first row, we depict an accepted $u_{obs}$ on surface boundary  and on the second row, we present an acceptence errors.  In MCMC algorithm, we use a random walk sampler with $\delta = 0.5$. We use  $\sigma_f = 0.02$ and set $\sigma_c = 2 \cdot  \sigma_f$. 

For \textit{MS} preconditioning, we have
\begin{itemize}
\item Case 1 with 47 accepted and 105 passed first stage.
\item Case 2 with 27 accepted and 62 passed first stage.
\item Case 3 with 25 accepted and 72 passed first stage.
\end{itemize}

For \textit{ML} preconditioning, we have
\begin{itemize}
\item Case 1 with 28 accepted and 132 passed first stage.
\item Case 2 with 18 accepted and 81 passed first stage.
\item Case 3 with 14 accepted and 45 passed first stage.
\end{itemize}

The results for three-dimensional problem are presented in Figures \ref{fig:ms-sol-3d} and \ref{fig:ml-sol-3d} for Case 1, 2 and 3 (see Section \ref{ss1}). 

For \textit{MS} preconditioning, we have
\begin{itemize}
\item Case 1 with 16 accepted and 67 passed first stage.
\item Case 2 with 30 accepted and 71 passed first stage.
\item Case 3 with 40 accepted and 124 passed first stage.
\end{itemize}

For \textit{ML} preconditioning, we have
\begin{itemize}
\item Case 1 with 17 accepted and 89 passed first stage.
\item Case 2 with 26 accepted and 63 passed first stage.
\item Case 3 with 39 accepted and 159 passed first stage.
\end{itemize}

In Figures \ref{fig:passed-2d-ml} and \ref{fig:passed-3d-ml}, we depict an  examples of accepted porosities with a reference porosity that we used to calculate observation data. Numerical results are shown for two-stage MCMC algorithm with \textit{ML} preconditioning for Case 1,2 and 3 in 2D and 3D formulations.

We observe that we can obtain a good acceptance rate with very cheap machine learning-based prediction. The acceptance rate is $N_{accepted}/N_{fine}$, $N_{accepted}$ and $N_{fine}$ are the number of accepted fields and number of expensive fine grid calculations (passed the first stage).  For the fast construction of the dataset that used for training, we used a multiscale solver with $M_+ = 2$. 

Finally, we discuss the advantage of the proposed algorithm. In the single-stage MCMC method with $N_{iter}$ iterations, time of calculations $T_{F}$ is equal to number of iterations multiply to time of solution of the fine grid system
\[
T_{F}  = N_{iter} \cdot t_{fine},
\]
where $t_{fine}$ is the time of fine grid system solution. 
Here for 2D system with $DOF_f = 30603$, we have $t_{fine} = 6.2$ seconds and  $t_{fine} = 158.7$ seconds for 3D problem  with $DOF_f = 37044$ (see Table \ref{tab:err-time}). 

% ms
For \textit{MS} preconditioning of the MCMC method, we have
\[
T_{MS} = N_{iter} \cdot t_{ms} + N_{fine} \cdot t_{fine},
\]
where $t_{ms}$ is the time of coarse grid system solution using GMsFEM and $ N_{fine}$ is the number of accepted on the first stage. 
For $M_+ = 2$, we have $t_{ms} = 0.83$ seconds for 2D problem with $DOF_c = 847$ and $t_{ms} = 9.6$ seconds for 3D problem with $DOF_c = 1728$ (see Table \ref{tab:err-time}).
In preconditioned MCMC method, we obtain that $N_{fine} << N_{iter}$, and therefore, we have a huge reduction in the solution time because the presented multiscale method provides a huge reduction of the system size, but it still takes some time for the solution.

% ml
For further reduction of the time, we proposed a machine learning-based technique with a super quick prediction.
For \textit{ML} preconditioning of the MCMC method, we have
\[
T_{ML} = N_{iter} \cdot t_{ml} + N_{fine} \cdot t_{fine},
\]
where $t_{ml}$ is the time of coarse grid system solution using trained neural networks. 
Becuase prediction time is very fast  i.e. $t_{ml} << 1$ second, therefore
\[
T_{ML} = N_{fine} \cdot t_{fine}.
\]

\section{Conclusion} 

Simulation of the poroelasticity is difficult due to the complex heterogeneities and uncertainty. 
In this work, we considered a Two-stage Markov Chain Monte Carlo method for geomechanical subsidence. We presented two techniques for preconditioning: (MS) multiscale method for model order reduction and (ML) machine learning technique.  
Numerical results are presented for two- and three-dimensional models to show the efficiency of the method as an expedited MCMC sampling method. 

Codes used in this manuscript are publicly available on Bitbucket at https://bitbucket.org/vmasha/ms-mcmc.

\bibliographystyle{plain}
\bibliography{lit}

\end{document}